\documentclass{amsart}
\usepackage{amsmath} 
\usepackage{amssymb}

\newtheorem{theorem}{Theorem}[subsection] 
\newtheorem{claim}{Claim}[theorem]
\newtheorem{lemma}[theorem]{Lemma} 
\newtheorem{proposition}[theorem]{Proposition} 
\newtheorem{corollary}[theorem]{Corollary} 

\theoremstyle{definition}
\newtheorem{definition}[theorem]{Definition}
\newtheorem{example}[theorem]{Example}
\newtheorem{problem}{Problem}[section]

\theoremstyle{remark}
\newtheorem{remark}[theorem]{Remark}

\newtheorem{conclusion}[theorem]{Conclusion}
\newtheorem{context}[theorem]{Context}

\numberwithin{equation}{section}
\newcommand{\forces}{\Vdash} 

\newcommand{\incomp}{\bot}
 
\newcommand{\lbv}{[\![} 
\newcommand{\rbv}{]\!]}
\newcommand{\bV}{{\bf V}} 
\newcommand{\lesdot}{\mathrel{\mathord{<}\!\!\raise 
0.8 pt\hbox{$\scriptstyle\circ$}}} 
\newcommand{\comp}{\circ} 
  
\newcommand{\con}{{\mathfrak c}}

\newcommand{\can}{2^{\textstyle \omega}} 
\newcommand{\fs}{2^{\textstyle <\!\omega}} 
\newcommand{\baire}{\omega^{\textstyle \omega}} 
\newcommand{\iso}{[\omega]^{\textstyle \omega}} 
 
\newcommand{\fseo}{\omega^{\textstyle <\!\omega}} 

\newcommand{\conc}{{}^\frown\!}
\newcommand{\lh}{{\rm lh}\/}
\newcommand{\rest}{{\restriction}}
\newcommand{\mrot}{{\rm root}\/} 
\newcommand{\suc}{{\rm succ}} 

\newcommand{\dom}{{\rm dom}} 
\newcommand{\rng}{{\rm rng}}

\newcommand{\nor}{{\rm {\bf nor}}\/} 
\newcommand{\pos}{{\rm pos}}
\newcommand{\val}{{\bf val}}
\newcommand{\dis}{{\bf dis}}
\newcommand{\dn}{{\rm dn}}
\newcommand{\up}{{\rm up}}
\newcommand{\CR}{{\rm CR}}

\newcommand{\FC}{{\rm FC}}
\newcommand{\PFC}{{\rm PFC}}
\newcommand{\PC}{{\rm PC}}

\newcommand{\TCR}{{\rm LTCR}}

\newcommand{\tree}{{\rm tree}}
\newcommand{\mtree}{{\rm mt}}

\newcommand{\POS}{{\rm POS}}
\newcommand{\cnor}{{{\mathcal C}(\nor)}}
\newcommand{\hn}{{\bf hn}}
\newcommand{\HN}{{\bf HN}}
\newcommand{\cHa}{{\mathcal H}(\omega_1)}
\newcommand{\Ipw}{{\mathcal I}_{\bP,\dot{W}}}
\newcommand{\Ipwz}{{\mathcal I}^0_{\bP,\dot{W}}}
\newcommand{\borel}{{\bf Borel}}
\newcommand{\lev}{{\rm lev}}
\newcommand{\cA}{{\mathcal A}}
\newcommand{\bB}{{\bf B}}
\newcommand{\cB}{{\mathcal B}}
\newcommand{\BB}{{\mathbb B}}
\newcommand{\cC}{{\mathcal C}}
\newcommand{\bC}{{\mathbb C}}
\newcommand{\bD}{{\bf D}}
\newcommand{\cD}{{\mathcal D}}
\newcommand{\bH}{{\bf H}}
\newcommand{\bF}{{\bf F}}
\newcommand{\cF}{{\mathcal F}}
\newcommand{\cG}{{\mathcal G}}
\newcommand{\cI}{{\mathcal I}}
\newcommand{\cK}{{\mathcal K}}

\newcommand{\bP}{{\mathbb P}}
\newcommand{\cP}{{\mathcal P}}
\newcommand{\gp}{{\mathfrak p}}
\newcommand{\bQ}{{\mathbb Q}}
\newcommand{\dbQ}{{\dot{\mathbb Q}}}
\newcommand{\mbR}{{\mathbb R}}
\newcommand{\bS}{{\bf S}}
\newcommand{\cS}{{\mathcal S}}
\newcommand{\cT}{{\mathcal T}}
\newcommand{\cU}{{\mathcal U}}
\newcommand{\cX}{{\mathcal X}}
\newcommand{\cY}{{\mathcal Y}}
\newcommand{\cZ}{{\mathcal Z}}

\newcount\skewfactor
\def\mathunderaccent#1#2 {\let\theaccent#1\skewfactor#2
\mathpalette\putaccentunder}
\def\putaccentunder#1#2{\oalign{$#1#2$\crcr\hidewidth
\vbox to.2ex{\hbox{$#1\skew\skewfactor\theaccent{}$}\vss}\hidewidth}}

\newcommand{\qcnor}{{\bQ^*_\cnor(K,\Sigma,\Sigma^\bot)}}
\newcommand{\qzero}{{\bQ^*_\emptyset(K,\Sigma,\Sigma^\bot)}}
\newcommand{\qinf}{{\bQ^*_\infty(K,\Sigma,\Sigma^\bot)}}
\newcommand{\qcF}{{\bQ^*_\cF(K,\Sigma,\Sigma^\bot)}}
\newcommand{\qscF}{{\bQ^*_\cF(K,\Sigma)}}
\newcommand{\qonef}{{\bQ^*_f(K,\Sigma,\Sigma^\bot)}}
\newcommand{\qonefx}{{\bQ^*_f(K,\Sigma)}}
\newcommand{\bQp}{{\bQ^\tree(\gp)}}
\begin{document}

\setcounter{page}{0}
\title[Sweet {\&} Sour]{Sweet {\&} Sour and other flavours of ccc
forcing notions}

\author{Andrzej Ros{\l}anowski}
\address{Department of Mathematics\\
 University of Nebraska at Omaha\\
 Omaha, NE 68182-0243, USA} 
\email{roslanow@member.ams.org}
\urladdr{http://www.unomaha.edu/$\sim$aroslano}
\thanks{The first author thanks the Hebrew University of Jerusalem for
  their hospitality during his visits to Jerusalem in Summer'00 and
  Summer'01. His research was also partially supported by a grant from
  the University Committee on Research of UNO}
\author{Saharon Shelah}
\address{Institute of Mathematics\\
 The Hebrew University of Jerusalem\\
 91904 Jerusalem, Israel\\
 and  Department of Mathematics\\
 Rutgers University\\
 New Brunswick, NJ 08854, USA}
\email{shelah@math.huji.ac.il}
\urladdr{http://www.math.rutgers.edu/$\sim$shelah}
\thanks{The research of the second author was partially supported by the
 Israel Science Foundation. Publication 672} 

\date{November 2001}

\begin{abstract}
We continue developing the general theory of forcing notions built with the
use of {\em norms on possibilities}, this time concentrating on ccc forcing
notions and classifying them.
\end{abstract}

\maketitle

\tableofcontents

\section{Introduction}
The present paper has three themes. First, we continue the investigations
started in Judah, Ros{\l}anowski and Shelah \cite{JRSh:373} and
Ros{\l}anowski and Shelah \cite{RoSh:470}, \cite{RoSh:628}, and we
investigate the method of {\em norms on possibilities} in the context of ccc
forcing notions, getting a number of constructions of nicely definable ccc
forcings. Most of them fall into the class of nep--forcing notions of Shelah
\cite{Sh:630}, \cite{Sh:669} (giving yet more examples to which the general
theory developed there can be applied).

The second theme of the paper is a part of the general program ``how special
are random and Cohen forcing notions (or: the respective ideals)''. Kunen
(see \cite[Question 1.2]{Ku84}) suspected that the null ideal and the meager
ideal on $\can$ can be somehow characterized by their combinatorial
properties, but in \cite{RoSh:628} we constructed $\sigma$--ideals (or
rather forcing notions) that have nice properties, however are different
from the two. (But see also Kechris and Solecki \cite{KeSo95} and Solecki
\cite{So98} for results in the opposite direction.) Shelah \cite{Sh:480}
shows that the two forcing notions may occupy special positions in the realm
of nicely definable forcing notions. In this realm we may classify forcing
notions using the methods of \cite{Sh:630}, \cite{Sh:669} and, for example,
declare that very Souslin (or generally $\omega$--nw--nep) ccc forcing
notions (see \ref{verysouslin}) are really nice. Both the Cohen forcing
notion and the random forcing notion and their FS iterations (and nice
subforcings) are all ccc $\omega$--nw--nep, and \cite[Problem 4.24]{Sh:666}
asked if we have more examples. It occurs that our method relatively easily
results in very Souslin ccc forcing notions (see \ref{allsimple}(3),
\ref{conclOne}(2), \ref{extra}, \ref{conclThree}(3)). 

The third theme is {\em sweet {\&} sour} and it is related to one of the
most striking differences between the random and the Cohen forcing notions
which appears when we consider the respective regularity properties of
projective set. In \cite{Sh:176}, Shelah proved that the Lebesgue
measurability of $\Sigma^1_3$ sets implies $\omega_1$ is inaccessible in
${\bf L}$, while one can construct (in ZFC) a forcing notion  $\bP$ such
that $\bV^\bP\models$ ``projective subsets of $\mbR$ have the Baire property
''. The latter construction involved a strong version of ccc, so called
``sweetness'' (see \ref{sweet}). The heart of the former result is that the
composition of two Amoeba for measure forcing notions is sour (see
\ref{sour}) over random. Also from a sequence of $\omega_1$ reals we can
define a non-measurable set, but not one without the Baire Property. 

It seems that sweet--sour properties of forcing notions could be used to
classify them as either close to Cohen or as more random--like. Again, our
methods result in examples for both cases.  

Let us postpone the discussion of the general context of this paper till
{\em Epilogue}, when we can easier refer to the definitions and notions
discussed in the paper. (But the curious reader may start reading this paper
from that section.)

We try to make this work self contained, citing the most important
definitions and results from \cite{RoSh:470}, \cite{RoSh:628} whenever
needed. However, at least superficial familiarity with those papers could be
of some help in reading this paper. 

\subsection{The content of the paper} 
Like in \cite{RoSh:470}, the basic intention of this paper is to present
``the general theory'' rather than particular examples. Therefore, we
extract those properties of an example we want to construct which are
responsible for the fact that it works and we separate ``the general
theory'' from its applications. But to make the paper more readable, in
most cases, we sacrifice generality for clarity.  

In the first section we uniformize and generalize the constructions of
\cite{JRSh:373} and \cite{RoSh:628}. We investigate the complexity of the
resulting forcing notions as well as properties like ``adding unbounded
reals'', ``preserving unbounded families'', etc.

The next section introduces more ways in which creatures (or
tree--creatures) can be used to build ccc forcing notions. We discuss
mixtures with randoms, some generalizations of the Amoeba for Category
forcing notion, as well as as ``artificial'' modifications of previously
introduced forcings.

The third part formalizes definitions of $\sigma$--ideals corresponding to
our forcing notions. 

The following section discusses sweet--sour properties of our forcing
notions. We recall the notions of sweetness and introduce yet another sweet
property, and we show that very often our constructions are (somewhat)
sweet. However, there are exceptions to this rule. So we define some strong
negations of sweetness (sourness) and we show how our schema may end up with
very sour results.

Finally, the last section is (in some sense) a continuation of the
introduction. We discuss the results of the paper and formulate some
problems.

\subsection{Notation} 
Most of our notation is standard and compatible with that of classical
textbooks on Set Theory (like Bartoszy\'nski and Judah \cite{BaJu95}). 
However in forcing we keep the convention that {\em a stronger condition is
the larger one}.   
\medskip

\noindent {\sc Basic Notation:} In this paper $\bH$ will stand for a function
with domain $\omega$ and such that $(\forall m\in\omega)(|\bH(m)|\geq 2)$. 
We usually assume that $0\in \bH(m)$ (for all $m\in\omega$); if it is not
the case then we fix an element of $\bH(m)$ and we use it whenever
appropriate notions refer to $0$. Moreover we demand $\bH\in {\mathcal H}(
\omega_1)$ (i.e., $\bH$ is hereditarily countable. 
\medskip

\noindent {\sc More Notation:}
\begin{enumerate}
\item $\mbR^{{\geq}0}$ stands for the set of non-negative reals. The integer
part of a real $r\in\mbR^{{\geq}0}$ is denoted by $\lfloor r\rfloor$.
\item For two sequences $\eta,\nu$ we write $\nu\vartriangleleft\eta$
whenever $\nu$ is a proper initial segment of $\eta$, and $\nu
\trianglelefteq\eta$ when either $\nu\vartriangleleft\eta$ or $\nu=\eta$. 
The length of a sequence $\eta$ is denoted by $\lh(\eta)$.
\item A {\em tree} is a family $T$ of finite sequences such that for some
$\mrot(T)\in T$ we have
\[(\forall\nu\in T)(\mrot(T)\trianglelefteq \nu)\quad\mbox{ and }\quad
\mrot(T)\trianglelefteq\nu\trianglelefteq\eta\in T\ \Rightarrow\ \nu\in T.\]
For a tree $T$, the family of all $\omega$--branches through $T$ is
denoted by $[T]$, and we let
\[\max(T)\stackrel{\rm def}{=}\{\nu\in T:\mbox{ there is no }\rho\in
T\mbox{ such that }\nu\vartriangleleft\rho\}.\]   
If $\eta$ is a node in the tree $T$ then 
\[\begin{array}{lcl}
\suc_T(\eta)&=&\{\nu\in T: \eta\vartriangleleft\nu\ \&\ \lh(\nu)=\lh(\eta)+1
\}\ \mbox{ and}\\
T^{[\eta]}&=&\{\nu\in T:\eta\trianglelefteq\nu\}.
  \end{array}\]
\item The quantifiers $(\forall^\infty n)$ and $(\exists^\infty n)$ are
abbreviations for  
\[(\exists m\in\omega)(\forall n>m)\quad\mbox{ and }\quad(\forall m\in\omega)
(\exists n>m),\] 
respectively.
\item For a set $X$,\ \ \ $[X]^{\textstyle{\leq}\omega}$,
$[X]^{\textstyle{<}\omega}$ and ${\mathcal P}(X)$ will stand for families of
countable, finite and all, respectively, subsets of the set $X$. The family
of $k$-element subsets of $X$ will be denoted by $[X]^{\textstyle k}$. The
set of all finite sequences with values in $X$ is called $X^{\textstyle
{<}\omega}$ (so domains of elements of $X^{\textstyle {<}\omega}$ are
integers). The collection of all {\em finite partial} functions from
$\omega$ to $X$ is $X^{\mathunderaccent\smile-3 \omega}$. 
\item For a relation $R$ (a set of ordered pairs), $\dom(R)$ stands for the
domain of $R$ and $\rng(R)$ denotes the range of $R$. 
\item The Cantor space $\can$ and the Baire space $\baire$ are the spaces of
all functions from $\omega$ to $2$, $\omega$, respectively, equipped with
natural (Polish) topology. 
\item For $f,g\in\baire$ we write $f<^* g$ ($f\leq^* g$, respectively)
whenever $(\forall^\infty n\in\omega)(f(n)<g(n))$ ($(\forall^\infty n\in
\omega)(f(n)\leq g(n))$, repectively). 
\item For a forcing notion $\bP$, $\Gamma_\bP$ stands for the canonical
$\bP$--name for the generic filter in $\bP$. With this one exception, all
$\bP$--names for objects in the extension via $\bP$ will be denoted with a
dot above (e.g.~$\dot{\tau}$, $\dot{X}$). The weakest element of $\bP$ will
be denoted by $\emptyset_\bP$ (and we will always assume that there is one,
and that there is no other condition equivalent to it).
\end{enumerate}

\section{Building Souslin ccc forcing notions}
In this section we will review methods for building ccc forcing notions
announced or present in some form in \cite{JRSh:373}, \cite{RoSh:470}, and
\cite{RoSh:628}.  

\subsection{Glue and cut --- the method of [RoSh:628]}
Here we re-present the method of building ccc forcing notions with use of
(semi--) creating triples from \cite{RoSh:628}. We will slightly modify the
definitions loosing some generality. However, we will gain more direct
connection to the method of \cite{RoSh:470} and (hopefully) a better clarity 
of arguments. Note that the main difference is that here we do not worry
about ``the permutation invariance'' of our forcing notions, so the
creatures get back their $m^t_\dn,m^t_\up$ (and they are like those of  
\cite{RoSh:470}). 

\begin{definition}
\label{triples}
Let $\bH:\omega\longrightarrow {\mathcal H}(\omega_1)$.
\begin{enumerate}
\item (See \cite[Def.\ 1.1.1, 1.2.1]{RoSh:470}) {\em A creature for $\bH$}
is a triple 
\[t=(\nor,\val,\dis)=(\nor[t],\val[t],\dis[t])\]
such that $\nor\in\mbR^{{\geq}0}$, $\dis\in {\mathcal H}(\omega_1)$, and for
some integers $m^t_\dn<m^t_\up<\omega$
\[\emptyset\neq\val\subseteq\{\langle u,v\rangle\in\prod_{i<m^t_\dn}\bH(i)
\times\prod_{i<m^t_\up}\bH(i): u\vartriangleleft v\}.\]
The set of all creatures for $\bH$ will be denoted by $\CR[\bH]$, and for
$m_0<m_1<\omega$ we let $\CR_{m_0,m_1}[\bH]=\{t\in\CR[\bH]:m^t_\dn=m_0\ \&\ 
m^t_\up=m_1\}$. 
\item (See \cite[Def.\ 1.1.4, 1.2.2, 1.2.5]{RoSh:470}) Let $K\subseteq\CR[
\bH]$. We say that a function $\Sigma:[K]^{\textstyle{<}\omega}
\longrightarrow\cP(K)$ is a {\em composition operation on $K$} whenever the
following conditions are satisfied. 
\begin{enumerate}
\item[(a)] If $\cS\in [K]^{\textstyle{<}\omega}$ and $\Sigma(\cS)\neq
\emptyset$, then for some enumeration $\cS=\{t_0,\ldots,t_k\}$ we have
$m^{t_i}_\up=m^{t_{i+1}}_\dn$ for all $i<k$ [from now on, whenever we write
$\Sigma(t_0,\ldots,t_k)$, we mean the enumeration in which $m^{t_i}_{\up}=
m^{t_{i+1}}_\dn$], and 
\item[(b)] if $s\in\Sigma(t_0,\ldots,t_k)$, then $m^s_\dn=m^{t_0}_\dn$ and
$m^s_\up=m^{t_k}_\up$, and
\item[(c)] $t\in\Sigma(t)$ for each $t\in K$, $\Sigma(\emptyset)=\emptyset$,
and 
\item[(d)] {[{\em transitivity\/}]} if $s_i\in \Sigma(t^i_0,\ldots,t^i_{k_i}
)$ (for $i\leq n$), then 
\[\Sigma(s_0,\ldots,s_n)\subseteq \Sigma(t^i_j: i\leq n\ \&\ j\leq k_i),\] 
\item[(e)] {[{\em niceness\ \&\ smoothness\/}]} if $s\in\Sigma(t_0,\ldots,
t_k)$, $m^{t_i}_\up=m^{t_{i+1}}_\dn$ (for $i<k$), then $\dom(\val[t_0])=
\dom(\val[s])$ and 
\[(\forall\langle u,v\rangle\in\val[s])(\forall i\leq k)(\langle v\rest
m^{t_i}_\dn,v\rest m^{t_i}_\up\rangle\in\val[t_i]).\]
\end{enumerate}
\item (See \cite[Def.\ 1.1]{RoSh:628}) A function $\Sigma^\bot:K
\longrightarrow [K]^{\textstyle{<}\omega}\setminus\{\emptyset\}$ is called
{\em a decomposition operation on $K$} if for each $t\in K$:
\begin{enumerate}
\item[(a)$^\bot$] if $\cS\in\Sigma^\bot(t)$, then for some enumeration
$\cS=\{s_0,\ldots,s_k\}$ we have $m^{s_i}_\up=m^{s_{i+1}}_\dn$ (for $i<k$)
[from now on, if we write $\{s_0,\ldots,s_k\}\in\Sigma^\bot(t)$, we mean the
enumeration in which $m^{s_i}_\up=m^{s_{i+1}}_\dn$], and 
\item[(b)$^\bot$] if $\{s_0,\ldots,s_k\}\in\Sigma^\bot(t)$ then $m^{s_0}_\dn
=m^t_\dn$, $m^{s_k}_\up=m^t_\up$, 
\item[(c)$^\bot$] $\{t\}\in\Sigma^\bot(t)$,
\item[(d)$^\bot$] {[{\em transitivity\/}]} if $\cS=\{s_0,\ldots,s_k\}\in
\Sigma^\bot(t)$ and $\cS_i\in\Sigma^\bot(s_i)$ (for $i\leq k$), then $\cS_0
\cup\ldots\cup\cS_k\in \Sigma^\bot(t)$,
\item[(e)$^\bot$] if $\{s_0,\ldots,s_k\}\in\Sigma^\bot(t)$, $m^{s_i}_\up=
m^{s_{i+1}}_\dn$ (for $i<k$), then 
\[\dom(\val[t])=\dom(\val[s_0])\ \mbox{ and }\ (\forall i<k)\big(\rng(
\val[s_i])\subseteq\dom(\val[s_{i+1}])\big),\]
and 
\[\{\langle u,v\rangle: u\in\dom(\val[s_0])\ \&\ u\vartriangleleft v\ \&\
(\forall i\leq k)(\langle v\rest m^{s_i}_\dn,v\rest m^{s_i}_\up\rangle\in
\val[s_i])\}\subseteq\val[t].\]
\end{enumerate}
\item If $K\subseteq\CR[\bH]$ and $\Sigma$ is a composition operation on
$K$, then $(K,\Sigma)$ is called {\em a creating pair for $\bH$}. If,
additionally, $\Sigma^\bot$ is a decomposition operation on $K$, then
$(K,\Sigma,\Sigma^\bot)$ is called {\em a $\otimes$--creating triple for
$\bH$}. 
\item If $t_0,\ldots,t_n\in K$ are such that $m^{t_i}_\up=m^{t_{i+1}}_\dn$
(for $i<n$) and $w\in\dom(\val[t_0])$, then we let
\[\pos(w,t_0,\ldots,t_n)\stackrel{\rm def}{=}\{v\!\in\!\!\prod_{j<m^{
t_n}_\up}\!\!\bH(j)\!: w\vartriangleleft v\ \&\ (\forall i\leq n)(\langle
v\rest m^{t_i}_\dn, v\rest m^{t_i}_\up\rangle\in\val[t_i])\}.\]
\end{enumerate}
\end{definition}

\begin{definition}
\label{varia}
Let $(K,\Sigma,\Sigma^\bot)$ be a $\otimes$--creating triple for $\bH$. We
say that 
\begin{enumerate}
\item $\Sigma^\bot$ is {\em trivial\/} if $\Sigma^\bot(t)=\big\{\{t\}\big\}$
for each $t\in K$;
\item $(K,\Sigma)$ is {\em simple\/} if $\Sigma(\cS)$ is non-empty for
singletons only; if additionally $\Sigma^\bot$ is trivial, then we say that
$(K,\Sigma,\Sigma^\bot)$ is {\em simple\/};
\item $K$ (or $(K,\Sigma)$ or $(K,\Sigma,\Sigma^\bot)$) is {\em local\/} if
$m^t_\up=m^t_\dn +1$ for each creature $t\in K$ (so then necessarily $(K,
\Sigma,\Sigma^\bot)$ is simple);
\item $K$ is {\em forgetful\/} if for every creature $t\in K$ we have
\[[\langle u,v\rangle\in \val[t]\ \&\ w\in\prod_{i<m^t_\dn}\bH(i)]\ \
\Rightarrow\ \ \langle w,w\conc v\rest [m^t_\dn,m^t_\up)\rangle\in
\val[t];\]
\item $K$ is {\em full\/} if $\dom(\val[t])=\prod\limits_{i<m^t_\dn}\bH(i)$
for each $t\in K$. 
\end{enumerate}
\end{definition}

\begin{definition}
\label{forcing}
[See {\cite[Def.\ 1.1.7, 1.2.6]{RoSh:470}, \cite[Def.\ 1.3]{RoSh:628}}] Let
$(K,\Sigma,\Sigma^\bot)$ be a $\otimes$--creating triple for $\bH$ and let
$\cnor$ be a property of $\omega$--sequences of creatures from $K$ (so
$\cnor$ can be thought of as a subset of $K^{\textstyle\omega}$). We define
a forcing notion $\qcnor$ as follows.
\medskip

\noindent {\bf A condition} in $\qcnor$ is a sequence $p=(w^p,t^p_0,t^p_1,
t^p_2,\ldots)$ such that 
\begin{enumerate}
\item[(a)] $t^p_i\in K$ and $m^{t^p_i}_\up=m^{t^p_{i+1}}_\dn$ (for $i<
\omega$),
\item[(b)] $w\in\dom(\val[t^p_0])$ and $\langle t^p_0,t^p_1,t^p_2,\ldots
\rangle\in\cnor$,
\item[(c)] $\pos(w^p,t^p_0,\ldots,t^p_i)\subseteq \dom(\val[t^p_{i+1}])$ for
each $i<\omega$. 
\end{enumerate}
$\qzero$ is defined similarly, but we skip the demand ``$\langle t^p_0,t^p_1,
\ldots\rangle\in\cnor$'' in clause (b) above (or we just let $\cnor=K^{
\textstyle\omega}$; it is perhaps unfortunate to use $\emptyset$ in this
context, but that  notation was established in \cite{RoSh:470}).

\noindent{\bf The relation $\leq$} on $\qcnor$ is given by: \quad $p\leq q$ 
\quad if and only if\quad $(w^q,t^q_0,t^q_1,t^q_2,\ldots)$ can be obtained
from $(w^p,t^p_0,t^p_1,t^p_2,\ldots)$ by applying finitely many times the
following operations (describing the operations, we say what are the results
of applying the operation to a condition $(w,t_0,t_1,t_2,\ldots)\in\qzero$).  
\smallskip

{\em Deciding the value} for $(w,t_0,t_1,t_2,\ldots)$:

\noindent a result of this operation is a condition $(w^*,t_n,t_{n+1},
t_{n+2},\ldots)\in\qzero$ such that $w^*\in\pos(w,t_0,\ldots,t_{n-1})$ for
some $n<\omega$. 
\smallskip

{\em Applying $\Sigma$} to $(w,t_0,t_1,t_2,\ldots)$:

\noindent a result of this operation is a condition $(w,t^*_0,t^*_1,t^*_2,
\ldots)\in\qzero$ such that for some increasing sequence $0=n_0<n_1<n_2<
\ldots<\omega$, for each $i<\omega$, we have $t^*_i\in\Sigma(t_{n_i},
\ldots,t_{n_{i+1}-1})$.
\smallskip

{\em Applying $\Sigma^\bot$} to $(w,t_0,t_1,t_2,\ldots)$:

\noindent a result of this operation is a condition $(w,t^*_0,t^*_1,t^*_2,
\ldots)\in\qzero$ such that for some increasing sequence $0=n_0<n_1<n_2<
\ldots<\omega$, for each $i<\omega$, we have $\{t^*_{n_i},\ldots,t^*_{n_{i+
1}-1}\}\in\Sigma^\bot(t_i)$.
\end{definition}

\begin{remark}
\label{side}
In the definition of the relation $\leq$ on $\qcnor$ we do not require that
the intermediate steps satisfy the norm condition $\cnor$. So it may happen
that the sequence of witnesses for $p\leq q$ (i.e., the results of the
respective operations) is not in $\qcnor$.

If $\Sigma^\bot$ is trivial we may omit it; note that then we are exactly in
the setting of \cite[\S 1.2]{RoSh:470}.
\end{remark}

\begin{definition}
\label{nordef}
We will consider the following norm conditions $\cnor$:
\begin{itemize}
\item A sequence $\langle t_i:i<\omega\rangle$ satisfies ${\mathcal
C}^\infty(\nor)$ if $\lim\limits_{i\to\infty}\nor[t_i]=\infty$

[the respective forcing notion is called $\qinf$].
\item Let $\cF\subseteq\baire$; a sequence $\langle t_i:i<\omega\rangle$
satisfies ${\mathcal C}^\cF(\nor)$ if 
\[(\exists f\in\cF)(\forall^\infty i\in\omega)(\nor[t_i]\geq
f(m^{t_i}_\dn))\]
[the respective forcing notion is denoted $\qcF$].
\item Let $f:\omega\times\omega\longrightarrow\omega$; a sequence $\langle
t_i:i<\omega\rangle$ satisfies ${\mathcal C}^f(\nor)$ if 
\[(\forall n\in\omega)(\forall^\infty i\in\omega)(\nor[t_i]\geq f(n,m^{
t_i}_\dn))\]
[the respective forcing notion is denoted $\qonef$].
\end{itemize}
We will consider the norm conditions ${\mathcal C}^\cF(\nor)$, ${\mathcal
C}^f(\nor)$ only for $h$--closed families $\cF$ and fast functions $f$, see
\ref{functions} below. Later we will introduce more methods for building ccc
forcing notions, including more norm conditions.
\end{definition}

\begin{definition}
\label{functions}
\begin{enumerate}
\item A function $f:\omega\times\omega\longrightarrow\omega$ is {\em fast\/}
if 
\[(\forall k,\ell\in\omega)(f(k,\ell)\leq f(k,\ell+1)\ \&\ 2\cdot f(k,\ell)<
f(k+1,\ell)).\]
\item A function $h:\omega\times\omega\longrightarrow\omega$ is {\em
regressive} if 
\[(\forall m\in\omega)\big((\forall k>1)(1\leq h(m,k)<k)\ \&\ (\forall
k<\ell<\omega)(h(m,k)\leq h(m,\ell))\big).\]
\item Let $h:\omega\times\omega\longrightarrow\omega$. We say that a family
$\cF\subseteq\baire$ is {\em $h$--closed\/} if for every $f\in\cF$ there is
$f^*\in\cF$ such that $(\forall^\infty n\in\omega)(f^*(n)\leq h(n,f(n)))$. 
\item A family $\cF\subseteq\baire$ is {\em $\geq^*$--directed\/} if
\[(\forall f_0,f_1\in\cF)(\exists f^*\in\cF)(\forall^\infty n\in\omega)(f^*
(n)\leq \min\{f_0(n),f_1(n)\}).\]
Similarly we define $\leq^*$--directed families (just reversing the
inequality).  
\end{enumerate}
\end{definition}

\begin{remark}
\label{rem1}
Let $f(n,m)=2^{2n}$ (for $n,m\in\omega$). Then the function $f$ is fast and
the norm conditions ${\mathcal C}^f(\nor)$ and ${\mathcal C}^\infty(\nor)$
agree (and thus $\bQ^*_f(K,\Sigma)=\bQ^*_\infty(K,\Sigma)$ for a local
creating pair $(K,\Sigma)$). In practical applications, when we consider the
norm condition ${\mathcal C}^f(\nor)$, the function $f$ is such that $f(n,m)
<f(n,m+1)$ (for all $n,m\in\omega$) and thus the norm condition ${\mathcal
C}^f(\nor)$ is stronger than ${\mathcal C}^\infty(\nor)$.
\end{remark}

\begin{proposition}
If $(K,\Sigma,\Sigma^\bot)$ is a $\otimes$--creating triple for $\bH$,
$\cnor\subseteq K^{\textstyle\omega}$, then $\qcnor$ is a forcing notion
(i.e., the relation $\leq$ of $\qcnor$ is transitive).
\end{proposition}

\begin{definition}
\label{candidates}
[See {\cite[Def.\ 1.2.4]{RoSh:470}}]
Let $(K,\Sigma,\Sigma^\bot)$ be a $\otimes$--creating triple for $\bH$. We
define {\em finite candidates} ($\FC$) and {\em pure finite candidates}
($\PFC$) with respect to $(K,\Sigma,\Sigma^\bot)$: 
\[\hspace{-0.1cm}\begin{array}{ll}
\FC(K,\Sigma,\Sigma^\bot)=&\{(w,t_0,\ldots,t_n):
w\in\dom(\val[t_0])\mbox{ and for each }i\leq n:\\
\ & \ t_i\in K, m^{t_i}_{\up} = m^{t_{i + 1}}_{\dn}\mbox{ and }
\pos(w,t_0,\ldots,t_{i})\subseteq\dom(\val[t_{i+1}])\},\\
\end{array}\]
\[\PFC(K,\Sigma,\Sigma^\bot)=\{(t_0,\ldots,t_n):(\exists w\!\in\!\dom(\val[
t_0]))((w,t_0, \ldots,t_n)\in \FC(K,\Sigma,\Sigma^\bot))\}.\]  
We have a natural partial order $\leq$ on $\FC(K,\Sigma,\Sigma^\bot)$ (like
in \ref{forcing}). [Note that $\Sigma,\Sigma^\bot$ have no influence on
$\FC(K,\Sigma,\Sigma^\bot)$, that is they are not present in the definition
of finite candidates, and we could have written $\FC(K,\Sigma)$ or $\FC(K)$. 
However, they come to game when the partial order $\leq$ on $\FC(K,\Sigma,
\Sigma^\bot)$ is considered.]

A sequence $\langle t_0,t_1,t_2,\ldots\rangle$ of creatures from $K$ is 
{\em a pure candidate} with respect to $(K,\Sigma,\Sigma^\bot)$ if  
\[(\forall i<\omega)(m^{t_i}_{\up}=m^{t_{i+1}}_{\dn})\quad\mbox{ and}\]
\[(\exists w\in\dom(\val[t_0]))(\forall i<\omega)(\pos(w,t_0,\ldots,t_i)
\subseteq\dom(\val[t_{i+1}])).\] 
The set of pure candidates with respect to $(K,\Sigma)$ is denoted by
$\PC(K,\Sigma,\Sigma^\bot)$. The partial order $\leq$ on $\PC(K,\Sigma,
\Sigma^\bot)$ is defined naturally.

For a norm condition ${\mathcal C}(\nor)$ the family of {\em ${\mathcal
C}(\nor)$--normed pure candidates} is
\[\PC_{{\mathcal C}(\nor)}(K,\Sigma,\Sigma^\bot)=\{\langle t_0,t_1,\ldots
\rangle\in\PC(K,\Sigma,\Sigma^\bot)\!:\langle t_0,t_1,\ldots,\rangle \mbox{
satisfies }{\mathcal C}(\nor)\}. \] 
\end{definition}

\begin{definition}
\label{POS}
Let $(K,\Sigma,\Sigma^\bot)$ be a $\otimes$--creating triple for $\bH$.
\begin{enumerate}
\item For a condition $p\in\qzero$ we let
\[\POS(p)\stackrel{\rm def}{=}\{u\in\bigcup_{n<\omega}\prod_{i<n}\bH(i):
(\exists\ell<\omega)(\exists v\in\pos(w^p,t^p_0,\ldots,t^p_\ell))(u
\trianglelefteq v)\}.\]
\item For a finite candidate $c=(w,t_0,\ldots,t_k)\in\FC(K,\Sigma,
\Sigma^\bot)$ we define 
\[\POS(c)\stackrel{\rm def}{=}\{u\in\bigcup_{n\leq m^{t_k}_\up} \prod_{i<n}
\bH(i): (\exists v\in\pos(w,t_0,\ldots,t_k))(u\trianglelefteq v)\}.\]
\end{enumerate}
\end{definition}

\begin{proposition}
\label{sidePOS}
Suppose $(K,\Sigma,\Sigma^\bot)$ is a $\otimes$--creating triple for $\bH$.
\begin{enumerate}
\item If $p,q\in\qzero$, $p\leq q$ then $\POS(p)\subseteq\POS(q)$, and {\em
if\/} $\lh(w^q)=m^{t^p_\ell}_\up$ for some $\ell<\omega$, {\em then\/}
$w^q\in\pos(w^p,t^p_0,\ldots,t^p_\ell)$.
\item The same holds if one replaces conditions from $\qzero$ by finite
candidates from $\FC(K,\Sigma,\Sigma^\bot)$. 
\end{enumerate}
\end{proposition}

\begin{proof}
Note that each of the three operations described in \ref{forcing} shrinks
$\POS$ (remember \ref{triples}(2e) and \ref{triples}(3e$^\bot$)). 
\end{proof}

\begin{definition}
\label{forlin}
[See {\cite[Def.\ 2.1]{RoSh:628}, \cite[Def.\ 2.1.7]{RoSh:470}}] Assume that
$(K,\Sigma,\Sigma^\bot)$ is a $\otimes$--creating triple for $\bH$.
\begin{enumerate}
\item We say that $(K,\Sigma)$ (or $(K,\Sigma,\Sigma^\bot)$) is {\em
linked\/} if for each $t_0,t_1\in K$ such that $\nor[t_0],\nor[t_1]>1$ and
$m^{t_0}_\dn=m^{t_1}_\dn$, $m^{t_0}_\up=m^{t_1}_\up$, there is $s\in\Sigma(
t_0)\cap\Sigma(t_1)$ with 
\[\nor[s]\geq\min\{\nor[t_0],\nor[t_1]\}-1.\]
Let $h:\omega\times\omega\longrightarrow\omega$. The pair $(K,\Sigma)$ is
said to be {\em $h$--linked} if for each $k>1$, and creatures $t_0,t_1\in K$
such that $\nor[t_0],\nor[t_1]\geq k$ and $m^{t_0}_\dn=m^{t_1}_\dn$,
$m^{t_0}_\up=m^{t_1}_\up$, there is $s\in\Sigma(t_0)\cap\Sigma(t_1)$ with
$\nor[s]\geq h(m^{t_0}_\dn,k)$. 
\item We say that $(K,\Sigma)$ (or $(K,\Sigma,\Sigma^\bot)$) is {\em gluing
\/} if it is full and for each $k<\omega$ there is $n_0=n_0(k)<\omega$ such
that for every $n\geq n_0$ and $(t_0,\ldots,t_n)\in\PFC(K,\Sigma)$, there is 
$s\in\Sigma(t_0,\ldots,t_n)$ such that 
\[\nor[s]\geq\min\{k,\nor[t_0],\ldots,\nor[t_n]\}.\]
\item We say that $(K,\Sigma^\bot)$ (or $(K,\Sigma,\Sigma^\bot)$) {\em has
the cutting property\/} if for every $t\in K$ with $\nor[t]>1$ and an
integer $m\in (m^t_\dn,m^t_\up)$, there are $s_0,s_1\in K$ such that
\begin{enumerate}
\item[$(\alpha)$] $m^{s_0}_\dn=m^t_\dn$, $m^{s_0}_\up=m=m^{s_1}_\dn$,
$m^{s_1}_\up=m^t_\up$, 
\item[$(\beta)$]  $\nor[s_\ell]\geq\min\{\nor[t]-1,m^t_\dn\}$ (for $\ell=0,
1$),
\item[$(\gamma)$] $\{s_0,s_1\}\in\Sigma^\bot(K)$.
\end{enumerate}
\end{enumerate}
\end{definition}

\begin{definition}
\label{linked}
A forcing notion $\bQ$ is {\em $\sigma$-$n$--linked\/} if there is a
partition $\langle A_i: i<\omega\rangle$ of $\bQ$ such that 
\[\mbox{if }\quad q_0,\ldots,q_{n-1}\in A_i,\ i\in\omega\quad\mbox{ then
}\quad (\exists q\in\bQ)(q_0\leq q\ \&\ \ldots\ \&\ q_{n-1}\leq q).\]
We say that $\bQ$ is {\em $\sigma$-$*$--linked\/} if it is
$\sigma$-$n$--linked for every $n\in\omega$.
\end{definition}

\begin{proposition}
\label{cccthm}
Let $\bH:\omega\longrightarrow{\mathcal H}(\omega_1)$ and let $(K,\Sigma,
\Sigma^\bot)$ be a $\otimes$--creating triple for $\bH$.
\begin{enumerate}
\item If $(K,\Sigma,\Sigma^\bot)$ is linked, gluing and has the cutting
property, then the forcing notion $\qinf$ is $\sigma$-$*$--linked.
\item If $f:\omega\times\omega\longrightarrow\omega$ is fast and $(K,\Sigma,
\Sigma^\bot)$ is local and linked, then the forcing notions $\qinf$ and
$\qonef$ are $\sigma$-$*$--linked. 
\item Assume that $h:\omega\times\omega\longrightarrow\omega$ is regressive 
and $\cF\subseteq (\omega\setminus 2)^{\textstyle\omega}$ is an $h$--closed
family which is either countable, or $\geq^*$--directed. Suppose $(K,\Sigma, 
\Sigma^\bot)$ is local and $h$--linked. Then the forcing notion $\qcF$ is
$\sigma$-$*$--linked.    
\end{enumerate}
\end{proposition}

\begin{proof}
Straightforward (and the proof of the first part is essentially the same as
that of \cite[Thm\ 2.4]{RoSh:628}; compare the proof of \ref{cl1}).
\end{proof}

\subsection{Tree--like conditions} Here we recall the setting of \cite[\S
1.3]{RoSh:470} and \cite{JRSh:373}. Since in getting the ccc we will have to
require that the tree--creating pair under considerations is local, we will
restrict our attention to that case only. So our definitions here are much
simpler than those in the general case, but we still try to keep the
notation and flavour of the tree case of \cite{RoSh:470}. 

\begin{definition}
\label{treecreature}
Let $\bH:\omega\longrightarrow {\mathcal H}(\omega_1)$.
\begin{enumerate}
\item {\em A local tree--creature for $\bH$\/} is a triple 
\[t=(\nor,\val,\dis)=(\nor[t],\val[t],\dis[t])\]
such that $\nor\in\mbR^{{\geq}0}$, $\dis\in {\mathcal H}(\omega_1)$, and for
some sequence $\eta\in\prod\limits_{i<n}\bH(i)$, $n<\omega$, we
have 
\[\emptyset\neq\val\subseteq\{\langle \eta,\nu\rangle: \eta\vartriangleleft
\nu\in \prod_{i\leq n}\bH(i)\}.\] 
For a tree--creature $t$ we let $\pos(t)\stackrel{\rm def}{=}\rng(\val[t])$.

The set of all local tree--creatures for $\bH$ will be denoted by
$\TCR[\bH]$, and for $\eta\in\bigcup\limits_{n<\omega}\prod\limits_{i<n}
\bH(i)$ we let $\TCR_\eta[\bH]=\{t\in\TCR[\bH]:\dom(\val[t])=\{\eta\}\}$. 
\item Let $K\subseteq\TCR[\bH]$. We say that a function $\Sigma:K
\longrightarrow\cP(K)$ is a {\em local tree composition on $K$\/} whenever
the following conditions are satisfied. 
\begin{enumerate}
\item[(a)] If $t\in\TCR_\eta[\bH]$, $\eta\in\prod\limits_{i<n}\bH(i)$,
$n<\omega$, then $\Sigma(t)\subseteq\TCR_\eta[\bH]$. 
\item[(b)] If $s\in\Sigma(t)$ then $\val[s]\subseteq\val[t]$.
\item[(c)] {[{\em transitivity\/}]} If $s\in\Sigma(t)$ then $\Sigma(s)
\subseteq\Sigma(t)$.
\end{enumerate}
\item If $K\subseteq\TCR[\bH]$ and $\Sigma$ is a (local) tree composition
operation on $K$ then $(K,\Sigma)$ is called {\em a (local) tree--creating
pair for $\bH$}. 
\end{enumerate}
\end{definition}

\begin{definition}
\label{treeforc}
Let $(K,\Sigma)$ be a (local) tree--creating pair for $\bH$.
\begin{enumerate}
\item We define the forcing notion $\bQ^{\tree}_1(K,\Sigma)$ as follows.
\medskip

\noindent {\bf A condition} is a system $p=\langle t_\eta:\eta\in T\rangle$ 
such that 
\begin{enumerate}
\item[(a)] $T\subseteq\bigcup\limits_{n\in\omega}\prod\limits_{i<n}\bH(i)$ is
a non-empty tree with $\max(T)=\emptyset$,
\item[(b)] $t_\eta\in\TCR_\eta[\bH]\cap K$ and $\pos(t_\eta)=\suc_T(\eta)$
(for $\eta\in T$),
\item[(c)${}_1$] for every $\eta\in[T]$ we have:
\[\mbox{the sequence }\langle\nor[t_{\eta\rest k}]:\lh(\mrot(T))\leq
k<\omega\rangle\mbox{ diverges to infinity.}\] 
\end{enumerate}

\noindent{\bf The order} is given by:

\noindent $\langle t^1_\eta: \eta\in T^1\rangle\leq\langle t^2_\eta:\eta\in
T^2\rangle$ if and only if 

\noindent $T^2\subseteq T^1$ and $t^2_\eta\in\Sigma(t^1_\eta)$ for each
$\eta\in T^2$.  
\medskip

\noindent If $p=\langle t_\eta:\eta\in T\rangle$, then we write
$\mrot(p)=\mrot(T)$, $T^p= T$, $t^p_\eta = t_\eta$ etc.
  
\item Similarly, we define forcing notions $\bQ^{\tree}_\cF(K,\Sigma)$ for
a family $\cF\subseteq\baire$ and $\bQ^{\tree}_f(K,\Sigma)$ for a function
$f:\omega\times\omega\longrightarrow\omega$, replacing the condition {\rm 
(c)${}_1$} by {\rm (c)${}_\cF$}, {\rm (c)${}_f$}, respectively, where:
\begin{enumerate}
\item[(c)${}_\cF$]  $(\exists f\in\cF)(\exists N<\omega)(\forall\eta\in T)(
\lh(\eta)\geq N\ \Rightarrow\ \nor[t_\eta]\geq f(\lh(\eta)))$,
\item[(c)${}_f$]  $(\forall n\in\omega)(\exists N<\omega)(\forall\eta\in T)
(\lh(\eta)\geq N\ \Rightarrow\ \nor[t_\eta]\geq f(n,\lh(\eta)))$.
\end{enumerate}
\item If $p\in\bQ^{\tree}_x(K,\Sigma)$ then, for $\eta\in T^p$, we let
$p^{[\eta]}=\langle t^p_\nu:\nu\in (T^p)^{[\eta]}\rangle$.
\end{enumerate}
\end{definition}

\begin{definition}
\label{treelin}
Assume that $(K,\Sigma)$ is a tree--creating pair for $\bH$.
\begin{enumerate}
\item We say that $(K,\Sigma)$ is {\em linked\/} if for each $\eta\in
\bigcup\limits_{n<\omega}\prod\limits_{i<n}\bH(i)$ and tree--creatures $t_0,
t_1\in K\cap\TCR_\eta[\bH]$ with $\nor[t_0],\nor[t_1]>1$, there is $s\in
\Sigma(t_0)\cap\Sigma(t_1)$ such that $\nor[s]\geq\min\{\nor[t_0],\nor[t_1]
\}-1$.
\item Let $h:\omega\times\omega\longrightarrow\omega$. The pair $(K,\Sigma)$
is {\em $h$--linked} if for each $t_0,t_1\in K\cap\TCR_\eta[\bH]$ such that
$\nor[t_0],\nor[t_1]\geq k$, $k>1$, there is $s\in\Sigma(t_0)\cap\Sigma(
t_1)$ with $\nor[s]\geq h(\lh(\eta),k)$. 
\end{enumerate}
\end{definition}

\begin{proposition}
\label{ccctree}
Let $\bH:\omega\longrightarrow{\mathcal H}(\omega_1)$ and let $(K,\Sigma)$
be a local tree--creating pair for $\bH$. 
\begin{enumerate}
\item If $f:\omega\times\omega\longrightarrow\omega$ is fast and
$(K,\Sigma)$ is linked, then the forcing notions $\bQ^{\tree}_1(K,\Sigma)$
and  $\bQ^{\tree}_f(K,\Sigma)$ are $\sigma$-$*$--linked. 
\item Assume that $h:\omega\times\omega\longrightarrow\omega$ is regressive 
and a family $\cF\subseteq (\omega\setminus 2)^{\textstyle\omega}$ is
$h$--closed and either countable, or $\geq^*$--directed. Suppose
$(K,\Sigma)$ is $h$--linked. Then the forcing notion $\bQ^{\tree}_\cF(K,
\Sigma)$ is $\sigma$-$*$--linked.    
\end{enumerate}
\end{proposition}

\begin{proof}
Straightforward.
\end{proof}

\subsection{The complexity of our forcing notions}

\begin{definition}
\label{verysouslin}
\begin{enumerate}
\item A forcing notion $(\bP,\leq_\bP)$ is {\em Souslin\/} ({\em Borel},
respectively) if $\bP$, $\leq_\bP$ and the incompatibility relation
$\bot_{\bP}$ are $\Sigma^1_1$ (Borel, respectively) subsets of $\mbR$ and
$\mbR\times\mbR$.  
\item A forcing notion $(\bP,\leq_\bP)$ is {\em very Souslin ccc\/} ({\em
very Borel ccc}, respectively), if it is Souslin (Borel, resp.), satisfies
the ccc and the notion 
\[\mbox{`` }\langle r_n:n<\omega\rangle\mbox{ is a maximal antichain ''}\]
is $\Sigma^1_1$ (Borel, resp.)
\end{enumerate}
\end{definition}
On Souslin forcing notions and their applications see Judah and Shelah
\cite{JdSh:292} and Goldstern and Judah \cite{GoJu92} (the results of these
two and many other papers on the topic are presented in Bartoszy\'nski and
Judah \cite{BaJu95}). A systematic treatment of definable forcing notions is
presented in \cite{Sh:630}, \cite{Sh:669} (note that very Souslin ccc
forcing notions are $\omega$--nw--nep). Here we are going to show that the
forcing notions built according to the schemes presented above typically are 
Borel ccc and (sometimes) even very Borel ccc. Thus we have tools for
constructing new ccc $\omega$--nw--nep forcing notions (the only examples 
known before were those coming from random forcing, the Cohen forcing and
their FS iterations; see \cite[\S 4]{Sh:666} for a discussion of this topic). 
Note that, by Shelah \cite{Sh:711}, ccc $\omega$--nw--nep forcing notions
cannot add dominating reals. Thus the forcing notions that are covered by
\ref{adddom} cannot be represented as very Souslin ccc forcing notions.

\begin{definition}
\label{regular}
A $\otimes$--creating triple $(K,\Sigma,\Sigma^\bot)$ for $\bH$ is {\em
regular\/} if the following condition is satisfied.
\begin{enumerate}
\item[$(\boxdot)$] Assume $(w,t_0,\ldots,t_n), (u,s_0,\ldots,s_m)\in\FC(K,
\Sigma,\Sigma^\bot)$ are such that 
\begin{itemize}
\item $m^{t_\ell}_\dn<m^{s_0}_\dn<m^{t_\ell}_\up\leq m^{s_0}_\up$ for some
$\ell\leq n$,  
\item $\nor[t_\ell]\geq 3$ (for the $\ell$ as above), and 
\item $(w,t_0,\ldots,t_n)\leq (u,s_0,\ldots,s_m)$, $m^{s_m}_\up\leq
m^{t_n}_\up$, and $\nor[s_0]\geq 3$.
\end{itemize}
Then there are $t',t''$ such that $\{t',t''\}\in \Sigma^\bot(t_\ell)$,
$m^{t'}_\up=m^{s_0}_\dn=m^{t''}_\dn$, $\nor[t'']\geq 2$ and $u\in\pos(w,t_0,
\ldots,t_{\ell-1},t')$. 
\end{enumerate}
\end{definition}

\begin{definition}
\label{finitary}
Let $\bH:\omega\longrightarrow{\mathcal H}(\omega_1)$ and let $(K,\Sigma)$
be either a creating pair for $\bH$ or a (local) tree--creating pair for
$\bH$. We say that $(K,\Sigma)$ is {\em really finitary\/} if the
following conditions are satisfied: 
\begin{enumerate}
\item[(a)] $\bH(n)$ is finite for all $n<\omega$ (so $\val[t]$ is finite for
all $t\in K$), and 
\item[(b)] for each $n\in\omega$, the set $\{t\in K:\rng(\val[t])\subseteq
\prod\limits_{i<n}\bH(i)\}$ is finite.
\end{enumerate}
\end{definition}

\begin{theorem}
\label{allsimple}
Let $\bH:\omega\longrightarrow{\mathcal H}(\omega_1)$.
\begin{enumerate}
\item Let $(K,\Sigma,\Sigma^\bot)$ be a $\otimes$--creating triple for
$\bH$ such that $K$ is countable. 
\begin{enumerate}
\item[(a)] If $(K,\Sigma,\Sigma^\bot)$ is regular, linked, gluing and has
the cutting property, then the forcing notion $\qinf$ is Souslin ccc.
\item[(b)] If $f:\omega\times\omega\longrightarrow\omega$ is fast and
$(K,\Sigma,\Sigma^\bot)$ is local and linked then $\qinf$ and $\qonef$ are
Borel ccc.  
\item[(c)] Assume that $h:\omega\times\omega\longrightarrow\omega$ is a
regressive function and $\cF\subseteq (\omega\setminus 2)^{\textstyle
\omega}$ is a countable $h$--closed family which is $\geq^*$--directed. If
$(K,\Sigma,\Sigma^\bot)$ is local and $h$--linked, then $\qcF$ is Borel ccc.
\end{enumerate}
\item Assume that $(K,\Sigma)$ is a local tree--creating pair for $\bH$ and
$K$ is countable. 
\begin{enumerate}
\item[(a)] If $f:\omega\times\omega\longrightarrow\omega$ is fast and
$(K,\Sigma)$ is linked, then $\bQ^{\tree}_f(K,\Sigma)$ is Borel ccc.
\item[(b)] Suppose that $h:\omega\times\omega\longrightarrow\omega$ is
regressive and $\cF\subseteq (\omega\setminus 2)^{\textstyle\omega}$ is a
countable $h$--closed family which is $\geq^*$--directed. If $(K,\Sigma)$ is
$h$--linked, then $\bQ^{\tree}_\cF(K,\Sigma)$ is Borel ccc.   
\end{enumerate}
\item If in 1(c) and 2(b) above the pair $(K,\Sigma)$ is really finitary,
then the respective forcing notions are very Borel ccc.
\end{enumerate}
\end{theorem}

\begin{proof}
1(a)\quad Let $\cX=(\bigcup\limits_{n<\omega}\prod\limits_{i<n}\bH(i))
\times K^{\textstyle\omega}$ be equipped with the product topology (of
countably many countable discrete spaces). So $\cX$ is a Polish space and it
should be clear that $\qzero$, $\qinf$ are its Borel subsets. To express
``$p\leq q$'' we have to say that {\em there is\/} a sequence $p=p_0,\ldots,
p_n=q$ of elements of $\qzero$ such that $p_{i+1}$ is obtained from $p_i$ by
one of the operations described in \ref{forcing}. Each of these operations
corresponds to a Borel subset of $\cX\times\cX$, so easily $\leq_{\qzero}$,
$\leq_{\qinf}$ are $\Sigma^1_1$ subsets of $\cX\times\cX$. The main
difficulty is to show that the incompatibility relation $\bot_{\qinf}$ is a
$\Sigma^1_1$ subset of $\cX\times\cX$. But this follows from the following
observation (note that this is the place where we use the assumption that
$(K,\Sigma,\Sigma^\bot)$ is regular). 
\begin{claim}
\label{cl1}
Conditions $p,q\in\qinf$ are compatible\quad if and only if\quad there
are $N,\ell,m<\omega$, $t_0',t_1',t_0'',t_1''\in K$ and $u$ such that
\begin{itemize}
\item $m^{t^p_\ell}_\dn\leq N<m^{t^p_\ell}_\up$,\quad $m^{t^q_m}_\dn\leq
N<m^{t^q_m}_\up$,
\item $\{t_0',t_1'\}\in\Sigma^\bot(t^p_\ell)$,\quad $\{t_0'',t_1''\}\in
\Sigma^\bot(t^q_m)$,\quad $m^{t_1'}_\dn=m^{t_1''}_\dn=N$,
\item $\nor[t_1']\geq 2$,\quad $\nor[t_1'']\geq 2$,
\item $u\in\pos(w^p,t^p_0,\ldots,t^p_{\ell-1},t_0')\cap\pos(w^q,t^q_0,
\ldots,t^q_{m-1},t_0'')$,
\item $(\forall n>\ell)(\nor[t^p_n]\geq 2)$ and $(\forall n>m)(\nor[t^q_n]
\geq 2)$.  
\end{itemize}
(If $N=m^{t^p_\ell}_\dn$ then $t'_0$ is not present; similarly on the $q$
side.) 
\end{claim}

\begin{proof}[Proof of the claim] 
First assume that conditions $p,q\in\qinf$ are compatible and let
$r\in\qinf$ be stronger than both $p$ and $q$. Passing to a stronger
condition we may demand that if $\ell,m$ are such that 
\[m^{t^p_\ell}_\dn\leq m^{t^r_0}_\dn<m^{t^p_\ell}_\up,\qquad m^{t^q_m}_\dn
\leq m^{t^r_0}_\dn<m^{t^q_m}_\up\]
then $m^{t^p_\ell}_\up\leq m^{t^r_0}_\up$, $m^{t^q_m}_\up\leq m^{t^r_0}_\up$
and that $\nor[t^r_0]\geq 5$ and
\[m^{t^p_n}_\up\geq\lh(w^r)\quad\Rightarrow\quad \nor[t^p_n]\geq 5,\quad
\mbox{and}\quad m^{t^q_n}_\up\geq\lh(w^r)\quad\Rightarrow\quad \nor[t^q_n]
\geq 5.\]
Now we may apply the regularity of $(K,\Sigma,\Sigma^\bot)$ (see
\ref{regular}) to get $t_0',t_1',t_0'',t_1''\in K$ such that 
\[\begin{array}{l}
\{t_0',t_1'\}\in\Sigma^\bot(t^p_\ell),\quad \{t_0'',t_1''\}\in\Sigma^\bot(
t^q_m),\quad\nor[t_1']\geq 2,\quad \nor[t_1'']\geq 2\quad\mbox{ and}\\
u=w^r\in\pos(w^p,t^p_0,\ldots,t^p_{\ell-1},t_0')\cap\pos(w^q,t^q_0,\ldots,
t^q_{m-1},t_0'').
  \end{array}\]
Put $N=\lh(w^r)$ and check that all demands are satisfied.

For the other implication suppose that $N,\ell,m,t_0',t_1',t_0'',t_1''$ and
$u$ are as in the second statement. Choose increasing sequences $\langle
n_i:i<\omega\rangle$ and $\langle k_i:i<\omega\rangle$ such that $n_0>\ell+
5$, $k_0>m+5$ and
\begin{itemize}
\item $(\forall n\geq n_i)(\nor[t^p_n]\geq i+5)$ and $(\forall n\geq k_i)(
\nor[t^q_n]\geq i+5)$,\quad and 
\item $m^{t^q_{k_i}}_\dn\leq m^{t^p_{n_i}}_\dn<m^{t^q_{k_i}}_\up$,\quad
$m^{t^q_{k_i+1}}_\up<m^{t^p_{n_{i+1}}}_\dn$.
\end{itemize}
Apply the cutting property to choose (for each $i<\omega$) $s_i',s_i''\in K$
such that 
\[\{s_i',s_i''\}\in\Sigma^\bot(t^q_{k_i}),\quad m^{s_i'}_\dn=m^{t^q_{
k_i}}_\dn,\quad m^{s_i''}_\dn=m^{t^p_{n_i}}_\dn,\quad\mbox{and}\quad 
\nor[s_i'],\nor[s_i'']\geq i+4.\]
(If $m^{t^p_{n_i}}_\dn=m^{t^q_{k_i}}_\dn$ then the $s_i'$ is not present.)
Next use gluing to choose $r_i,s_i$ so that 
\[\begin{array}{ll}
r_0\in\Sigma(t'_1,t^p_{\ell+1},\ldots,t^p_{n_0-1}), &s_0\in\Sigma(t_1'',
t^q_{m+1},\ldots,t^q_{k_0-1},s_0')\\
r_{i+1}\in\Sigma(t^p_{n_i},\ldots,t^p_{n_{i+1}-1}), &s_{i+1}\in\Sigma(
s_i'',t^q_{k_i+1},\ldots,t^q_{k_{i+1}-1},s_{i+1}'),\\
\nor[r_i],\nor[s_i]\geq i+2. &\ \\
  \end{array}\]
Since $(K,\Sigma)$ is linked we may choose $t_i\in\Sigma(r_i)\cap\Sigma(
s_i)$ such that $\nor[t_i]\geq i+1$. Now look at $(u,t_0,t_1,\ldots)$. It is
a condition in $\qinf$ stronger than both $p$ and $q$. 
\end{proof}
1(b,c) and 2(a,b)\quad Similarly (and much easier).
\medskip

\noindent 3.\quad Let $h\in\baire$ be a regressive function and let
$\cF\subseteq (\omega\setminus 2)^{\textstyle\omega}$ be a countable
$h$--closed family which is $\geq^*$--directed. Suppose that $(K,\Sigma)$ is
a local, $h$--linked and really finitary creating pair (because of the
``local'' $\Sigma^\bot$ can be omitted as it is trivial). We are going to
show that ``being a (countable) pre-dense subset of $\qscF$'' is a Borel
property.

Let $\cX=(\bigcup\limits_{n<\omega}\prod\limits_{i<n}\bH(i))\times
K^{\textstyle\omega}$, $\cX^{\textstyle\omega}$ and $\cY=\cP(\FC(K,\Sigma))$
be equipped with the natural (product) Polish topologies (note that
$\FC(K,\Sigma)$ is a countable set). For $\bar{p}=\langle p_n: n<\omega
\rangle\in\cX^{\textstyle\omega}$, $\bar{p}\subseteq\qscF$, $w\in
\bigcup\limits_{m<\omega}\prod\limits_{i<m}\bH(i)$ and $f\in\cF$ we define
\[N^{\bar{p}}(n)=\min\{m^{t^{p_n}_i}_\dn: (\forall j\geq i)(\nor[t^{p_n}_j]
\geq 2)\},\]
and
\[\begin{array}{ll}
\cT^{\bar{p}}_{w,f}=&\{(w,t_0,\ldots,t_k)\in\FC(K,\Sigma): (\forall i\leq
k)(\nor[t_i]\geq f(m^{t_i}_\dn))\mbox{ and }\ \\  
&\quad (\forall n<\omega)(N^{\bar{p}}(n)\leq m^{t_k}_\up\ \ \Rightarrow\ \
\pos(w,t_0,\ldots,t_k)\cap \POS(p_n)=\emptyset)\}.
  \end{array}\]
Note that $(\qscF)^{\textstyle\omega}$ is a Borel subset of $\cX^{\textstyle
\omega}$ and the functions
\[\bar{p}\mapsto N^{\bar{p}}:(\qscF)^{\textstyle\omega}\longrightarrow
\baire\quad\mbox{ and }\quad \bar{p}\mapsto\cT^{\bar{p}}_{w,f}:(\qscF)^{
\textstyle\omega}\longrightarrow\cY\]
are Borel. Now, each $\cT^{\bar{p}}_{w,f}$ is essentially a finitary tree,
so 
\[\cT^{\bar{p}}_{w,f}\mbox{ is well founded\quad \ \ if and only if\quad \ \
} \cT^{\bar{p}}_{w,f}\mbox{ is finite.}\]
Consequently, for each $w$ and $f$, the set 
\[\{\bar{p}\in (\qscF)^{\textstyle\omega}: \cT^{\bar{p}}_{w,f}\mbox{ is well
founded }\}\]
is Borel. Since there are countably many possibilities for $w$ and $f$, we
easily finish the proof using the following observation.
\begin{claim}
\label{cl2}
Let $\bar{p}=\langle p_n:n<\omega\rangle\in (\qscF)^{\textstyle\omega}$. 
Then\\
$\bar{p}$ is pre-dense in $\qscF$\quad if and only if\\
for each $w\in\bigcup\limits_{m<\omega}\prod\limits_{i<m}\bH(i)$ and
$f\in\cF$ the tree $\cT^{\bar{p}}_{w,f}$ is well--founded.
\end{claim}

\begin{proof}[Proof of the claim] 
Suppose that, for some $w$ and $f$, the tree $\cT^{\bar{p}}_{w,f}$ has an
$\omega$--branch and let $q=(w,t_0,t_1,\ldots)$ be such a
branch. Necessarily $q\in\qscF$ (as witnessed by $f$). If follows from the
definition of $\cT^{\bar{p}}_{w,f}$ that $\POS(q)\cap\POS(p_n)$ is finite
for each $n\in\omega$ and therefore $q\;\bot_{\qscF}\; p_n$ (remember
\ref{sidePOS}). 

Now assume that $\bar{p}$ is not pre-dense in $\qscF$ and let $q\in\qscF$ be
a condition incompatible with all $p_n$. We may demand that for some
$f\in\cF$ we have $(\forall i\in\omega)(\nor[t^q_i]\geq f(m^{t^q_i}_\dn))$. 
It should be clear that $q$ determines an $\omega$--branch in the tree
$\cT^{\bar{p}}_{w^q,f}$ (remember that $(K,\Sigma)$ is $h$--linked).
\end{proof}

Similarly we deal with the respective variant of 2(b).
\end{proof}

\subsection{Unbounded and dominating reals}

\begin{lemma}
\label{lemdom}
Let $\bV\subseteq\bV^*$ be universes of ${\rm ZFC}^*$. Assume that $\langle
f_i:i<\omega\rangle\in\bV$ and $g\in\bV^*$ are such that
\begin{enumerate}
\item[(a)] $g\in\baire$, $f_i\in\baire$, $f_{i+1}<^*f_i$ for all
$i\in\omega$,
\item[(b)] $(\forall i\in\omega)(\exists^\infty k\in\omega)(g(k)<f_i(k))$,
\item[(c)] if $h\in\baire\cap\bV$ is such that $(\forall i\in\omega)(h<^*
f_i)$, then $h<^*g$.
\end{enumerate}
Then $\baire\cap\bV$ is bounded in $\bV^*$.
\end{lemma}

\begin{proof}
It follows from the assumptions (a), (b) that we may find an infinite set
$K=\{k_0,k_1,k_2,\ldots\}\in\bV^*\cap\iso$ such that for each $i\in\omega$
we have
\begin{enumerate}
\item[$(*)$] $f_0(k_i)>f_1(k_i)>\ldots>f_i(k_i)>g(k_i)$.
\end{enumerate}
Let $\varphi\in\baire\cap\bV^*$ be such that $(\forall n\in\omega)(|K\cap
(n,\varphi(n))|>2^{n+1})$. We claim that the function $\varphi$ dominates
$\baire\cap\bV$, i.e.
\[(\forall f\in\baire\cap\bV)(\forall^\infty n\in\omega)(f(n)<\varphi(n)).\]
If not, then we may choose an increasing sequence $\langle n_i:i<\omega
\rangle\in\bV$ of integers such that $n_0=0$ and
\begin{enumerate}
\item[(i)]  $(\exists^\infty i\in\omega)(|(n_i,n_{i+1})\cap K|>2^{n_i})$,
\item[(ii)] $(\forall i\in\omega)(\forall n\geq n_{i+1})(f_{i+1}(n)<
f_i(n))$. 
\end{enumerate}
Define $h\in\baire\cap\bV$ by $h\rest [n_i,n_{i+1})=f_i\rest [n_i,n_{i+1})$
(for $i\in\omega$). It follows from (ii) that $h<^*f_i$ for each
$i\in\omega$, so we may apply the assumption (c) to conclude that $h<^*g$. 
But look at the clauses $(*)$ and (i) above. Whenever $|(n_i,n_{i+1})\cap K|
>2^{n_i}$, there is $\ell\in\omega$ such that $k_\ell\in (n_i,n_{i+1})$,
$\ell>i$ and 
\[f_0(k_\ell)>\ldots>f_i(k_\ell)=h(k_\ell)>\ldots>f_\ell(k_\ell)>g(k_\ell),\]
so easily we get a contradiction.
\end{proof}

\begin{definition}
\label{reducible}
Let $(K,\Sigma)$ be a creating pair or a (local) tree--creating pair. 
\begin{enumerate}
\item (See \cite[Def.\ 5.1.6]{RoSh:470}) We say that $(K,\Sigma)$ is {\em
reducible\/} if for each $t\in K$ with $\nor[t]\geq 3$, there is $s\in\Sigma
(t)$ such that $\frac{\nor[t]}{2}\leq\nor[s]\leq\nor[t]-1$.
\item The pair $(K,\Sigma)$ is {\em normal\/} if it is reducible, linked and
if 
\begin{enumerate}
\item[$(\boxplus)$] for each $s,t\in K$:
\[\nor[s]<\nor[t]\quad\Rightarrow\quad (\forall u\in\dom(\val[t]))(\exists
v)(\langle u,v\rangle\in\val[t]\setminus\val[s]).\]
\end{enumerate}
\item A creating pair $(K,\Sigma)$ is {\em semi--normal\/} if it is linked, 
and for each $n\in\omega$ and $t\in K$ such that $\nor[t]>2^{n+2}$,
$m^t_\up-m^t_\dn>2^{2^{n+4}}$, there is a sequence $\langle s_\ell:\ell\leq
n\rangle\subseteq K$ satisfying
\begin{enumerate}
\item[$(\alpha)$] $s_0=t$, $s_{\ell+1}\in\Sigma(s_\ell)$, $2^{n+1-\ell}<
\nor[s_{\ell+1}]\leq 2^{n-\ell+2}$ (for $\ell<n$), and
\item[$(\beta)$]  if $s\in K$, $m^s_\dn=m^t_\dn$, $m^s_\up=m^t_\up$ and
$\nor[s]>2^{n-\ell+3}$, $\ell<n$, then
\[(\forall u\in\dom(\val[s]))(\exists v)(\langle u,v\rangle\in\val[s]
\setminus \val[s_{\ell+1}]).\]
\end{enumerate}
\item A $\otimes$--creating triple $(K,\Sigma,\Sigma^\bot)$ (or just
$(K,\Sigma)$) is {\em super--gluing\/} if it is gluing and for every $s_0,
\ldots, s_k\in K$ and $N\in\omega$ such that $m^{s_i}_{\up}\leq m^{s_{i+
1}}_{\dn}$ for $i<k$ and $m^{s_k}_\up\leq N$, there is $s\in K$ satisfying:
\begin{itemize}
\item $m^s_\dn=m^{s_0}_\dn$, $m^s_\up=N$, $\dom(\val[s])=\dom(\val[s_0])$, and
\item $\nor[s]\geq\min\{\nor[s_i]:i\leq k\}-1$, and
\item $(\forall\langle u,v\rangle\in\val[s])(\forall i\leq k)(\langle v\rest
m^{s_i}_\dn,v\rest m^{s_i}_\up\rangle\in\val[s_i])$.
\end{itemize}
\end{enumerate}
\end{definition}

\begin{remark}
Note that ``normal'' implies ``semi--normal''. What we really need in the
proofs of \ref{adddom}(1,2) is semi--normality (or rather a suitable variant
of it). However, the normality is more natural and only in the case of
$\otimes$--creating triples (which are gluing and have the cutting property)
the natural norms are semi--normal but not normal; see {\bf Examples} at the
end of this section.
\end{remark}

\begin{theorem}
\label{adddom}
\begin{enumerate}
\item Let $f:\omega\times\omega\longrightarrow\omega$ be a fast function,
and let $(K,\Sigma)$ be a local creating pair (a local tree--creating pair,
respectively). Assume that $(K,\Sigma)$ is normal and $\qonefx\neq\emptyset$
($\bQ^\tree_f(K,\Sigma)\neq\emptyset$, resp.). Then the forcing notion
$\qonefx$ ($\bQ^\tree_f(K,\Sigma)$, resp.) adds a dominating real.
\item If $(K,\Sigma)$ is a normal (local) tree--creating pair and
$\bQ^\tree_1(K,\Sigma)\neq\emptyset$, then the forcing notion
$\bQ^\tree_1(K,\Sigma)\neq\emptyset$ adds a dominating real.
\item Assume that $(K,\Sigma,\Sigma^\bot)$ is a semi--normal
$\otimes$--creating triple which is super--gluing and has the cutting
property (and $\qinf\neq \emptyset$). Then the forcing notion $\qinf$ adds a
dominating real. 
\end{enumerate}
\end{theorem}

\begin{proof}
In all cases we will use Lemma \ref{lemdom} for functions $f_i\in
\prod\limits_{n\in\omega}(n+1)$ defined by $f_i\rest [0,i)\equiv 0$,
$f_i(n)=n-i$ for $n\geq i$ (for $i\in\omega$) and a suitably chosen name
$\dot{g}$ for a function in $\prod\limits_{n\in\omega}(n+1)$.

\noindent (1)\qquad We consider the case when $(K,\Sigma)$ is a local
creating pair only. 

Let $p\in\qonefx$. Using the normality of $(K,\Sigma)$, choose an increasing
sequence $\langle m_n:n<\omega\rangle\subseteq\omega$ and a sequence
$\langle s^\ell_n:\ell\leq n,\ n<\omega\rangle\subseteq K$ such that for
each $n\in\omega$ and $\ell<n$: 
\begin{enumerate}
\item[(a)] $s^0_n=t^p_{m_n}$, $\nor[t^p_{m_n}]>f(n+2,m^{t^p_{m_n}}_\dn)$,
\item[(b)] $s^{\ell+1}_n\in\Sigma(s^\ell_n)$, $f(n-\ell+1,m^{s^{\ell+1}_n
}_\dn)<\nor[s^{\ell+1}_n]\leq f(n-\ell+2,m^{s^{\ell+1}_n}_\dn)$.
\end{enumerate}
Let $\dot{W}$ be the name for $\qonefx$--generic real, i.e.
\[\forces_{\qonefx}\dot{W}=\bigcup\{w^q:q\in\Gamma_{\qonefx}\}\]
(see \cite[Def.\ 1.1.13, Prop.\ 1.1.14]{RoSh:470}). Let $\dot{g}$ be a
$\qonefx$--name for a function in $\prod\limits_{n\in\omega}(n+1)$ defined
by 
\[p\forces\mbox{`` }(\forall n\in\omega)(\forall \ell\leq n)(\dot{g}(n)=
\ell\ \Leftrightarrow\ \langle\dot{W}\rest m^{s^\ell_n}_\dn,\dot{W}\rest
m^{s^\ell_n}_\up\rangle\in\val[s^\ell_n]\setminus\val[s^{\ell+1}_n])
\mbox{ ''}\]
(if $\langle\dot{W}\rest m^{s^n_n}_\dn,\dot{W}\rest m^{s^n_n}_\up\rangle
\in\val[s^n_n]$ then $\dot{g}(n)=n$).
\begin{claim}
\label{cl3}
\begin{enumerate}
\item[$(\alpha)$] $p\forces_{\qonefx}(\forall i\in\omega)(\exists^\infty n
\in\omega)(\dot{g}(n)<f_i(n))$,
\item[$(\beta)$]  Assume that $h\in\baire$ is such that $h<^* f_i$ for all
$i\in\omega$. Then $p\forces$`` $h<^*\dot{g}$ ''.
\end{enumerate}
\end{claim}

\begin{proof}[Proof of the claim] $(\alpha)$\quad Let $i\in\omega$ and let
$p_0\geq p$, $N\in\omega$. Take $n>N+i$ so that for some $k<\omega$ we have
$t^{p_0}_k\in\Sigma(t^p_{m_n})$ and $\nor[t^{p_0}_k]>f(i+3,
m^{t^{p_0}_k}_\dn)$. Since $(K,\Sigma)$ is normal, we may find $w\in\pos(
w^{p_0},t^{p_0}_0,\ldots,t^{p_0}_k)$ such that $\langle w\rest
m^{t^{p_0}_k}_\dn, w\rest m^{t^{p_0}_k}_\up\rangle\notin \val[s^{n-i}_n]$
(remember $\nor[s^{n-i}_n]\leq f(i+3,m^{s^{n-i}_n}_\dn)$). Clearly the
condition $q=(w,t^{p_0}_{k+1},t^{p_0}_{k+2},\ldots)$ forces that
$\dot{g}(n)<n-i=f_i(n)$. 
\smallskip

\noindent $(\beta)$\quad Let $h\in\baire$ be such that $h<^* f_i$ for all
$i\in\omega$ and let $p_0\geq p$. Let $\ell$ be such that $t^{p_0}_0\in
\Sigma(t^p_\ell)$ (so $t^{p_0}_k\in\Sigma(t^p_{k+\ell})$ for each $k$). We
may assume that if $m_n\geq\ell$ then $h(n)<n-5$ and that $(\forall
k<\omega)(\nor[t^{p_0}_k]> f(5,m^{t^{p_0}_k}_\dn))$. For each $k\in\omega$
choose $t_k$ as follows:
\begin{itemize}
\item if $k+\ell\notin\{m_n:n\in\omega\}$, then $t_k=t^{p_0}_k$,
\item if $k+\ell=m_n$, $n\in\omega$, then $t_k\in\Sigma(t^{p_0}_k)\cap
\Sigma(s^{h(n)+1}_n)$ is such that $\nor[t_k]\geq\min\{\nor[t^{p_0}_k],
\nor[s^{h(n)+1}_n]\}-1$
\end{itemize}
(remember that $(K,\Sigma)$ is linked). Since $\nor[s^{h(n)+1}_n]>f(n-h(n)+ 
1, m^{s^0_n}_\dn)$ we easily see that $q=(w^{p_0},t_0,t_1,\ldots)\in
\qonefx$, and clearly it is a condition stronger than $p_0$. As $q\forces_{
\qonefx} (\forall n\in\omega)(m_n\geq\ell\ \Rightarrow\ \dot{g}(n) >h(n))$,
the claim follows. 
\end{proof}
Now, the first clause of the theorem is an immediate consequence of
\ref{cl3} and \ref{lemdom}. 
\medskip

\noindent (2)\qquad The proof is similar to the one above. Let
$p\in\bQ^\tree_1(K,\Sigma)$. Choose fronts $F_n$ of $T^p$ (for $n\in\omega$)
such that for each $n$:
\begin{itemize}
\item $(\forall\eta\in F_{n+1})(\exists\nu\in F_n)(\nu\vartriangleleft
\eta)$,
\item $(\forall\eta\in F_n)(\nor[t^p_\eta]>2^{n+2})$
\end{itemize}
(clearly possible; see \cite[Prop.\ 1.3.8]{RoSh:470}). For each $n\in
\omega$ and $\eta\in F_n$ choose a sequence $\langle s^\ell_\eta:\ell\leq
n\rangle\subseteq K$ such that
\[s^0_\eta=t^p_\eta,\quad s^{\ell+1}_\eta\in\Sigma(s^\ell_\eta)\quad \mbox{
and }\quad 2^{n-\ell+1}<\nor[s^{\ell+1}_\eta]\leq 2^{n-\ell+2}.\]
Let $\dot{W}$ be the name for the $\bQ^\tree_1(K,\Sigma)$--generic real and
let $\dot{g}$ be a $\bQ^\tree_1(K,\Sigma)$--name for a real in
$\prod\limits_{n\in\omega}(n+1)$ such that (the condition $p$ forces that)
if $\eta=\dot{W}\rest m\in F_n$ (for some $m,n\in\omega$) and $\dot{W}\rest
(m+1)\in\pos(s^\ell_\eta)\setminus\pos(s^{\ell+1}_\eta)$, then $\dot{g}(n)= 
\ell$. 
\begin{claim}
\label{cl4}  
\begin{enumerate}
\item[$(\alpha)$] $p\forces_{\bQ^\tree_1(K,\Sigma)}(\forall i\in\omega)(
\exists^\infty n\in\omega)(\dot{g}(n)<f_i(n))$,
\item[$(\beta)$]  If $h\in\baire$ is such that $h<^* f_i$ for all
$i\in\omega$, then $p\forces_{\bQ^\tree_1(K,\Sigma)}$`` $h<^*\dot{g}$ ''.
\end{enumerate}
\end{claim}

\begin{proof}[Proof of the claim] $(\alpha)$\quad Like \ref{cl3}$(\alpha)$.
\smallskip

\noindent $(\beta)$\quad Let $q\geq p$. We may assume that for some $m>2$ we
have: $\mrot(q)\in F_m$, $\nor[t^q_\nu]>8$ for all $\nu\in T^q$, and $h(n)<
n-5$ for all $n\geq m$. We build inductively a tree $T\subseteq T^q$ and a
system $\langle t_\eta:\eta\in T\rangle$ as follows. We declare that
$\mrot(q)=\mrot(T)\in T$. Suppose we have declared that $\eta\in T$. If
$\eta\notin\bigcup\limits_{n\in\omega}F_n$, then we let $t_\eta=t^q_\eta$
and we declare $\pos(t_\eta)\subseteq T$. If $\eta\in F_n$ for some $n\geq
m$, then we choose $t_\eta\in\Sigma(t^q_\eta)\cap \Sigma(s^{h(n)+1}_\eta)$
such that $\nor[t_\eta]\geq\min\{\nor[t^q_\eta],2^{n-h(n)+1}\}-1$, and we
declare $\pos(t_\eta)\subseteq T$.\\
Finally, we let $q^*=\langle t_\eta:\eta\in T\rangle$ and we notice that
$q^*\forces h<^*\dot{g}$.
\end{proof}
\medskip

\noindent (3)\qquad Let $p\in\qinf$. By ``gluing'', we may assume that
$m^{t^p_k}_\up-m^{t^p_k}_\dn> 2^{2^{k+4}}$ for each $k<\omega$. Using
semi--normality we may choose an increasing sequence $\langle m_n:n<\omega
\rangle$ and a sequence $\langle s^\ell_n: \ell\leq n,\ n<\omega\rangle
\subseteq K$ such that 
\begin{enumerate}
\item[(a)] $s^0_n=t^p_{m_n}$, $\nor[t^p_{m_n}]>2^{n+2}$, $s^{\ell+1}_n\in
\Sigma(s^\ell_n)$, $2^{n-\ell+1}<\nor[s^{\ell+1}_n]\leq 2^{n-\ell+2}$ (for
$\ell<n<\omega$), 
\item[(b)] if $t\in K$, $m^t_\dn=m^{s^0_n}_\dn$, $m^t_\up=m^{s^0_n}_\up$,
and $\nor[t]>2^{n-\ell+3}$, $\ell<n<\omega$, then
\[(\forall u\in\dom(\val[t]))(\exists v)(\langle u,v\rangle\in\val[t]
\setminus\val[s_n^{\ell+1}]).\]
\end{enumerate}
Now we define the name $\dot{g}$ like before, so
\[p\forces_{\qinf}\mbox{`` }\dot{g}(n)=\ell\ \Leftrightarrow\ \langle\dot{W}
\rest m^{s^\ell_n}_\dn,\dot{W}\rest m^{s^\ell_n}_\up\rangle\in\val[s^\ell_n]
\setminus\val[s^{\ell+1}_n]\mbox{ ''}.\]
\begin{claim}
\label{cl5}  
\begin{enumerate}
\item[$(\alpha)$] $p\forces_{\qinf}(\forall i\in\omega)(\exists^\infty n\in
\omega)(\dot{g}(n)<f_i(n))$,
\item[$(\beta)$]  If $h<^* f_i$ for all $i\in\omega$, then
$p\forces_{\qinf}$`` $h<^*\dot{g}$ ''. 
\end{enumerate}
\end{claim}

\begin{proof}[Proof of the claim] $(\alpha)$\quad Suppose $q\geq p$, and 
$i,N<\omega$. Passing to a stronger condition (using ``gluing and cutting'')
we may assume that
\begin{itemize}
\item $\nor[t^q_0]>2^{i+4}$ and
\item $m^{t^q_0}_\dn=m^{t^p_{m_n}}_\dn$, $m^{t^q_0}_\up=m^{t^p_{m_n}}_\up$
for some $n>N+i+1$.
\end{itemize}
Choose $w\in\pos(w^q,t^q_0)$ such that $\langle w^q,w\rangle\notin s^{n-i}_n$
and look at the condition $q'=(w,t^q_1,t^q_2,\ldots)$.

\noindent $(\beta)$\quad Let $q\geq p$. Passing to a stronger condition if
necessary, we may assume that for some increasing sequence $\langle N_k:k<
\omega\rangle\subseteq\omega$ we have:
\begin{itemize}
\item $m^{t^q_k}_\dn=m^{t^p_{m_{N_k}}}_\dn$, $m^{t^q_k}_\up=m^{t^p_{
m_{N_{k+1}}}}_\dn$, $\nor[t^q_k]>5$ for all $k<\omega$,
\item if $n\geq N_0$ then $h(n)<n-5$.
\end{itemize}
Using ``super--gluing'' choose creatures $s_k\in K$ (for $k\in\omega$) such
that 
\begin{itemize}
\item $m^{s_k}_\dn=m^{t^p_{m_{N_k}}}_\dn$,\quad $m^{s_k}_\up=m^{t^p_{m_{
N_{k+1}}}}_\dn$,\quad and 
\item $\nor[s_k]\geq\min\{\nor[s^{h(n)+1}_n]: N_k\leq n<N_{k+1}\}-1$,\quad
and
\item $(\forall\langle u,v\rangle\in\val[s_k])(\forall n\in [N_k,N_{k+1}))( 
\langle v\rest m^{s^0_n}_\dn,v\rest m^{s^0_n}_\up\rangle \in
\val[s^{h(n)+1}_n])$.
\end{itemize}
Apply ``linked'' to choose creatures $t_k\in\Sigma(s_k)\cap\Sigma(t^q_k)$
such that 
\[1+\nor[t_k]\geq\min\{\nor[s_k],\nor[t^q_k]\}\geq \min\{\nor[t^q_k],2^{n-
h(n)}: N_k\leq n<N_{k+1}\}.\]
Then $q^*=(w^q,s_0,s_1,s_2,\ldots)$ is a condition in $\qinf$, $q^*\geq q$
and it forces that $(\forall n\geq N_0)(h(n)<\dot{g}(n))$.
\end{proof}
\end{proof}

\begin{remark}
Note that \ref{adddom}(1) applies to forcing notions $\bQ^*_\infty(K,
\Sigma)$ too, see \ref{rem1}.
\end{remark}

\begin{definition}
\label{nicefor}
A ccc forcing notion $\bP$ is {\em nice\/} if there is a partition $\langle
P_m:m<\omega\rangle$ of $\bP$ such that
\begin{enumerate}
\item[$(\bigstar)$] if $\langle p_n:n<\omega\rangle\subseteq\bP$ is a
maximal antichain in $\bP$, $m\in\omega$,\\
then there is $N<\omega$ such that 
\[(\forall p\in P_m)(\exists n<N)(p,p_n\mbox{ are compatible }).\]
\end{enumerate}
\end{definition}

\begin{theorem}
[Miller \cite{Mi81}, Brendle and Judah {\cite{BrJd93}}; see {\cite[Thm.\
6.5.11]{BaJu95}}] 
\label{presunb}
If $\cF\subseteq\baire$ is an unbounded family and $\bP$ is a nice ccc
forcing notion, then 
\[\forces_{\bP}\mbox{`` the family $\cF$ is unbounded ''.}\]
\end{theorem}

\begin{remark}
Since no dominating reals can be added at limit stages of FS iterations of
ccc forcing notions (see \cite[Con.\ VI.3.17]{Sh:f}), it follows from
\ref{presunb} that FS iterations of nice forcing notions do not add
dominating reals.
\end{remark}

\begin{definition}
\label{Cohen}
\begin{enumerate}
\item Let $(K,\Sigma)$ be a local creating pair (or a local tree--creating
pair) for $\bH$. We say that $(K,\Sigma)$ is {\em Cohen--producing\/} if for
each $n\in\omega$ there is a set $A_n\subseteq\bH(n)$ such that
\begin{quotation}
if $t\in K$, $\nor[t]>1$, $u\in\dom(\val[t])$, $\lh(u)=n$,

\noindent then there are $v_0,v_1$ such that $\langle u,v_0\rangle, \langle
u,v_1\rangle\in\val[t]$ and $v_1(n)\in A_n$ and $v_0(n)\notin A_n$.
\end{quotation}
\item A creating pair $(K,\Sigma)$ is {\em of the BCB--type\/} if it is
local, forgetful and satisfies the following condition:
\begin{enumerate}
\item[$(\circledast^{\rm BCB})$]\quad for every sequence $\langle s_n:n<
\omega\rangle$ of creatures from $K$ with $m^{s_n}_\dn=i$, $\nor[s_n]\geq 2$
(for all $n$), there are $a_0,\ldots,a_m\in\bH(i)$ and an increasing
sequence $\langle n_k:k<\omega\rangle\subseteq\omega$ such that 
\[(\forall a\in\bH(i)\setminus\{a_0,\ldots,a_m\})(\forall^\infty k\in\omega)
(\forall u\in\dom(\val[s_{n_k}]))(\langle u,u\conc\langle a\rangle\rangle\in 
\val[s_{n_k}]).\]
\end{enumerate}
\item A local tree--creating pair $(K,\Sigma)$ is {\em of the
${\rm BCB}^\tree$--type\/} if 
\begin{enumerate}
\item[$(\circledast^{\rm BCB}_\tree)$]\quad for every $\eta\in
\bigcup\limits_{n<\omega}\prod\limits_{i<n}\bH(i)$ and a sequence $\langle
s_n: n<\omega\rangle\subseteq\TCR_\eta\cap K$ such that $\nor[s_n]\geq 2$, 
there are $a_0,\ldots,a_m\in \bH(\lh(\eta))$ and an increasing sequence
$\langle n_k:k<\omega\rangle\subseteq\omega$ such that
\[(\forall a\in\bH(\lh(\eta))\setminus\{a_0,\ldots,a_m\})(\forall^\infty
k\in\omega)(\eta\conc\langle a\rangle\in\pos(s_{n_k})).\]
\end{enumerate}
\end{enumerate}
\end{definition}

\begin{remark}
\begin{enumerate}
\item Note that if $\bH(i)$ is finite for each $i\in\omega$ then any local
forgetful creating pair (local tree creating pair) is of the BCB--type
(BCB$^\tree$--type, respectively). 
\item The difference between BCB and BCB$^\tree$ is not serious, the two
notions are just fitted to their contexts.
\item ``BCB'' stands for ``bounded -- co-bounded''. The ``bounded'' part
reflects what is stated in (1) above, and the ``co-bounded'' is supposed to
point out the analogy to the co-bounded topology on $\omega$ in the case
when each $\bH(i)$ is infinite; compare Miller \cite{Mi81}.
\end{enumerate}
\end{remark}

\begin{theorem}
\label{unbounded}
Assume that $h:\omega\times\omega\longrightarrow\omega$ is a regressive
function and $\cF\subseteq (\omega\setminus 2)^{\textstyle\omega}$ is a
countable $h$--closed family which is $\geq^*$--directed. 
\begin{enumerate}
\item If $(K,\Sigma)$ is a local Cohen--producing $h$--linked creating pair
(tree--creating pair, respectively), then the forcing notion $\bQ^*_{\cF}(K,  
\Sigma)$ ($\bQ^\tree_{\cF}(K,\Sigma)$, resp.) adds a Cohen real.
\item If $(K,\Sigma)$ is an $h$--linked tree--creating pair of the {\rm
BCB$^\tree$}--type, then the forcing notion $\bQ^\tree_{\cF}(K,\Sigma)$ is
nice. 
\item If a creating pair $(K,\Sigma)$ is $h$--linked and of the {\rm
BCB}--type, then the forcing notion $\bQ^*_{\cF}(K,\Sigma)$ is nice.
\end{enumerate}
\end{theorem}

\begin{proof} 
(1)\qquad Let $(K,\Sigma)$ be a creating pair for $\bH$ and let sets $A_n
\subseteq\bH(n)$ (for $n\in\omega$) witness that it is  Cohen--producing.
Let $\dot{c}$ be a $\bQ^*_{\cF}(K,\Sigma)$--name for a real in $\can$
defined by  
\[\forces_{\bQ^*_{\cF}(K,\Sigma)}(\forall n\in\omega)(\dot{c}(n)=1\
\Leftrightarrow\ \dot{W}(n)\in A_n)\]
(where $\dot{W}$ is the name for the $\bQ^*_{\cF}(K,\Sigma)$--generic
real). Suppose that a condition $p\in\bQ^*_{\cF}(K,\Sigma)$ is such that
$\nor[t^p_n]>1$ for all $n\in\omega$. Let $\sigma:[\lh(w^p),N]
\longrightarrow 2$, $\lh(w^p)\leq N<\omega$. It should be clear that there
is $w\in\pos(w^p,t^p_0,\ldots,t^p_{N-\lh(w^p)})$ such that 
\[(\forall n\in [\lh(w^p),N])(w(n)\in A_n\ \Leftrightarrow\ \sigma(n)=1).\]
Hence easily $\dot{c}$ is a name for a Cohen real. 
\medskip

\noindent(2)\qquad Let $\nu\in\bigcup\limits_{n\in\omega}\prod\limits_{i<n}
\bH(i)$, $f\in\cF$ and let
\[P_{\nu,f}\stackrel{\rm def}{=}\{p\in\bQ^\tree_{\cF}(K,\Sigma): \mrot(p)=
\nu\ \&\ (\forall\eta\in T^p)(\nor[t^p_\eta]\geq f(\lh(\eta)))\}.\]
Suppose that $\langle p_n:n<\omega\rangle\subseteq\bQ^\tree_{\cF}(K,\Sigma)$
is a maximal antichain such that, for each $n\in\omega$ and $\eta\in
T^{p_n}$, we have $\nor[t^{p_n}_\eta]\geq 2$. 
\begin{claim}
\label{cl6}
There is $N<\omega$ such that 
\[(\forall q\in P_{\nu,f})(\exists n<N)(q,p_n\mbox{ are compatible }).\]
\end{claim}
\begin{proof}[Proof of the claim]
Assume not. Then we may choose a sequence $\langle q_k:k<\omega\rangle
\subseteq P_{\nu,f}$ such that for each $n<k<\omega$ the conditions $q_k$
and $p_n$ are incompatible. We inductively build a tree $T$ and a system
$\langle s_\eta:\eta\in T\rangle$ together with a sequence $\langle X_n,
Y_n: n<\omega\rangle$ so that 
\begin{enumerate}
\item[$(\alpha)$] $X_n\subseteq\bigcup\limits_{m\in\omega}\prod\limits_{i
<m}\bH(i)$, $Y_{n+1}\subseteq Y_n\in\iso$,
\item[$(\beta)$]  $(\forall\eta\in X_n)(\forall^\infty k\in Y_n)(\eta\in
T^{q_k})$,
\item[$(\gamma)$] $\nor[s_\eta]\geq f(\lh(\eta))$, $T=\bigcup\limits_{n\in
\omega} X_n$.
\end{enumerate}
Fix a bijection $\#:\bigcup\limits_{m\in\omega}\prod\limits_{i<m}\bH(i)
\longrightarrow\omega$ such that $\eta_0\vartriangleleft\eta_1$ implies
$\#(\eta_0)<\#(\eta_1)$. 

We declare that $\nu=\mrot(T)$, $X_0=\{\nu\}$, $Y_0=\omega$.

Suppose we have arrived to the $(n+1)^{\rm th}$ stage of the construction
and $X_n,Y_n$ have been already defined so that the clauses
$(\alpha),(\beta)$ above are satisfied. Let $\eta\in X_n$ be such that 
\[\#(\eta)=\min\{\#(\eta'):\eta'\in X_n\}.\]
Let $Y'_n\in\iso$ consist of these $k\in Y_n$ that $\eta\in T^{q_k}$
(remember $(\beta)$). Apply $(\circledast^{\rm BCB}_\tree)$ of
\ref{Cohen}(3) to the sequence $\langle t^{q_k}_\eta:k\in Y'_n\rangle$ to
choose an infinite set $Y_{n+1}\subseteq Y'_n$ such that, letting $k^*=
\min(Y_{n+1})$ and $s_\eta=t^{q_{k^*}}_\eta$ we have
\[(\forall\eta'\in\pos(s_\eta))(\forall^\infty k\in Y_{n+1})(\eta'\in\pos(
t^{q_k}_\eta)).\]
Finally, we let $X_{n+1}=(X_n\setminus\{\eta\})\cup\pos(s_\eta)$. This
finishes the description of the inductive step.

After the construction is carried out we let $q^*=\langle s_\eta:\eta\in
T\rangle$. It should be clear that $q^*\in\bQ^\tree_{\cF}(K,\Sigma)$ (and
even $q^*\in P_{\nu,f}$). Consequently we find $n<\omega$ such that the
conditions $p_n$ and $q^*$ are compatible. 

Suppose that $\nu\trianglelefteq\mrot(p_n)$. Then necessarily $\mrot(p_n)
\in T$. It follows from our construction (remember clause $(\beta)$) that we 
may find $k>n$ such that $\mrot(p_n)\in T^{q_k}$. But then, using the
assumption that $(K,\Sigma)$ is $h$--linked and $\cF$ is $h$--closed and 
$\geq^*$--directed, we immediately get that the conditions $p_n,q_k$ are
compatible, contradicting the choice of the $q_k$. Similarly, if $\mrot(p_n)
\vartriangleleft\nu$ then taking any $k>n$ we get that the conditions
$q_k,p_n$ are compatible, again a contradiction. 
\end{proof}

Since the conditions of the form used above are dense in $\bQ^\tree_{\cF}(
K,\Sigma)$ one easily concludes that the forcing notion $\bQ^\tree_{\cF}(
K,\Sigma)$ is nice.
\medskip

\noindent (3)\qquad Similarly.
\end{proof}

\subsection{Examples}
Our first example recalls the forcing notion of \cite[\S 3]{RoSh:628}. Let
us start with presenting the main tool for this type of constructions -- 
norms determined by Hall's Marriage Theorem.

\begin{definition}
\label{hallnorm}
Let $\bH:\omega\longrightarrow\cHa$.
\begin{enumerate}
\item Let $\cK^\bH$ be the collection of all finite non-empty families
$\Delta$ of finite partial functions $f$ such that $\emptyset\neq\dom(f)
\subseteq\omega$ and $f(n)\in\bH(n)$ for all $n\in\dom(f)$. For integers
$m_0<m_1$ let  
\[\cK^\bH_{m_0,m_1}\stackrel{\rm def}{=}\{\Delta\in\cK^\bH: (\forall
f\in\Delta)(\dom(f)\subseteq [m_0,m_1))\}.\]
\item Let $\Delta\in\cK^\bH$, $k\in\omega$. A function $F:\Delta
\longrightarrow[\omega]^{\textstyle k}$ is {\em a $k$--selector for
$\Delta$} if  
\[(\forall f,f'\in\Delta)\big(F(f)\subseteq\dom(f)\ \mbox{ and }\ 
f\neq f'\ \Rightarrow F(f)\cap F(f')=\emptyset\big).\]
\item For $\Delta_0,\Delta_1\in\cK^\bH$ we write $\Delta_0\preceq\Delta_1$
whenever
\[(\forall f\in\Delta_0)(\exists g\in\Delta_1)(g\subseteq f).\]
\item We define {the Hall norms} of a set $\Delta\in\cK^\bH$ as follows:
\[\begin{array}{l}
\hn^+(\Delta)\stackrel{\rm def}{=}\max\{k+1:k\in\omega\mbox{ and there is
an $k$--selector for }\Delta\},\\
\hn(\Delta)\stackrel{\rm def}{=}\max\{k+1:k\in\omega\mbox{ and for every }
\Delta'\subseteq\Delta\mbox{ there is }\Delta''\subseteq\Delta'\mbox{ such
that}\\ 
\qquad\qquad\qquad\qquad\qquad\mbox{elements of $\Delta''$ have pairwise
disjoint domains and}\\
\qquad\qquad\qquad\qquad\qquad|\bigcup\limits_{f\in\Delta''}\dom(f)|\geq
k\cdot |\Delta'|\},\\ 
\HN(\Delta)\stackrel{\rm def}{=}\max\{\hn(\Delta'):\Delta\preceq
\Delta'\}. 
  \end{array}\]
\end{enumerate}
\end{definition}

\begin{lemma}
\label{basichn}
\begin{enumerate}
\item If $\Delta\in\cK^\bH$ and $k_0\in\omega$ then 
\[\hn^+(\Delta)>k_0\ \mbox{ if and only if }\ (\forall\Delta'\subseteq
\Delta)(|\bigcup\{\dom(f): f\in\Delta'\}|\geq k_0\cdot |\Delta'|)\]
and $1\leq\hn(\Delta)\leq\hn^+(\Delta)\leq\HN(\Delta)$.
\item If $\Delta_0,\Delta_1\in\cK^\bH$ and ${\bf xx}\in\{\hn,\hn^+,\HN\}$
then  
\[{\bf xx}(\Delta_0\cup\Delta_1)\geq \lfloor\min\{
\frac{{\bf xx}(\Delta_0)}{2},\frac{{\bf xx}(\Delta_1)}{2}\}\rfloor.\]
\item If $m^0_0<m^0_1\leq m^1_0<m^1_1\leq\ldots\leq m^k_0<m^k_1<\omega$,
$\Delta_i\in\cK^\bH_{m^i_0,m^i_1}$ (for $i\leq k$) and ${\bf
xx}\in\{\hn,\hn^+,\HN\}$ then  
\[{\bf xx}(\bigcup\limits_{i\leq k}\Delta_i)=\min\{{\bf xx}(\Delta_i):i\leq
k\}.\]
\item Suppose that $m_0<m<m_1<\omega$ and $\Delta\in\cK^\bH_{m_0,m_1}$. Let 
\[\begin{array}{lcl}
\Delta_0&=&\{f\rest [m_0,m): f\in\Delta\ \&\ |\dom(f)\cap [m_0,m)|\geq
\frac{1}{2}|\dom(f)|\},\\
\Delta_1&=&\{f\rest [m,m_1): f\in\Delta\ \&\ |\dom(f)\cap [m,m_1)|\geq
\frac{1}{2}|\dom(f)|\}.
  \end{array}\]
Then, for $i<2$, either $\Delta_i=\emptyset$ or $\hn(\Delta_i)\geq
\frac{1}{2}\hn(\Delta)$. 
\end{enumerate}
\end{lemma}

\begin{proof}
(1)\quad It follows from Hall's Theorem (see Hall \cite{Ha35}) and the
definitions of the norms. 

\noindent (2)--(4)\quad Straightforward (compare \cite[Claim
3.1.2]{RoSh:628}). 
\end{proof}

\begin{example}
\label{bas628}
Let $\bH:\omega\longrightarrow\cHa$, $|\bH(n)|\geq 2$ for all
$n\in\omega$.\\ 
We construct a $\otimes$--creating triple $(K_\bH,\Sigma_\bH,
\Sigma^\bot_\bH)$ for $\bH$ which: 
\begin{enumerate}
\item is semi--normal (see \ref{reducible}(3)), forgetful (see
\ref{varia}(4)) and super--gluing (see \ref{reducible}(4)),   
\item has the cutting property (see \ref{forlin}(3)),
\item is really finitary (see \ref{finitary}) provided $\bH(n)$ is finite
for each $n\in\omega$. 
\end{enumerate}
\end{example}

\begin{proof}[Construction] 
Let $K_\bH$ consist of all creatures $t\in\CR[\bH]$ such that 
\begin{itemize}
\item $\dis[t]=(m^t_\dn,m^t_\up,\Delta_t)$ for some $\Delta_t\in
\cK^\bH_{m_0,m_1}\cup\{\emptyset\}$ such that, if $\Delta_t\neq\emptyset$,
$\hn^+(\Delta_t)>1$, 
\item $\val[t]=\{\langle u,v\rangle\in\prod\limits_{i<m^t_\dn}\bH(i)\times
\prod\limits_{i<m^t_\up}\bH(i): u\vartriangleleft v\ \&\ (\forall f\in
\Delta_t)(f\nsubseteq v)\}$,
\item if $\Delta_t=\emptyset$ then $\nor[t]=m^t_\dn+1$, otherwise $\nor[t]=
\log_8(\hn(\Delta_t))$. 
\end{itemize}
(Note that $\hn^+(\Delta_t)>1$ implies $\val[t]\neq\emptyset$.) For $t_0,
\ldots,t_n\in K_\bH$ such that $m^{t_i}_\up=m^{t_{i+1}}_\dn$ (for $i<n$) let 
\[\Sigma_\bH(t_0,\ldots,t_n)=\{t\in K_\bH: m^t_\dn=m^{t_0}_\dn\ \&\ m^t_\up
= m^{t_n}_\up\ \&\ \bigcup_{i\leq n}\Delta_{t_i}\subseteq\Delta_t\}.\]
It should be clear that $\Sigma_\bH$ is a composition operation on $K_\bH$,
$K_\bH$ is countable and forgetful, and if each $\bH(n)$ is finite then
$K_\bH$ is really finitary.

For a creature $t\in K_\bH$ we define $\Sigma^\bot_\bH(t)$ as follows. It
consists of all finite sets $\{s_0,\ldots,s_n\}\subseteq K_\bH$ (a suitable 
enumeration) such that  
\begin{itemize}
\item $m^t_\dn=m^{s_0}_\dn<m^{s_0}_\up=m^{s_1}_\dn<\ldots<m^{s_{n-1}}_\up
=m^{s_n}_\dn<m^{s_n}_\up=m^t_\up$,\quad and 
\item $(\forall f\in\Delta_t)(\exists\ell\leq n)(f\rest [m^{s_\ell}_\dn,
m^{s_\ell}_\up)\in \Delta_{s_\ell})$.
\end{itemize}
It is clear that $\Sigma^\bot_\bH$ is a decomposition operation on $K_\bH$,
so $(K_\bH,\Sigma_\bH,\Sigma^\bot_\bH)$ is a $\otimes$--creating triple.

It follows from \ref{basichn}(2) that $(K_\bH,\Sigma_\bH)$ is linked, and
using \ref{basichn}(3) one easily shows that it is super--gluing. Similarly,
$(K_\bH,\Sigma_\bH,\Sigma^\bot_\bH)$ has the cutting property by
\ref{basichn}(4). 

Note that 
\begin{enumerate}
\item[$(*)$] if $f\in\prod\limits_{i=m_0}^{m_1-1}\bH(i)$, then $\hn(\{f\})= 
\HN(\{f\})=m_1-m_0+1$. 
\end{enumerate}
Hence, using \ref{basichn}(2), we may easily conclude that $(K_\bH,
\Sigma_\bH)$ is reducible. However, it is not normal -- one can build
$s,t\in K_\bH$ such that $\nor[s]<\nor[t]$ but $\val[t]\subseteq \val[s]$
(which is in some sense paradoxical, and this is why we modify this example
in \ref{bas628bis}).  

\begin{claim}
\label{cl7}
$(K_\bH,\Sigma_\bH,\Sigma^\bot_\bH)$ is semi--normal.
\end{claim}

\begin{proof}[Proof of the claim]
Let $n\in\omega$, $t\in K_\bH$ be such that $\nor[t]>2^{n+2}$, $m^t_\up-
m^t_\dn>2^{2^{n+4}}$. We may assume that $\Delta_t\neq\emptyset$ (remember
$(*)$). We choose inductively a sequence $\langle\Delta_\ell,A_\ell:\ell
\leq n\rangle$ such that  
\begin{enumerate}
\item $\Delta_\ell\in \cK^\bH_{m^t_\dn,m^t_\up}$, $A_\ell\subseteq
\prod\limits_{i=m^t_\dn}^{m^t_\up-1}\bH(i)$, $A_\ell\preceq\Delta_\ell$,
\item $\Delta_0=\Delta_t$, 
\item $8^{2^{n+2-\ell}-1}\leq \hn(\Delta_{\ell+1})=\HN(\Delta_{\ell+1})=
\HN(A_{\ell+1})<8^{2^{n+2-\ell}}$ (for $\ell<n$),
\item $(\forall f\in A_\ell)(\forall g\in\bigcup\limits_{k<\ell}\Delta_k)(
g\nsubseteq f)$.
\end{enumerate}
There are no problems for $\ell=0$ (note that there are practically no
restrictions on $A_0$, and $\Delta_0$ is determined). So suppose that
we have arrived to a stage $\ell+1\leq n$ of the construction. It follows
from \ref{basichn}(1,2) that $\hn(\bigcup\limits_{k\leq\ell}\Delta_k)>
8^{2^{n+3-\ell}-2}$. Let 
\[A^*=\{f\in\prod_{i=m^t_\dn}^{m^t_\up-1}\bH(i): (\forall g\in 
\bigcup_{k\leq\ell}\Delta_k)(g\nsubseteq f)\}.\]
Necessarily $\HN(A^*)\leq 2$. Now, using $(*)$ and \ref{basichn}(2), we may
choose $A_{\ell+1}\subseteq A^*$ such that $8^{2^{n+2-\ell}-1}\leq
\HN(A_{\ell+1})<8^{2^{n+2-\ell}}$. Next, we pick $\Delta_{\ell+1}\in
\cK^\bH_{m^t_\dn,m^t_\up}$ satisfying $A_{\ell+1}\preceq\Delta_{\ell+1}$ and
$\HN(A_{\ell+1})=\hn(\Delta_{\ell+1})=\HN(\Delta_{\ell+1})$. This finishes
the construction.

For $\ell\leq n$ let $s_\ell\in\Sigma_\bH(t)$ be such that $\Delta_{s_\ell}=
\bigcup\limits_{k\leq\ell}\Delta_k$. It should be clear that 
\begin{enumerate}
\item[$(\alpha)$] $s_0=t$, $s_{\ell+1}\in\Sigma(s_\ell)$, $2^{n+1-\ell}<
\nor[s_{\ell+1}]\leq 2^{n-\ell+2}$ (for $\ell<n$).
\end{enumerate}
We claim that additionally
\begin{enumerate}
\item[$(\beta)$]  if $s\in K_\bH$, $m^s_\dn=m^t_\dn$, $m^s_\up=m^t_\up$ and
$\nor[s]>2^{n-\ell+3}$, $\ell<n$, then
\[(\forall u\in\dom(\val[s]))(\exists v)(\langle u,v\rangle\in\val[s]
\setminus \val[s_{\ell+1}])\]
\end{enumerate}
(what will finish the proof of the claim). So suppose $s\in K_\bH$,
$m^s_\dn=m^t_\dn$, $m^s_\up=m^t_\up$ and $\nor[s]>2^{n-\ell+3}$, $\ell<n$. 
Let $u\in\prod\limits_{i<m^t_\dn}\bH(i)$. If $\Delta_s=\emptyset$, then we
may take $v\in\prod\limits_{i<m^t_\up}\bH(i)$ such that $u\vartriangleleft
v$ and $(\exists f\in\Delta_t)(f\subseteq v)$, so clearly $\langle u,v
\rangle\in\val[s]\setminus \val[s_{\ell+1}]$. Assume now that $\Delta_s\neq
\emptyset$, so $\hn(\Delta_s)>8^{2^{n-\ell+3}}$. Since $\HN(A_{\ell+1})<
8^{2^{n+2-\ell}}<\hn(\Delta_s)$, there is $f\in A_{\ell+1}$ such that
$(\forall g\in\Delta_s)(g\nsubseteq f)$. Clearly $\langle u,u\conc f\rangle
\in \val[s]\setminus\val[s_{\ell+1}]$.
\end{proof}
Finally note that the forcing notion
$\bQ^*_\infty(K_\bH,\Sigma_\bH,\Sigma^\bot_\bH)$ is not trivial.  
\end{proof}

\begin{conclusion}
\label{notCon1}
Let $\bH:\omega\longrightarrow\cHa$ be such that $|\bH(n)|\geq 2$ for all
$n<\omega$. Then the forcing notion $\bQ^*_\infty(K_\bH,\Sigma_\bH,
\Sigma^\bot_\bH)$ (where $(K_\bH,\Sigma_\bH,\Sigma^\bot_\bH)$ is as defined
in \ref{bas628}) is $\sigma$-$*$--linked and it adds a dominating
real. Consequently it is not $\omega$--nw--nep (by \cite{Sh:711}).  
\end{conclusion}

If one looks at $\val[t]$ for $t\in K_\bH$ in \ref{bas628}, then it is clear 
that $\HN$ is more appropriate to determine the norms of creatures. We
presented \ref{bas628} as it is a direct relative of the forcing notion of
\cite[\S 3]{RoSh:628}. However, it seems that the following example presents
a nicer member of this family.

\begin{example}
\label{bas628bis}
Let $\bH:\omega\longrightarrow\cHa$, $|\bH(n)|\geq 2$ for all
$n\in\omega$.\\ 
We construct a $\otimes$--creating triple $(K_{\ref{bas628bis}},
\Sigma_{\ref{bas628bis}},\Sigma^\bot_{\ref{bas628bis}})$ for $\bH$ which: 
\begin{enumerate}
\item is almost normal (see the construction), regular (see \ref{regular}),
forgetful and super--gluing,    
\item has the cutting property,
\item is really finitary provided $\bH(n)$ is finite for each $n\in\omega$. 
\end{enumerate}
\end{example}

\begin{proof}[Construction] 
It is similar to \ref{bas628}, but instead of $\hn$ we use $\HN$. So
$K_{\ref{bas628bis}}$ consists of $t\in\CR[\bH]$ such that 
\begin{itemize}
\item $\dis[t]=(m^t_\dn,m^t_\up,\Delta_t)$ for some $\Delta_t\in
\cK^\bH_{m_0,m_1}\cup\{\emptyset\}$ such that, if $\Delta_t\neq\emptyset$,
$\hn^+(\Delta_t)>1$, 
\item $\val[t]=\{\langle u,v\rangle\in\prod\limits_{i<m^t_\dn}\bH(i)\times
\prod\limits_{i<m^t_\up}\bH(i): u\vartriangleleft v\ \&\ (\forall f\in
\Delta_t)(f\nsubseteq v)\}$,
\item if $\Delta_t=\emptyset$ then $\nor[t]=\log_8(m^t_\up-m^t_\dn)+2
m^t_\dn+1$, otherwise $\nor[t]=\log_8(\HN(\Delta_t))$.  
\end{itemize}
For $t_0,\ldots,t_n\in K_{\ref{bas628bis}}$ such that $m^{t_i}_\up=
m^{t_{i+1}}_\dn$ (for $i<n$) let 
\[\Sigma_{\ref{bas628bis}}(t_0,\ldots,t_n)=\{t\in K_{\ref{bas628bis}}:
m^t_\dn=m^{t_0}_\dn\ \&\ m^t_\up = m^{t_n}_\up\ \&\ \bigcup_{i\leq n}
\Delta_{t_i}\subseteq\Delta_t\}.\]
For $t\in K_{\ref{bas628bis}}$, $\Sigma^\bot_{\ref{bas628bis}}(t)$ consists
of all finite sets $\{s_0,\ldots,s_n\}\subseteq K_{\ref{bas628bis}}$ such
that   
\begin{itemize}
\item $m^t_\dn=m^{s_0}_\dn<m^{s_0}_\up=m^{s_1}_\dn<\ldots<m^{s_{n-1}}_\up
=m^{s_n}_\dn<m^{s_n}_\up=m^t_\up$,\quad and 
\item $(\forall f\in\Delta_t)(\exists\ell\leq n)(f\rest [m^{s_\ell}_\dn,
m^{s_\ell}_\up)\in \Delta_{s_\ell})$.
\end{itemize}
Clearly $(K_{\ref{bas628bis}},\Sigma_{\ref{bas628bis}},
\Sigma^\bot_{\ref{bas628bis}})$ is a forgetful super--gluing
$\otimes$--creating triple. It is linked, and it is really finitary provided 
$\bH(n)$ is finite for each $n\in\omega$.    

\begin{claim}
\label{cl8}
Suppose $s,t\in K_{\ref{bas628bis}}$ are such that $\nor[s]<\nor[t]$.
Then 
\[(\forall u\in\dom(\val[t]))(\exists v)(\langle u,v\rangle\in\val[t]
\setminus\val[s]).\]
\end{claim}

\begin{proof}[Proof of the claim]
We may assume that $m^s_\dn=m^t_\dn$, $m^s_\up=m^t_\up$ (otherwise
trivial). It follows from the assumptions (and the definition of $\nor[s]$)
that $\Delta_s\neq\emptyset$. If $\Delta_t=\emptyset$, then the conclusion
is immediate, so assume $\Delta_t\neq\emptyset$. Thus $\Delta_s,\Delta_t\in 
\cK^\bH$ and $\HN(\Delta_s)<\HN(\Delta_t)$. Choose $\Delta\in\cK^\bH$ such
that $\Delta_t\preceq\Delta$ and $\hn(\Delta)=\HN(\Delta_t)$. Note that, by
\ref{basichn}(1), we have $\hn(\Delta)=\HN(\Delta)=\hn^+(\Delta)$. 
Consequently, we may choose $\Delta'\in\cK^\bH$ such that elements of
$\Delta'$ have pairwise disjoint domains, $\Delta\preceq\Delta'$ and
$\hn(\Delta')=\hn(\Delta)$. Now, for some $f\in\Delta_s$ we have $(\forall
g\in\Delta')(g\nsubseteq f)$. By the choice of $\Delta'$ we may build
$f^*\in\prod\limits_{i=m^s_\dn}^{m^s_\up-1}\bH(i)$ such that $f\subseteq
f^*$ and $(\forall g\in\Delta')(g\nsubseteq f^*)$. Since $\Delta_t\preceq
\Delta'$ we are done.
\end{proof}

Using \ref{cl8} we see that $(K_{\ref{bas628bis}},\Sigma_{\ref{bas628bis}})$ 
is almost normal in the following sense: the reducibility demand (see
\ref{reducible}(1)) holds for those $t\in K_{\ref{bas628bis}}$ for which
$\Delta_t\neq\emptyset$. However, this is enough to carry out, e.g., the
proof of \ref{adddom}(3) with almost no changes.  

\begin{claim}
\label{cl9}
$(K_{\ref{bas628bis}},\Sigma_{\ref{bas628bis}},\Sigma^\bot_{\ref{bas628bis}} 
)$ is regular and has the cutting property.
\end{claim}

\begin{proof}[Proof of the claim]
First we show the regularity. So suppose that 
\[(w,t_0,\ldots,t_n), (u,s_0,\ldots,s_m)\in \FC(K_{\ref{bas628bis}},
\Sigma_{\ref{bas628bis}},\Sigma^\bot_{\ref{bas628bis}})\]
are such that $(w,t_0,\ldots,t_n)\leq (u,s_0,\ldots,s_m)$, $m^{s_m}_\up
\leq m^{t_n}_\up$, $\nor[s_0]\geq 3$, and $\nor[t_\ell]\geq 3$, where
$\ell\leq n$ is such that $m^{t_\ell}_\dn<m^{s_0}_\dn<m^{t_\ell}_\up\leq
m^{s_0}_\up$. It follows from the definition of $\leq$ (see \ref{forcing})
that   
\[(\forall f\in \Delta_{t_\ell})\big(f\rest m^{s_0}_\dn\nsubseteq u\ \mbox{
or }\ (\exists g\in \Delta_{s_0})(g\subseteq f)\big).\]
Let $t',t''\in K_{\ref{bas628bis}}$ be such that $m^{t'}_\dn=m^{t_\ell}_\dn$,
$m^{t'}_\up=m^{t''}_\dn=m^{s_0}_\dn$, $m^{t''}_\up=m^{t_\ell}_\up$ and 
\[\begin{array}{lcl}
\Delta_{t'}&=&\{f\rest [m^{t_\ell}_\dn,m^{s_0}_\dn): f\in\Delta_{t_\ell}\ \&\ 
f\rest m^{s_0}_\dn\nsubseteq u\},\\
\Delta_{t''}&=&\{f\rest [m^{s_0}_\dn,m^{t_\ell}_\up): f\in\Delta_{t_\ell}\
\&\  f\rest m^{s_0}_\dn\subseteq u\}.
  \end{array}\]
Then $\{t',t''\}\in\Sigma^\bot_{\ref{bas628bis}}(t_\ell)$ and clearly
$u\in\pos(w,t_0,\ldots,t_{\ell-1},t')$. The only thing left is 
to show that the norm of $t''$ is at least 2. If $\Delta_{t''}=\emptyset$,
then it is clearly true as $m^{t''}_\dn\geq 1$. So suppose that
$\Delta_{t''}\neq \emptyset$ and we have to argue that
$\HN(\Delta_{t''})\geq 64$. But this is clear as $\Delta_{t''}\preceq
\Delta_{s_0}$.

Now let us show that $(K_{\ref{bas628bis}},\Sigma_{\ref{bas628bis}},
\Sigma^\bot_{\ref{bas628bis}})$ has the cutting property. Let $t\in
K_{\ref{bas628bis}}$, $\nor[t]>1$, $m^t_\dn<m<m^t_\up$. Choose $\Delta\in
\cK^\bH_{m^t_\dn,m^t_\up}$ such that elements of $\Delta$ have pairwise
disjoint domains, $\Delta_t\preceq\Delta$ and $\HN(\Delta_t)=\hn(\Delta)$
(like in the proof of \ref{cl8}). Put 
\[\begin{array}{lcl}
\Delta^0&=&\{f\rest [m^t_\dn,m): f\in\Delta\ \&\ |\dom(f)\cap [m^t_\dn,
m)|\geq\frac{1}{2}|\dom(f)|\},\\
\Delta^1&=&\{f\rest [m,m^t_\up): f\in\Delta\ \&\ |\dom(f)\cap [m,
m^t_\up)|\geq\frac{1}{2}|\dom(f)|\}.
  \end{array}\]
Let $s_0,s_1\in K_{\ref{bas628bis}}$ be such that $m^{s_0}_\dn=m^t_\dn$,
$m^{s_0}_\up=m=m^{s_1}_\dn$, $m^{s_1}_\up=m^t_\up$ and 
\[\begin{array}{lcl}
\Delta_{s_0}&=&\{g\rest [m^t_\dn,m): g\in\Delta_t\ \&\ (\exists
f\in\Delta^0)(f\subseteq g)\},\\
\Delta_{s_1}&=&\{g\rest [m,m^t_\up): g\in\Delta_t\ \&\ (\exists
f\in\Delta^1)(f\subseteq g)\}.
  \end{array}\]
Now check.
\end{proof}
\end{proof}

\begin{conclusion}
\label{notCon2}
The forcing notion $\bQ^*_\infty(K_{\ref{bas628bis}},
\Sigma_{\ref{bas628bis}},\Sigma^\bot_{\ref{bas628bis}})$
is $\sigma$-$*$--linked Souslin ccc and it adds a dominating
real. Consequently it adds a Cohen real (by \cite{Sh:480}) and it is not
$\omega$--nw--nep (by \cite{Sh:711}). 
\end{conclusion}

Hall's norms are of special interest because of ``gluing and cutting'', but
we may use them to build local creating pairs to.

\begin{example}
\label{loc628}
Let $\bH^*:\omega\longrightarrow\cHa$, $2\leq |\bH^*(i)|<\omega$ for all
$i\in\omega$. Suppose that $\bar{n}=\langle n_k:k<\omega\rangle\subseteq
\omega$ is an increasing sequence such that $\lim\limits_{k\to\infty}
n_{k+1}-n_k=\infty$ and let $\bH=\bH^*[\bar{n}]:\omega\longrightarrow\cHa$
be defined by $\bH(k)=\prod\limits_{i=n_k}^{n_{k+1}-1}\bH^*(i)$. 

We construct a really finitary creating pair $(K_{\ref{loc628}},
\Sigma_{\ref{loc628}})$ for $\bH$ which is local, forgetful, normal and
Cohen--producing (see \ref{Cohen}(1)). 
\end{example}

\begin{proof}[Construction] 
Let $K_{\ref{loc628}}$ consist of creatures $t\in\CR[\bH]$ such that 
\begin{itemize}
\item $\dis[t]=(m_t,\Delta_t)$ for some $m_t<\omega$ and $\Delta_t\in
\cK^{\bH^*}_{n_{m_t},n_{m_t+1}}\cup\{\emptyset\}$ such that, if $\Delta_t
\neq\emptyset$, $\hn^+(\Delta_t)>1$, 
\item $\val[t]=\{\langle u,v\rangle\in\prod\limits_{i<m_t}\bH(i)\times
\prod\limits_{i\leq m_t}\bH(i): u\vartriangleleft v\ \&\ (\forall f\in
\Delta_t)(f\nsubseteq v(m_t))\}$,
\item if $\Delta_t\neq\emptyset$, then $\nor[t]=\log_8(\HN(\Delta_t))$, 
otherwise $\nor[t]=\log_8(n_{m_t+1}-n_{m_t}+1)$.
\end{itemize}
The operation $\Sigma_{\ref{loc628}}$ gives non-empty results for singletons 
only and 
\[\Sigma_{\ref{loc628}}(t)=\{s\in K_{\ref{loc628}}:m^t_\dn=m^s_\dn\ \&\
\Delta_t\subseteq \Delta_s\}.\]
Clearly, $(K_{\ref{loc628}},\Sigma_{\ref{loc628}})$ is a really finitary
creating pair which is local, forgetful and linked (remember
\ref{basichn}(2)).   
\begin{claim}
\label{cl10}
$(K_{\ref{loc628}},\Sigma_{\ref{loc628}})$ is normal and
Cohen--producing.
\end{claim}

\begin{proof}[Proof of the claim]
First note that $(K_{\ref{loc628}},\Sigma_{\ref{loc628}})$ is
reducible (remember $(*)$ of the construction for \ref{bas628}; note that
here there are no problems caused by $\Delta_t=\emptyset$). Next note that
(the proof of) \ref{cl8} applies here too.

To show that $(K_{\ref{loc628}},\Sigma_{\ref{loc628}})$ is
Cohen--producing fix (for $k\in\omega$) an $a_k\in\bH^*(n_k)$ and let 
\[A_k=\{f\in\prod_{i=n_k}^{n_{k+1}-1}\bH^*(i):
f(n_k)=a_k\}\subseteq\bH(k).\]
Suppose that $t\in K_{\ref{loc628}}$, $\nor[t]>1$, $m^t_\dn=k$ and let
$u\in\dom(\val[t])$. If $\Delta_t=\emptyset$, then we easily choose $v_0,v_1 
\in\rng(\val[t])$ extending $u$ and such that $v_0(k)\notin A_k$, $v_1(k)\in
A_k$ (remember $|\bH^*(n_k)|\geq 2$). So suppose that $\Delta_t\neq
\emptyset$ and thus $\HN(\Delta_t)>8$. Then we may find $\Delta'\in
\cK^{\bH^*}_{n_k,n_{k+1}}$ such that $\Delta_t\preceq\Delta'$, elements of
$\Delta'$ have pairwise disjoint domains all of size $>2$. Now we easily
build $v_0,v_1\in\dom(\rng[t])$, both extending $u$ and such that $v_0(k)
\notin A_k$, $v_1(k)\in A_k$.
\end{proof}
\end{proof}

\begin{conclusion}
\label{conclOne}
Assume $\bH^*,\bar{n},\bH$ and $(K_{\ref{loc628}},
\Sigma_{\ref{loc628}})$ are as in \ref{loc628}. 
\begin{enumerate}
\item Suppose $f:\omega\times\omega\longrightarrow\omega$ is a fast function
such that 
\[(\forall k\in\omega)(\forall^\infty\ell\in\omega)(f(k,\ell)<\log_8(n_{\ell
+1}-n_\ell))\]
(e.g.~$f(k,\ell)=2^{2k}$). Then $\bQ^*_f(K_{\ref{loc628}},
\Sigma_{\ref{loc628}})$ is a non-trivial $\sigma$-$*$--linked Borel ccc
forcing notion which adds a dominating real (so it adds a Cohen real and it
is not $\omega$--nw--nep).
\item Let $h(n,m)=\max\{0,m-1\}$ (so $h:\omega\times\omega\longrightarrow
\omega$ is a regressive function). Suppose that $\cF\subseteq
(\omega\setminus 2)^{\textstyle\omega}$ is a countable $h$--closed
$\geq^*$--directed family such that $(\forall f\in\cF)(\forall^\infty\ell
\in\omega)(f(\ell)<\log_8(n_{\ell+1}-n_\ell))$ (e.g.~$\cF=\{f_k:k<\omega\}$,
$f_k(\ell)=\max\{2,\lfloor\frac{1}{k}\log_8(n_{\ell+1}-n_\ell)\rfloor\}$). 
Then $\bQ^*_\cF(K_{\ref{loc628}},\Sigma_{\ref{loc628}})$ is a
non-trivial $\sigma$--$*$--linked very Borel ccc forcing notion, it adds a 
Cohen real and it is nice (so it preserves unbounded families).
\end{enumerate}
\end{conclusion}

Our next examples generalize (in some sense) the Eventually Different Real
Forcing of Miller \cite{Mi81}.

\begin{example}
\label{EDRF}
Let $\bH:\omega\longrightarrow\cHa$, $|\bH(k)|\geq 4$. For $k<\omega$, let
$N_k$ be $|\bH(k)|$ if $\bH(k)$ is finite, and $2^{k+2}$ otherwise. Assume
$\lim\limits_{k\mapsto\infty}N_k=\infty$. Let $h:\omega\times\omega
\longrightarrow\omega$ be given by
\[h(k,n)=\left\{\begin{array}{ll}
n-1   &\mbox{ if }n\geq N_k,\\
2n-N_k&\mbox{ if }\frac{7}{8}N_k< n<N_k,\\
1     &\mbox{ otherwise.}
		\end{array}\right.\]
(Note that $h$ is regressive.) We construct an $h$--linked creating pair
$(K_{\ref{EDRF}},\Sigma_{\ref{EDRF}})$ for $\bH$ which is local, forgetful,
Cohen--producing and of the {\rm BCB}--type (see \ref{Cohen}).
\end{example}

\begin{proof}[Construction]
Let $K_{\ref{EDRF}}$ be the collection of all $t\in\CR[\bH]$ such that 
\begin{itemize}
\item $\dis[t]=(k_t,E_t)$ for some $k_t<\omega$ and $E_t\subseteq\bH(k_t)$
such that $0<|E_t|<N_{k_t}$,
\item $\val[t]=\{\langle u,v\rangle\in\prod\limits_{i<k_t}\bH(i)\times
\prod\limits_{i\leq k_t}\bH(i): u\vartriangleleft v\ \&\ v(k_t)\notin E_t\}$,
\item if $|E_t|\geq\frac{1}{4}N_{k_t}$ then $\nor[t]=1$; otherwise $\nor[t]=
N_{k_t}-|E_t|$. 
\end{itemize}
The operation $\Sigma_{\ref{EDRF}}$ is natural: it gives non-empty results
for singletons only and 
\[\Sigma_{\ref{EDRF}}(t)=\{s\in K_{\ref{EDRF}}: k_s=k_t\ \&\ E_t\subseteq
E_s\}.\]
It should be clear that $(K_{\ref{EDRF}},\Sigma_{\ref{EDRF}})$ is a local
forgetful creating pair for $\bH$. 

To show that it is $h$--linked suppose that $k>1$, $t_0,t_1\in
K_{\ref{EDRF}}$, $\nor[t_0],\nor[t_1]\geq k$ and $\ell=m^{t_0}_\dn=
m^{t_1}_\dn$. Then $|E_{t_0}|,|E_{t_1}|<\frac{1}{4}N_\ell$ and thus $0<|
E_{t_0}\cup E_{t_1}|<N_\ell$. Let $s\in K_{\ref{EDRF}}$ be such that $k_s
=\ell$, $E_s=E_{t_0}\cup E_{t_1}$. Clearly $s\in\Sigma_{\ref{EDRF}}(t_0)
\cap\Sigma_{\ref{EDRF}}(t_1)$. If $h(\ell,k)=1$ then clearly $\nor[s]\geq
h(\ell,k)$, so suppose $h(\ell,k)>1$. Necessarily $\frac{7}{8} N_\ell<k$, so
$|E_{t_0}|,|E_{t_1}|<\frac{1}{8} N_\ell$ and therefore $|E_s|<\frac{1}{4}
N_\ell$. Hence 
\[\nor[s]=N_\ell-|E_s|\geq N_\ell-|E_{t_0}|-|E_{t_1}|\geq 2k-N_\ell= h(\ell,
k).\]

Let us show now that $(K_{\ref{EDRF}},\Sigma_{\ref{EDRF}})$ is
Cohen--producing. For each $n\in\omega$ choose a set $A_n\subseteq\bH(n)$
such that $|A_n|=\lfloor\frac{1}{2}|\bH(n)|\rfloor$ if $\bH(n)$ is finite,
and $A_n$ is infinite co-infinite if $\bH(n)$ is infinite. Suppose that
$t\in K_{\ref{EDRF}}$, $\nor[t]>1$. Then $E_t\subseteq\bH(k_t)$, $|E_t|<
\frac{1}{4} N_{k_t}$, so we may choose $a_1\in A_{k_t}\setminus E_{k_t}$ and
$a_0\in\bH(k_t)\setminus (A_{k_t}\cup E_{k_t})$, and we easily finish.

Finally, let us argue that $(K_{\ref{EDRF}},\Sigma_{\ref{EDRF}})$ is of the
{\rm BCB}--type, To this end suppose that $\langle s_n:n<\omega\rangle
\subseteq K_{\ref{EDRF}}$, $m^{s_n}_\dn=k_{s_n}=\ell$, $\nor[s_n]\geq 2$.
Then $|E_{s_n}|<\frac{1}{4} N_\ell$ for each $n$.\\
If $\bH(\ell)$ is finite, then the demand in $(\circledast^{\rm BCB})$ of
\ref{Cohen}(2) is trivially satisfied (just take $\{a_0,\ldots,a_m\}=
\bH(\ell)$).\\
So suppose that $\bH(\ell)$ is infinite and to simplify notation let
$\bH(\ell)=\omega$. Let $E_{s_n}=\{b^n_0,\ldots,b^n_{k_n-1}\}$ be the
increasing enumeration; $k_n=|E_{s_n}|$. We may find an infinite set $Y
\subseteq\omega$ and $k^*\leq k<\omega$ such that 
\begin{itemize}
\item $k_n=k$ for each $n\in Y$,
\item $\langle b^n_i:n\in Y\rangle$ is constant for each $i<k^*$,
\item $\langle b^n_i:n\in Y\rangle$ is strictly increasing for each $i\in
[k^*,k)$.
\end{itemize}
Suppose $a\in\bH(\ell)\setminus\{b^n_i:i<k^*\}$ for some (equivalently: all)
$n\in Y$. Then, for sufficiently large $n\in Y$, for every $i\in [k^*,k)$ we
have $b^n_i>a$. Consequently $(\forall^\infty n\in Y)(a\notin E_{s_n})$, so
we may easily finish.
\end{proof}

\begin{conclusion}
\label{conclTwo}
Let $\bH$, $\langle N_k: k<\omega\rangle$, $h$ and $(K_{\ref{EDRF}},
\Sigma_{\ref{EDRF}})$ be as in \ref{EDRF}. Suppose that $\cF\subseteq
(\omega\setminus 2)^{\textstyle \omega}$ is a countable $h$--closed
$\geq^*$--directed family such that for some $f\in\cF$ we have
$(\forall^\infty k\in\omega)(f(k)<N_k)$ (e.g.~$\cF=\{f_\ell:\ell<\omega\}$,
$f_\ell(k)=N_k-2^\ell$ if $N_k>2^{\ell+1}$, $f_\ell(k)=2$ otherwise).  Then
$\bQ^*_\cF(K_{\ref{EDRF}},\Sigma_{\ref{EDRF}})$ is a non-trivial
$\sigma$--$*$--linked Borel ccc forcing notion which adds a Cohen real
and is nice (so it preserves unbounded families). 
\end{conclusion}

If $\bH$ and $\cF$ are as in \ref{conclTwo}, and $\bH(k)$ is finite for each
$k$, then we may use \ref{allsimple}(3) to get that the forcing notion
$\bQ^*_\cF(K_{\ref{EDRF}},\Sigma_{\ref{EDRF}})$ is very Borel ccc. We may
prove the same conclusion without the additional assumption on $\bH$ (see
\ref{extra} below). Unfortunately, this proof is very specific for
$\bQ^*_\cF(K_{\ref{EDRF}},\Sigma_{\ref{EDRF}})$ and it does not generalize
to cover more forcing notions of the form $\bQ^*_\cF(K,\Sigma)$

\begin{proposition}
\label{extra}
Assume that $\bH$, $\langle N_k:k<\omega\rangle$, $h$, $K=K_{\ref{EDRF}}$
and $\Sigma=\Sigma_{\ref{EDRF}}$ are as in \ref{EDRF}, and $\cF\subseteq
(\omega\setminus 2)^{\textstyle \omega}$ is a countable $h$--closed
$\leq^*$--directed family such that 
\[(\exists f\in\cF)(\forall^\infty k\in\omega)(f(k)<N_k).\]
Then the forcing notion $\bQ^*_\cF(K_{\ref{EDRF}},\Sigma_{\ref{EDRF}})$ is
very Borel ccc. 
\end{proposition}

\begin{proof}
The only thing that should be shown is that being a maximal antichain is a
Borel relation (remember \ref{conclTwo} and so \ref{allsimple}(1c)). Put
$Z=\bigcup\limits_{i<\omega}\prod\limits_{j<i}\bH(i)\setminus \{\langle
\rangle\}$, $\cX=(Z\cup\{\langle\rangle\})\times K^{\textstyle\omega}$, and
$\cY=\cX^{\textstyle\omega}$ and $\cZ=\omega^{\textstyle Z}\times
\omega^{\textstyle Z}$. 

Then $\cX,\cY,\cZ$ are Polish spaces (each equipped with the respective
product topology), $\bQ^*_\emptyset(K,\Sigma)$ and $\bQ^*_\cF(K,\Sigma)$ are
$\Pi^0_1$ and $\Sigma^0_2$ subsets of $\cX$, respectively, and for our
conclusion it is enough to show that 
\[\begin{array}{ll}
P\stackrel{\rm def}{=}\Big\{\langle p_n:n<\omega\rangle\in\cY:&
  (\forall n<\omega)(p_n\in \bQ^*_\cF(K,\Sigma))\ \mbox{ and }\\
& \{p_n:n<\omega\}\mbox{ is pre-dense in }\bQ^*_\cF(K,\Sigma)\ \Big\}
  \end{array}\]
is a Borel subset of $\cY$.

Plainly we may assume that 
\[(\forall f\in\cF)(\forall k<\omega)(2\leq f(k)<N_k)\]
(as we may modify suitable the family $\cF$ without changing the forcing
notion). Now, for $f\in\cF$ we define
\[C_f=\{(h_0,h_1)\in\cZ:(\forall\eta\in Z)(h_0(\eta)<\lh(\eta)\ \&\
h_1(\eta)<N_{h_0(\eta)}-f(h_0(\eta)))\}.\]
Clearly each $C_f$ is a compact subset of $\cZ$.

\begin{claim}
\label{cl25}
Suppose that $\bar{p}=\langle p_n:n<\omega\rangle\in\cY$, where $p_n\in
\bQ^*_\cF(K,\Sigma)$ are such that $\nor[t^{p_n}_k]>1$ (for all
$k,n<\omega$). Then the following are equivalent:
\begin{enumerate}
\item[(A)$_{\bar{p}}$] There is a condition $p\in \bQ^*_\cF(K,\Sigma)$
incompatible with every $p_n$ (for $n<\omega$).
\item[(B)$_{\bar{p}}$] There are $eta\in Z$ and $f\in\cF$ and $(h_0,h_1)\in
C_f$ such that
\begin{enumerate}
\item[(i)]  for every $\nu\in Z$, if $\eta\vartriangleleft \nu$ then
$h_0(\nu)\geq \lh(\eta)$, and  
\item[(ii)] for every $\nu_0,\nu_1\in Z$, if $\eta\vartriangleleft \nu_0$,
$\eta\vartriangleleft \nu_1$ and $\nu_0,\nu_1\in\bigcup\limits_{n<\omega}
\POS(p_n)$ and $h_0(\nu_0)=h_0(\nu_1)$ and $h_1(\nu_0)=h_1(\nu_1)$, then
$\nu_0(h_0(\nu_0))=\nu_1(h_0(\nu_1))$. 
\end{enumerate}
\end{enumerate}
\end{claim}

\begin{proof}[Proof of the claim]
Assume (A)$_{\bar{p}}$ and pick $p\in \bQ^*_\cF(K,\Sigma)$ and $f\in\cF$
such that 
\[(\forall k<\omega)(\nor[t^p_k]=N_{\lh(w^p)+k}-|E_{t^p_k}|=f(\lh(w^p)+
k)),\]
and $p$ is incompatible with all $p_n$ (for $n<\omega$). For $\ell=\lh(w^p)
+k$, $k<\omega$, let $\langle x^\ell_m:m<N_\ell-f(\ell)\rangle$ be an
enumeration of $E_{t^p_k}$. Note that, as $p\bot p_n$, if $\nu\in\POS(p_n)$
and $w^p\vartriangleleft \nu$, then for some $h_0(\nu)\in [\lh(w^p),
\lh(\nu))$ and $h_1(\nu)<N_{h_0(\nu)}-f(h_0(\nu))$ we have $\nu(h_0(\nu))=
x^{h_0(\nu)}_{h_1(\nu)}$. Letting $\eta=w^p$ we may now easily define
$(h_0,h_1)$ so that $\eta,f,(h_0,h_1)$ witness (B)$_{\bar{p}}$.

Suppose now that $\eta,f$ and $(h_0,h_1)$ witness (B)$_{\bar{p}}$. Let
$w^p=\eta$ and for $\ell=\lh(w^p)+k$, $k<\omega$ and $m<N_\ell-f(\ell)$ let
$x^\ell_m\in\bH(\ell)$ be such that 
\begin{enumerate}
\item[$(*)$] if $\nu\in \bigcup\limits_{n<\omega}\POS(p_n)$, $w^p
\vartriangleleft\nu$ and $h_0(\nu)=\ell$ and $h_1(\nu)=m$,\\
then $x^\ell_m=\nu(\ell)$.
\end{enumerate}
(The choice of the $x^\ell_m$'s is possible by (B)$_{\bar{p}}$(ii).) Now,
for each $k<\omega$ pick $t^p_k\in K$ so that 
\[m^{t^p_k}_\dn=\lh(w^p)+k\quad \mbox{ and }\quad E_{t^p_k}=\{x^{\lh(w^p)
+k}_m: m<N_{\lh(w^p)+k}-f(\lh(w^p)+k)\}.\]
Notice that for sufficiently large $k$ we have $f(\lh(w^p)+k)>\frac{7}{8}
N_{\lh(w^p)+k}$, so for those $k$ we will also have $\nor[t^p_k]= f(\lh(w^p)
+k)$. Hence $p\in \bQ^*_\cF(K,\Sigma)$ and easily it is a condition
incompatible with all $p_n$'s.
\end{proof}

For $\eta\in Z$ and $f\in\cF$, let $B^{\eta,f}$ consist of all $(h_0,h_1,
\bar{p})\in\cZ\times\cY$ such that 
\begin{itemize}
\item $(h_0,h_1)\in C_f$, $\bar{p}=\langle p_n:n<\omega\rangle\in\cY$, where
$p_n\in\bQ^*_\emptyset(K,\Sigma)$ are such that $\nor[t^{p_n}_k]>1$ for all
$k,n<\omega$, and  
\item $(\forall \nu\in Z)(\eta\vartriangleleft\nu\ \Rightarrow\ \lh(\eta)
\leq h_0(\nu))$, and 
\item for all $\nu_0,\nu_1\in Z\cap\bigcup\limits_{n<\omega}\POS(p_n)$ we
have 
\[h_0(\nu_0)=h_0(\nu_1)\ \&\ h_1(\nu_0)=h_1(\nu_1)\quad\Rightarrow\quad
\nu_0(h_0(\nu_0))=\nu_1(h_1(\nu_1)).\] 
\end{itemize}
It should be clear that $B^{\eta,f}$ is a closed subset of $C_f\times \cY$
and hence (as $C_f$ is compact) the set 
\[A^{\eta,f}=\{\bar{p}\in\cY: (\exists (h_0,h_1)\in C_f)((h_0,h_1,\bar{p})
\in B^{\eta,f})\}\]
is a closed subset of $\cY$, and $A^{\eta,f}\cap \big(\bQ^*_\cF(K,\Sigma)
\big)^{\textstyle\omega}$ is Borel. Now, pick a Borel function
$\pi:\cY\longrightarrow\cY$ such that 

if $\bar{p}=\langle p_n:n<\omega\rangle\in \big(\bQ^*_\cF(K,\Sigma)\big)^{
\textstyle\omega}$ and $\bar{q}=\langle q_n:n<\omega\rangle =\pi(\bar{p})$, 

then
\begin{itemize}
\item for each $n<\omega$, $q_n\in \bQ^*_\cF(K,\Sigma)$ and $(\forall
k<\omega)(\nor[t^{q_n}_k]>1)$, 
\item for each $n<\omega$, for some $m<\omega$ and $k<\omega$ we have
\[w^{q_n}\in\pos(w^{p_m},t^{p_m}_0,\ldots,t^{p_m}_k)\quad\mbox{ and }\quad
q_n=(w^{q_n},t^{p_m}_{k+1},t^{p_m}_{k+2},\ldots),\]
\item if $m<\omega$, $k<\omega$, $w\in\pos(w^{p_m},t^{p_m}_0,\ldots,
t^{p_m}_k)$, and $(\forall\ell>k)(\nor[t^{p_m}_\ell]>1)$,\\
then $(w,t^{p_m}_{k+1},t^{p_m}_{k+2},\ldots )\in\{q_n:n<\omega\}$.
\end{itemize}
Clearly, $\bar{p}\in \big(\bQ^*_\cF(K,\Sigma)\big)^{\textstyle\omega}$ is
pre-dense if and only if so is $\pi(\bar{p})$. Hence, for $\bar{p}\in
\big(\bQ^*_\cF(K,\Sigma)\big)^{\textstyle\omega}$ we have

$\bar{p}$ is pre-dense\qquad if and only if 

there are $\eta\in Z$ and $f\in\cF$ such that $\pi(\bar{p})\in A^{\eta,f}$\\
(remember \ref{cl25}). Since both $Z$ and $\cF$ are countable, the proof of
the proposition is completed. 
\end{proof}

The construction presented in \ref{EDRF} is a particular case of a more
general method of building linked creating pairs from some of the examples
presented in \cite{RoSh:470}. First let us recall the following definition.

\begin{definition}
[See {\cite[Def. 5.2.5]{RoSh:470}}]
\label{additive}
Let $(K,\Sigma)$ be a creating pair. We say that a creature $t\in K$ is {\em
$(n,m)$--additive} if for all $t_0,\ldots,t_{n-1}\in \Sigma(t)$ such that
$\nor[t_i]\leq m$ (for $i<n$) there is $s\in\Sigma(t)$ such that   
\[t_0,\ldots,t_{n-1}\in\Sigma(s)\quad\mbox{ and }\quad\nor[s]\leq\max\{
\nor[t_\ell]: \ell<n\}+1.\]
\end{definition}

\begin{example}
\label{dual}
Suppose that $(K,\Sigma)$ is a local and forgetful creating pair for
$\bH$, and it satisfies the demand $(\boxplus)$ of \ref{reducible}(2). Let
$\bar{t}^*=\langle t^*_0,t^*_1,t^*_2,\ldots\rangle\in\PC_\infty(K,\Sigma)$
be such that each $t^*_n$ is $(2,\nor[t^*_n])$--additive and $m^{t_0}_\dn=
0$. We construct a local linked creating pair $(K^c_{\bar{t}^*},\Sigma^c_{
\bar{t}^*})$ ({\em the $\bar{t}^*$--dual of $(K,\Sigma)$}).   
\end{example}

\begin{proof}[Construction]
For $n<\omega$ and a creature $t\in\Sigma(t^*_n)$ let a creature $t^c$ be
such that 
\begin{itemize}
\item $\nor[t^c]=\max\{0,\nor[t^*_n]-\nor[t]\}$, 
\item $\val[t^c]=\big(\prod\limits_{i<n}\bH(i)\times\prod\limits_{i\leq n}
\bH(i)\big)\setminus\val[t]$,  
\item $\dis[t^c]=(\dis[t],c)$.
\end{itemize}
[The creature $t^c$ is defined only if $\val[t^c]\neq\emptyset$.] Let
$K^c_{\bar{t}^*}$ be the collection of all (correctly defined) $t^c$ (for
$t\in\Sigma(t^*_n)$, $n<\omega$). For $t^c\in K^c_{\bar{t}^*}$ (defined as
above for $t\in\Sigma(t^*_n)$) we let 
\[\Sigma^c_{\bar{t}^*}(t^c)=\{s^c: t\in\Sigma(s)\ \&\ s\in\Sigma(t^*_n)\}.\]
\end{proof}

The examples of local creating pairs have their (local) tree--creating
variants too. They can be constructed like the following example.

\begin{example}
\label{loctree}
Let $\bH\in\baire$ be a strictly increasing function such that $\bH(0)>2$.
We construct a really finitary, normal (local) tree--creating pair
$(K_{\ref{loctree}},\Sigma_{\ref{loctree}})$ for $\bH$ which is
Cohen--producing.
\end{example}

\begin{proof}[Construction]
The family $K_{\ref{loctree}}$ consists of tree--creatures $t\in\TCR[\bH]$
such that
\begin{itemize}
\item $\dis[t]=(m_t,\eta_t,E_t)$ such that $m_t<\omega$, $\eta_t\in
\prod\limits_{i<m_t}\bH(i)$ and $\emptyset\neq E_t\subseteq\bH(m_t)$,
$E_t\neq\bH(m_t)$,
\item $\val[t]=\{\langle\eta_t,\nu\rangle: \eta_t\vartriangleleft\nu\in
\prod\limits_{i\leq m_t}\bH(i)\ \&\ \nu(m_t)\notin E_t\}$,
\item $\nor[t]=\log_4(\frac{\bH(m_t)}{|E_t|})$.
\end{itemize}
The tree composition $\Sigma_{\ref{loctree}}$ is natural: it gives non empty
results for singletons only and then
\[\Sigma_{\ref{loctree}}(t)=\{s\in K_{\ref{loctree}}: \eta_s=\eta_t\ \&\
E_t\subseteq E_s\}.\]
Now check.
\end{proof}

\begin{conclusion}
\label{conclThree}
Suppose $\bH\in\baire$ is strictly increasing, $\bH(0)>4$.
\begin{enumerate}
\item The forcing notion $\bQ^\tree_1(K_{\ref{loctree}},
\Sigma_{\ref{loctree}})$ is non-trivial, $\sigma$--$*$--linked and it adds a
dominating real.
\item Assume that $f:\omega\times\omega\longrightarrow\omega$ is a fast
function such that 
\[(\forall n,m<\omega)(f(n,m)<\log_4(\bH(m))).\]
Then the forcing notion $\bQ^\tree_f(K_{\ref{loctree}},
\Sigma_{\ref{loctree}})$ is non-trivial, $\sigma$--$*$--linked, Borel ccc, 
and it adds a dominating real.
\item Suppose that $\cF\subseteq (\omega\setminus 2)^{\textstyle\omega}$ is
a countable $\geq^*$--directed family such that 
\[\begin{array}{l}
(\exists f\in\cF)(\forall^\infty n\in\omega)(f(n)<\log_4(\bH(n)))\quad
\mbox{ and }\\ 
(\forall f\in\cF)(\exists g\in\cF)(\forall^\infty n\in\omega)(g(n)<f(n)-1).
  \end{array}\]
Then the forcing notion $\bQ^\tree_\cF(K_{\ref{loctree}},
\Sigma_{\ref{loctree}})$ is non-trivial, $\sigma$--$*$--linked, very Borel
ccc, and it adds a Cohen real and is nice.
\end{enumerate}
\end{conclusion}

\section{More constructions}
In this section we introduce more schemes for building ccc forcing notions
as well as more norm conditions that can be used in conjunctions with the
methods presented in the previous section.

\subsection{Mixtures with random}

\begin{definition}
\label{mixing}
Let $\bH:\omega\longrightarrow\cHa$. We say that $(K,\Sigma,\bF)$ is {\em a
mixing triple for $\bH$} if 
\begin{enumerate}
\item[(a)] $(K,\Sigma)$ is a (local) tree--creating pair for $\bH$,
\item[(b)] for each $\eta\in\bigcup\limits_{m\in\omega}\prod\limits_{i<m}
\bH(i)$ there is $t^*_\eta\in K$ such that 
\[(\forall t\in\TCR_\eta[\bH] \cap K)(t\in \Sigma(t^*_\eta)),\]
\item[(c)] $\bF=\langle F_\eta:\eta\in\bigcup\limits_{m\in\omega}
\prod\limits_{i<m}\bH(i)\rangle$, where for each $\eta\in\bigcup\limits_{m
\in\omega}\prod\limits_{i<m}\bH(i)$:
\item[(d)] $F_\eta:[0,1]^{\textstyle \pos(t^*_\eta)}\longrightarrow [0,1]$,
\item[(e)] if $\langle r_\nu:\nu\in\pos(t^*_\eta)\rangle,\langle r_\nu':
\nu\in\pos(t^*_\eta)\rangle\in [0,1]^{\textstyle \pos(t^*_\eta)}$, $r_\nu
\leq r_\nu'$ for all $\nu\in\pos(t^*_\eta)$, then $F_\eta(r_\nu:\nu\in\pos(
t^*_\eta))\leq F_\eta(r_\nu':\nu\in\pos(t^*_\eta))$,
\item[(f)] if $\langle r_\nu:\nu\in\pos(t^*_\eta)\rangle\in
[0,1]^{\textstyle \pos(t^*_\eta)}$, $\varepsilon>0$, then there are $r_\nu'
<r_\nu$ (for $\nu\in\pos(t^*_\eta)$) such that for each $\langle r_\nu'':
\nu\in\pos(t^*_\eta)\rangle\in [0,1]^{\textstyle \pos(t^*_\eta)}$ satisfying
$r_\nu'<r_\nu''\leq r_\nu$ (for $\nu\in\pos(t^*_\eta)$) we have
\[F_\eta(r_\nu:\nu\in\pos(t^*_\eta))-\varepsilon<F_\eta(r_\nu'':\nu\in
\pos(t^*_\eta)),\]
\item[(g)] if $r_\nu\geq\varepsilon>0$ for $\nu\in\pos(t^*_\eta)$ then
$F_\eta(r_\nu:\nu\in\pos(t^*_\eta))\geq\varepsilon$.
\end{enumerate}
\end{definition}

\begin{definition}
\label{mixfor}
Let $(K,\Sigma,\bF)$ be a mixing triple for $\bH$.
\begin{enumerate}
\item Let $T^*=T^*_{K,\Sigma}\subseteq\bigcup\limits_{m\in\omega}
\prod\limits_{i<m}\bH(i)$ be a tree such that 
\[\mrot(T^*)=\langle\rangle\quad\mbox{ and }\quad (\forall\eta\in T^*)(
\suc_{T^*}(\eta)=\pos(t^*_\eta)).\]
\item If $X\subseteq\pos(t^*_\eta)$, $\eta\in T^*$ and $\langle r_\nu:\nu\in
X\rangle\subseteq [0,1]$, then we define $F_\eta(r_\nu:\nu\in X)$ as
$F_\eta(r^*_\nu:\nu\in\pos(t^*_\eta))$, where 
\[r^*_\nu=\left\{\begin{array}{ll}
r_\nu&\mbox{ if }\nu\in X,\\
0    &\mbox{ if }\nu\in\pos(t^*_\eta)\setminus X.
		 \end{array}\right.\]
\item Suppose that $p=\langle t^p_\eta:\eta\in T^p\rangle\in
\bQ^\tree_\emptyset(K,\Sigma)$ and $A\subseteq T^p$ is a front of $T^p$. We
let $T[p,A]=\{\eta\in T^p:(\exists\rho\in A)(\eta\trianglelefteq\rho)\}$,
and we define $\mu^\bF_{p,A}=\mu_{p,A}:T[p,A]\longrightarrow [0,1]$ by
downward induction as follows:
\begin{itemize}
\item if $\eta\in A$ then $\mu_{p,A}(\eta)=1$,
\item if $\mu_{p,A}(\nu)$ has been defined for all $\nu\in\pos(t^p_\eta)$,
$\eta\in T[p,A]\setminus A$, then we put $\mu_{p,A}(\eta)=F_\eta(\mu_{p,A}
(\nu):\nu\in\pos(t^p_\eta))$.
\end{itemize}
\item For $p=\langle t^p_\eta:\eta\in T^p\rangle\in\bQ^\tree_\emptyset(K,
\Sigma)$ we define 
\[\mu^\bF(p)=\inf\{\mu_{p,A}(\mrot(p)):A\mbox{ is a front of }T^p\}.\]
\item Let $\bQ^\mtree_\emptyset(K,\Sigma,\bF)=\{p\in\bQ^\tree_\emptyset(K,
\Sigma):\mu^\bF(p)>0\}$ be equipped with the partial order inherited from
$\bQ^\tree_\emptyset(K,\Sigma)$. Similarly we define forcing notions
$\bQ^\mtree_1(K,\Sigma,\bF)$, $\bQ^\mtree_f(K,\Sigma,\bF)$, $\bQ^\mtree_\cF(
K,\Sigma,\bF)$ (for suitable $f$ and $\cF$). 
\end{enumerate}
\end{definition}

\begin{definition}
\label{ccccomplete}
A mixing triple $(K,\Sigma,\bF)$ is {\em ccc--complete\/} if
\begin{enumerate}
\item[(a)] for each $\eta\in T^*_{K,\Sigma}$ and $A\subseteq\pos(t^*_\eta)$, 
there is a unique tree--creature $t_A\in\Sigma(t^*_\eta)$ such that
$\pos(t_A) =A$,
\item[(b)] if $\eta\in T^*_{K,\Sigma}$, $A\subseteq B\subseteq\pos(
t^*_\eta)$, then $t_A\in\Sigma(t_B)$ and $\nor[t_A]\leq\nor[t_B]$,
\item[(c)] if $r_\nu=r_\nu'+r_\nu''$, $r_\nu,r_\nu',r_\nu''\in [0,1]$ (for
$\nu\in\pos(t^*_\eta)$, $\eta\in T^*_{K,\Sigma}$), then
\[F_\eta(r_\nu:\nu\in\pos(t^*_\eta))=F_\eta(r_\nu':\nu\in\pos(t^*_\eta))+
F_\eta(r_\nu'':\nu\in\pos(t^*_\eta)).\]
\end{enumerate}
\end{definition}

\begin{lemma}
\label{mixlem}
Let $(K,\Sigma,\bF)$ be a ccc--complete mixing triple for $\bH$. Suppose
that $p_0,\ldots p_m\in\bQ^\mtree_\emptyset(K,\Sigma,\bF)$ are such that
$\sum\limits_{\ell\leq m}\mu^\bF(p_\ell)>1$ and $\mrot(p_0)=\ldots=
\mrot(p_m)$. Then for some $\ell<n\leq m$ the conditions $p_\ell,p_n$ are
compatible in $\bQ^\mtree_\emptyset(K,\Sigma,\bF)$. 
\end{lemma}

\begin{proof}
Let $\nu=\mrot(p_0)=\ldots=\mrot(p_m)$. For each $\ell<n\leq m$ such that
$[T^{p_\ell}]\cap [T^{p_n}]\neq\emptyset$ choose a tree $T_{\ell,n}\subseteq
T^*_{K,\Sigma}$ satisfying  
\[\max(T_{\ell,n})=\emptyset,\quad\mrot(T_{\ell,n})=\nu\quad \mbox{and}\quad
[T_{\ell,n}]=[T^{p_\ell}]\cap [T^{p_n}].\] 
Let $p_{\ell,n}\in\bQ^\tree_\emptyset(K,\Sigma)$ be such that $T^{p_{\ell,
n}}=T_{\ell,n}$ (defined if $[T^{p_\ell}]\cap [T^{p_n}]\neq\emptyset$,
$\ell<n\leq m$). Our aim is to show that for some $\ell<n\leq m$,
$p_{\ell,n}$ is defined and belongs to $\bQ^\mtree_\emptyset(K,\Sigma,\bF)$
(i.e., $\mu^\bF(p_{\ell,n})>0$). So suppose that for each $\ell<n\leq m$,
either $[T^{p_\ell}]\cap [T^{p_n}]=\emptyset$ or $\mu^\bF(p_{\ell,n})=0$. 
Let $\varepsilon=2^{-(m+1)}(\sum\limits_{\ell\leq m}\mu^\bF(p_\ell)-1)>0$
and let $p\in\bQ^\tree_\emptyset(K,\Sigma)$ be such that $\mrot(p)=\nu$ and
$T^p=T^{p_0}\cup\ldots\cup T^{p_m}$ (clearly $\min\{\mu^\bF(p_\ell):\ell\leq
m\}\leq\mu^\bF(p)\leq 1$). Choose a front $A$ of $T^p$ such that for each
$\ell<n\leq m$, if $p_{\ell,n}$ is defined and $A_{\ell,n}=A\cap
T^{p_{\ell,n}}$, then $\mu_{p_{\ell,n},A_{\ell,n}}(\nu)<\varepsilon$, and if
$[T^{p_\ell}]\cap [T^{p_n}]=\emptyset$ then $T^{p_\ell}\cap T^{p_n}\subseteq
T[p,A]\setminus A)$. 
\begin{claim}
\label{cl11}
For each $\eta\in T[p,A]$ we have
\begin{enumerate}
\item[$(\otimes)$] \quad $\mu_{p,A}(\eta)\geq\sum\limits_{\ell\leq m}
\mu_{p_\ell,A^\ell}(\eta)-\sum\limits_{\ell<n\leq m}\mu_{p_{\ell,n},
A_{\ell,n}}(\eta)$,
\end{enumerate}
where $A^\ell=A\cap T^{p_\ell}$ (for $\ell\leq m$), and if $\ell<n\leq m$
and $p_{\ell,n}$ is not defined or $\eta\notin T[p_{\ell,n},A_{\ell,n}]$
then we stipulate $\mu_{p_{\ell,n},A_{\ell,n}}(\eta)=0$ (and similarly
$\mu_{p_\ell,A^\ell}(\eta)=0$ if $\eta\notin T[p_\ell,A^\ell]$).
\end{claim}

\begin{proof}[Proof of the claim] 
We show this by downward induction on $\lh(\eta)$. 

First suppose that $\eta\in A$. Let $k=|\{\ell\leq m:\eta\in A^\ell\}|=
\sum\limits_{\ell\leq m}\mu_{p_\ell,A^\ell}(\eta)$. Then $\sum\limits_{\ell<
n\leq m}\mu_{p_{\ell,n},A_{\ell,n}}(\eta)=\binom{k}{2}$
and 
\[\sum_{\ell\leq m}\mu_{p_\ell,A^\ell}(\eta)-\sum_{\ell<n\leq m}\mu_{
p_{\ell,n},A_{\ell,n}}(\eta)=k-\frac{k(k-1)}{2}\leq 1=\mu_{p,A}(\eta).\] 

Suppose now that $(\otimes)$ has been shown for all $\rho\in\pos(t^p_\eta)$,
$\eta\in T[p,A]\setminus A$. Let 
\[X=\{\rho\in\pos(t^p_\eta): \sum_{\ell\leq m}\mu_{p_\ell,A^\ell}(\rho)\geq
\sum_{\ell<n\leq m}\mu_{p_{\ell,n},A_{\ell,n}}(\rho)
\}\]
and $Y=\pos(t^p_\eta)\setminus X$. It follows from the inductive hypothesis
and \ref{ccccomplete}(c) that
\[F_\eta(\mu_{p,A}(\rho):\rho\in X)\geq \sum_{\ell\leq m}F_\eta(\mu_{p_\ell,
A^\ell}(\rho):\rho\in X)-\!\sum_{\ell<n\leq m} F_\eta(\mu_{p_{\ell,
n},A_{\ell,n}}(\rho):\rho\in X).\]
[Note that though \ref{ccccomplete}(c) guarantees the additivity of $F_\eta$ 
only when $r_\nu=r_\nu'+r_\nu''$, $r_\nu,r_\nu',r_\nu''\in [0,1]$, we can
first prove that 
\[F_\eta(\frac{1}{M}\cdot r_\nu:\nu\in\pos(t^*_\eta))=\frac{1}{M}F_\eta(
r_\nu:\nu\in\pos(t^*_\nu)).\]
Next, we may reduce the needed additivity to the one postulated in
\ref{ccccomplete}(c) by dividing all terms by suitably large $M$.] Now, by
\ref{mixing}(e),  
\[\sum_{\ell\leq m}F_\eta(\mu_{p_\ell,A^\ell}(\rho):\rho\in Y)\leq
\sum_{\ell<n\leq m} F_\eta(\mu_{p_{\ell,n},A_{\ell,n}}
(\rho):\rho\in Y),\]
and hence
\[F_\eta(\mu_{p,A}(\rho):\rho\in Y)\geq \sum_{\ell\leq m}F_\eta(\mu_{p_\ell,
A^\ell}(\rho):\rho\in Y)-\sum_{\ell<n\leq m} F_\eta(\mu_{p_{\ell,n},
A_{\ell,n}}(\rho):\rho\in Y).\]
Since
\[\begin{array}{lcl}
\mu_{p,A}(\eta)&=&F_\eta(\mu_{p,A}(\rho):\rho\in X)+F_\eta(\mu_{p,A}(\rho):
\rho\in Y),\\
\mu_{p_\ell,A^\ell}(\eta)&=&F_\eta(\mu_{p_\ell,A^\ell}(\rho):\rho\in X)+
F_\eta(\mu_{p_\ell,A^\ell}(\rho):\rho\in Y),\\
\mu_{p_{\ell,n},A_{\ell,n}}(\eta)&=&F_\eta(\mu_{p_{\ell,n},A_{\ell,n}}(
\rho):\rho\in X)+F_\eta(\mu_{p_{\ell,n},A_{\ell,n}}(\rho):\rho\in Y),
  \end{array}\]
we may easily finish.
\end{proof}
Now we apply \ref{cl11} to $\nu=\mrot(p_0)$. We get then
\[1\geq \mu_{p,A}(\nu)\geq \sum_{\ell\leq m}\mu_{p_\ell,A^\ell}(\nu)-
\sum_{\ell<n\leq m}\mu_{p_{\ell,n},A_{\ell,n}}(\nu)>\sum_{\ell\leq m}
\mu^\bF(p_\ell)-2^{m+1}\cdot\varepsilon=1,\]
a contradiction.
\end{proof}

\begin{corollary}
\label{mixcor}
Let $(K,\Sigma,\bF)$ be a ccc--complete mixing triple for $\bH$.
\begin{enumerate}
\item The forcing notion $\bQ^\mtree_\emptyset(K,\Sigma,\bF)$ satisfies the
ccc. 
\item If $f:\omega\times\omega\longrightarrow\omega$ is a fast function and
$(K,\Sigma)$ is linked, then the forcing notions
$\bQ^\mtree_f(K,\Sigma,\bF)$ and $\bQ^\mtree_1(K,\Sigma,\bF)$ are ccc.
\item If $h:\omega\times\omega\longrightarrow\omega$ is regressive,
$(K,\Sigma)$ is $h$--linked and $\cF\subseteq (\omega\setminus
2)^{\textstyle \omega}$ is either countable or $\geq^*$--directed, then
$\bQ^\mtree_\cF(K,\Sigma,\bF)$ is ccc.
\end{enumerate}
\end{corollary}

\begin{remark}
\begin{enumerate}
\item Forcing notions determined by mixing triples are in some sense
mixtures of the random real forcing with forcings determined by
tree--creating pairs. The ``$\mtree$'' in $\bQ^\mtree_*(K,\Sigma,\bF)$
stands for ``{\bf m}easured {\bf t}ree''.
\item Because of \ref{mixcor}(1) (and the proof of \ref{mixlem}) we can be
very generous as far as the demands on the norms are concerned, and still we
may easily ensure that the resulting forcing notion satisfies the ccc. For
example, if $(K,\Sigma,\bF)$ is a ccc--complete mixing triple, $(K,\Sigma)$
is semi--linked in the sense that the demand of \ref{treelin}(1) is
satisfied whenever $\lh(\eta)$ is even, and 
\[\bQ^\mtree_{1/2}(K,\Sigma,\bF)=\{p\in\bQ^\mtree_\emptyset(K,\Sigma,\bF):
(\forall\eta\in [T^p])(\lim_{k\to\infty}\nor[t^p_{\eta\rest 2k}]=\infty)
\},\]
then $\bQ^\mtree_{1/2}(K,\Sigma,\bF)$ is ccc too.
\item This type of constructions (i.e., mixture--like) for not-ccc case will
be presented in \cite{RoSh:736} and \cite[\S 2]{RoSh:670}.
\end{enumerate}
\end{remark}

Let us finish this subsection with showing that the forcing notions
$\bQ^\mtree_\emptyset(K,\Sigma,\bF)$ tend to have many features of the
random real forcing. 

\begin{definition}
\label{measure}
Let $(K,\Sigma,\bF)$ be a mixing triple, $p\in\bQ^\mtree_\emptyset(K,\Sigma,
\bF)$.
\begin{enumerate}
\item A function $\mu:T^p\longrightarrow [0,1]$ is {\em a
semi--$\bF$--measure on $p$} if 
\[(\forall\eta\in T^p)\big(\mu(\eta)\leq F_\eta(\mu(\nu):\nu\in
\pos(t^p_\eta))\big).\]
\item If above the equality holds (for each $\eta\in T^p$), then $\mu$ is
called {\em an $\bF$--measure}.
\end{enumerate}
\end{definition}

\begin{proposition}
\label{usesemi}
Assume $(K,\Sigma,\bF)$ is a mixing triple, $p\in\bQ^\tree_\emptyset(K,
\Sigma)$.
\begin{enumerate}
\item If $\mu:T^p\longrightarrow [0,1]$ is semi--$\bF$--measure on $p$, then
for each $\eta\in T^p$ we have $\mu(\eta)\leq\mu^\bF(p^{[\eta]})$. 
\item The mapping $\eta\mapsto \mu^\bF(p^{[\eta]}):T^p\longrightarrow [0,1]$
is an $\bF$--measure on $p$.
\item If there is a semi--$\bF$--measure $\mu$ on $p$ such that $\mu(
\mrot(p))>0$, then $p\in\bQ^\mtree_\emptyset(K,\Sigma,\bF)$.
\end{enumerate}
\end{proposition}

\begin{proof}
Straightforward.
\end{proof}

\begin{proposition}
\label{mixcos}
Suppose that $(K,\Sigma,\bF)$ is a ccc--complete mixing triple, and $p_0,
\ldots,p_m\in\bQ^\mtree_\emptyset(K,\Sigma,\bF)$ are such that
$\mrot(p_0)=\ldots=\mrot(p_m)$. Let $p\in\bQ^\mtree_\emptyset(K,\Sigma,\bF)$
be such that $T^p=T^{p_0}\cup\ldots\cup T^{p_m}$. Then:
\begin{enumerate}
\item $\mu^\bF(p)\leq \sum\limits_{\ell\leq m}\mu^\bF(p_\ell)$,
\item if $[T^{p_\ell}]\cap [T^{p_n}]=\emptyset$ for $\ell<n\leq m$ (or just
$p_0,\ldots, p_m$ are pairwise incompatible in $\bQ^\mtree_\emptyset(K,
\Sigma,\bF)$), then $\mu^\bF(p)=\sum\limits_{\ell\leq m}\mu^\bF(p_\ell)$,
\item $\{p_0,\ldots,p_m\}$ is pre-dense above $p$.
\end{enumerate}
\end{proposition}

\begin{proof}
Like \ref{mixlem}.
\end{proof}

\begin{lemma}
\label{getstr}
Let $(K,\Sigma,\bF)$ be a ccc--complete mixing triple. 
\begin{enumerate}
\item Suppose that conditions $p,q\in\bQ^\mtree_\emptyset(K,\Sigma,\bF)$ are
such that $\mrot(p)=\mrot(q)$, $p\leq q$, $\mu^\bF(q)<\mu^\bF(p)$, and let
$0<\varepsilon<1$. Then there is a condition $r\in\bQ^\mtree_\emptyset(K,
\Sigma,\bF)$ stronger than $p$ and incompatible with $q$ and such that 
$\mrot(r)=\mrot(p)$, $\mu^\bF(r)\geq (1-\varepsilon)\cdot(\mu^\bF(p)-
\mu^\bF(q))$. 
\item If $p\in\bQ^\mtree_\emptyset(K,\Sigma,\bF)$, $\varepsilon>0$ then
there is $\eta\in T^p$ such that $\mu^\bF(p^{[\eta]})>(1-\varepsilon)$.
\end{enumerate}
\end{lemma}

\begin{proof}
1)\quad Let $\nu=\mrot(p)$. Choose a front $A$ of $T^p$ such that 
\[\mu_{q,A\cap T^q}(\nu)<\mu^\bF(q)+\varepsilon\cdot(\mu^\bF(p)-
\mu^\bF(q))<\mu^\bF(p)\]
(so necessarily $A\setminus T^q\neq\emptyset$). Take $r_0,r_1\in
\bQ^\tree_\emptyset(K,\Sigma)$ such that 
\[\begin{array}{lcl}
T^{r_0}&=&\{\eta\in T^p: (\exists\rho\in A\setminus T^q)(\eta\trianglelefteq
\rho\mbox{ or } \rho\vartriangleleft\eta)\}\\
T^{r_1}&=&\{\eta\in T^p: (\exists\rho\in A\cap T^q)(\eta\trianglelefteq
\rho\mbox{ or } \rho\vartriangleleft\eta)\}.
  \end{array}\]
Clearly $\mrot(r_i)=\nu$, $A\cap T^{r_i}$ is a front of $T^{r_i}$, 
$A=(A\cap T^{r_0}) \cup (A\cap T^{r_1})$ and $T^q\subseteq T^{r_1}$. Hence,
$r_1\in \bQ^\mtree_\emptyset(K,\Sigma,\bF)$ and $\mu^\bF(q)\leq\mu^\bF(r_1)<
\mu^\bF(q)+\varepsilon\cdot(\mu^\bF(p)-\mu^\bF(q))$. Now, using
\ref{ccccomplete}(c), we may conclude that $\mu^\bF(r_0)\geq (1-\varepsilon)
\cdot (\mu^\bF(p)-\mu^\bF(q))$, finishing the proof.
\medskip

\noindent 2)\quad Straightforward.
\end{proof}

\begin{proposition}
\label{semidec}
Assume that $(K,\Sigma,\bF)$ is a ccc--complete mixing triple, $p\in
\bQ^\mtree_\emptyset(K,\Sigma,\bF)$. Let $m<\omega$, $\varepsilon>0$ and let
$\dot{\tau}$ be a $\bQ^\mtree_\emptyset(K,\Sigma,\bF)$--name such that
$p\forces\dot{\tau}<m$. Then there are $X\subseteq m$ and conditions $q_\ell
\in\bQ^\mtree_\emptyset(K,\Sigma,\bF)$ (for $\ell\in X$) such that
\begin{enumerate}
\item[$(\alpha)$] $q_\ell\forces\dot{\tau}=\ell$,
\item[$(\beta)$]  $\mrot(q_\ell)=\mrot(p)$,
\item[$(\gamma)$] $\sum\limits_{\ell\in X}\mu^\bF(q_\ell)\geq
(1-\varepsilon) \mu^\bF(p)$.
\end{enumerate}
\end{proposition}

\begin{proof}
Let $\nu=\mrot(p)$. For each $\ell<m$ define $\mu_\ell:T^p\longrightarrow
[0,1]$ by:
\[\mu_\ell(\eta)=\sup\{\mu^\bF(q):q\in\bQ^\mtree_\emptyset(K,\Sigma,\bF)\
\&\ q\geq p^{[\eta]}\ \&\ \mrot(q)=\eta\ \&\ q\forces\dot{\tau}=\ell\}\] 
(with the convention that $\sup\emptyset=0$). It follows from
\ref{mixing}(f) that each $\mu_\ell$ is an $\bF$--measure on $p$.  

\begin{claim}
\label{cl12}
$\mu^\bF(p)=\sum\limits_{\ell<m}\mu_\ell(\nu)$.
\end{claim}

\begin{proof}[Proof of the claim] 
First note that, by a suitable modification of \ref{mixlem}, we have
$\mu^\bF(p)\geq\sum\limits_{\ell<m}\mu_\ell(\nu)$. So suppose that $\mu^\bF
(p)>\sum\limits_{\ell<m}\mu_\ell(\nu)$, and let $n\in\omega$ be such that
$(1-\frac{2}{n})\cdot\mu^\bF(p)>\sum\limits_{\ell<m}\mu_\ell(\nu)$. Let 
\[S=\{\eta\in T^p: (1-\frac{1}{n})\cdot\mu^\bF(p^{[\eta]})>
\sum\limits_{\ell<m}\mu_\ell(\eta)\},\]
and let $T\subseteq S$ be a tree such that $\mrot(T)=\nu$ and $\eta\in T$,
$\rho\in\suc_{T^p}(\eta)\cap S$ imply $\rho\in T$. Clearly $\max(T)=
\emptyset$, so we may choose $r\in\bQ^\tree_\emptyset(K,\Sigma)$ so that
$T=T^r$. If $\mu^\bF(r)>0$, then ($r\in\bQ^\mtree_\emptyset(K,\Sigma,\bF)$
and) we may choose a condition $q\in \bQ^\mtree_\emptyset(K,\Sigma,\bF)$
stronger than $r$ which decides the value of $\dot{\tau}$. By
\ref{getstr}(2) we find $\eta\in T^q$ such that $\mu^\bF(q^{[\eta]})>(1-
\frac{1}{n})\geq (1-\frac{1}{n})\mu^\bF(p^{[\eta]})$, what contradicts $T^q
\subseteq T^r\subseteq S$. Therefore $\mu^\bF(r)=0$, and like in
\ref{getstr}(1) we may build a condition $q\geq p$ and a front $A$ of $T^q$
such that  
\begin{itemize}
\item $\mrot(q)=\nu$, $A\cap T^r=\emptyset$,
\item $T^q=\{\eta\in T^p: (\exists\rho\in A)(\eta\trianglelefteq\rho\mbox{
or }\rho\vartriangleleft\eta)\}$, 
\item $(1-\frac{1}{n})\mu^\bF(p^{[\eta]})\leq\sum\limits_{\ell<m}\mu_\ell(
\eta)\leq\mu^\bF(p^{[\eta]})$ for each $\eta\in A$, and 
\item $\mu^\bF(q)>(1-\frac{1}{n})\mu^\bF(p)$. 
\end{itemize}
Now, we may easily conclude that $\sum\limits_{\ell<m}\mu_\ell(\nu)\geq
(1-\frac{2}{n})\cdot \mu^\bF(p)$, getting a contradiction.
\end{proof} 
The conclusion of the proposition follows immediately from \ref{cl12} and
the definition of $\mu_\ell$'s (so we take $X=\{\ell<m: \mu_\ell(\nu)>0\}$
and suitable $q_\ell$'s for $\ell\in X$). 
\end{proof}

\begin{proposition}
Suppose that $(K,\Sigma,\bF)$ is a ccc--complete mixing triple and
$\dot{\tau}$ is a $\bQ^\mtree_\emptyset(K,\Sigma,\bF)$--name for an ordinal.
Let $p\in\bQ^\mtree_\emptyset(K,\Sigma,\bF)$, $0<\varepsilon<1$. Then there
is a condition $q\geq p$ and a front $A$ of $T^q$ such that
\begin{enumerate}
\item[$(\alpha)$] $\mrot(q)=\mrot(p)$, $\mu^\bF(q)>(1-\varepsilon)
\mu^\bF(p)$, 
\item[$(\beta)$] for each $\eta\in A$ the condition $q^{[\eta]}$ decides the
value of $\dot{\tau}$.
\end{enumerate}
\end{proposition}

\begin{proof}
Let
\[\begin{array}{ll}
B=\{\eta\in T^p:&\mbox{for some }p^*\geq p\mbox{ we have: }\
\mrot(p^*)=\eta,\ \ \\
&\mu^\bF(p^*)\geq(1-\frac{\varepsilon}{2})\mu^\bF(p^{[\eta]}),\ \mbox{ and } 
p^*\mbox{ decides $\dot{\tau}$ on a front}\}. 
  \end{array}\]
It follows from \ref{getstr} that, if $q\in\bQ^\mtree_\emptyset(K,\Sigma,
\bF)$ is a condition stronger than $p$, then $T^q\cap B\neq\emptyset$. If
$\mrot(p)\in B$, then we are clearly done, so suppose $\mrot(p)\notin B$. 
Note that if $\eta\in T^p\setminus B$ then $\suc_{T^p}(\eta)\setminus
B\neq\emptyset$. Thus 
\[T\stackrel{\rm def}{=}\{\eta\in T^p: (\forall\nu\trianglelefteq\eta)(\nu
\notin B)\}\]
is a tree with $\max(T)=\emptyset$, $T\cap B=\emptyset$, $\mrot(T)=\mrot(p
)$. This $T$ determines a condition $r\in\bQ^\tree_\emptyset(K,\Sigma)$,
$\mrot(r)=\mrot(p)$. It follows from the previous remark that $\mu^\bF(r)
=0$. Take a front $A$ of $T^p$ such that $\mu_{r,A\cap T^r}(\mrot(r))<
\frac{1}{4}\cdot\varepsilon\cdot\mu^\bF(p)$ and $A\subseteq T^r\cup B$. For
each $\nu\in A\setminus T^r$ fix a condition $q_\nu$ such that
$p^{[\nu]}\leq q_\nu$, $\mu^\bF(q_\nu)\geq (1-\frac{\varepsilon}{2})\mu^\bF(
p^{[\nu]})$, and $q_\nu$ decides $\dot{\tau}$ on a front. Let $q$
be such that 
\[T^q=\{\eta\in T^p:(\exists \nu\in A\setminus T^r)(\eta\in T^{q_\nu}\mbox{
or }\nu\trianglelefteq\eta)\}.\] 
It should be clear that $q$ is a condition as required
\end{proof}

\subsection{Exotic norm conditions}
The norm conditions introduced in the first section have their counterparts
in the non-ccc case (as presented in \cite{RoSh:470}). Here we formulate
more norm conditions which may be used to build ccc forcing notions from
linked creating pairs (or tree--creating pairs), and which seem to be very
ccc--specific. Let us start with a norm condition that allows us to include
into our framework the ``Mathias with ultrafilter'' forcing notion.

\begin{definition}
\label{defult}
\begin{enumerate}
\item A local creating pair $(K,\Sigma)$ for $\bH$ is {\em strongly
linked\/} if it is full (see \ref{varia}(5)), linked and 
\begin{enumerate}
\item[$(\otimes)^{\rm sl}$] there are $t^{\min}_\ell\in K$ (for $\ell<
\omega$) such that $m^{t^{\min}_\ell}_\dn=\ell$ and if $t\in K$, $m^t_\dn
=\ell$, then $t^{\min}_\ell\in\Sigma(t)$.
\end{enumerate}
If additionally for each $\ell<\omega$ we have
\[(\forall u\in\dom(\val[t^{\min}_\ell]))(\exists !\,v\in\rng(\val[
t^{\min}_\ell]))(u\vartriangleleft v),\]
then we say that $(K,\Sigma)$ is {\em strongly$^+$ linked}.
\item A local tree--creating pair $(K,\Sigma)$ for $\bH$ is {\em (tree-)
strongly linked\/} if it is linked and 
\begin{enumerate}
\item[$(\otimes)^{\rm sl}_\tree$] there are $t^{\min}_\eta\in K\cap
\TCR_\eta[\bH]$ (for $\eta\in\bigcup\limits_{m\in\omega}\prod\limits_{i<m}
\bH(i)$) such that if $t\in K\cap\TCR_\eta[\bH]$, $\eta\in\bigcup\limits_{m
\in\omega}\prod\limits_{i<m}\bH(i)$, then $t^{\min}_\eta\in\Sigma(t)$.
\end{enumerate}
If, additionally, $|\pos(t^{\min}_\eta)|=1$ for each $\eta$ then we say that
$(K,\Sigma)$ is {\em strongly$^+$ linked}. 
\item Let $D$ be a non-principal ultrafilter on $\omega$. We define norm
conditions $\cC(D)$ and $\cC^\tree(D)$ (for the contexts of creating pairs
and tree--creating pairs, respectively) and the corresponding forcing
notions $\bQ^*_D(K,\Sigma)$, $\bQ^\tree_D(K,\Sigma)$ as follows.
\begin{itemize}
\item A sequence $\langle t_i: i<\omega\rangle$ satisfies $\cC(D)$ if
for some $\ell<\omega$ we have: 
\[(\forall i<\omega)(m^{t_i}_\dn=\ell+i)\quad\mbox{ and }\quad \lim_{D}
\langle\nor[t_{j-\ell}]: \ell\leq j<\omega\rangle=\infty.\] 
For a local creating pair $(K,\Sigma)$, $\bQ^*_D(K,\Sigma)$ is the forcing
notion 
\[\bQ^*_{\cC(D)}(K,\Sigma)=\{p\in\bQ^*_\emptyset(K,\Sigma):\langle t^p_i:
i<\omega\rangle\mbox{ satisfies }\cC(D)\}.\]  
\item A system $\langle t_\eta:\eta\in T\rangle\subseteq\TCR[\bH]$ satisfies
$\cC^\tree(D)$ if $T$ is a tree, $t_\eta\in\TCR_\eta[\bH]$ and
$\pos(t_\eta)=\suc_T(\eta)\subseteq\prod\limits_{i\leq\lh(\eta)}\bH(i)$ for
each $\eta\in T$, and $(\forall\eta\in [T])(\lim_D\langle\nor[t_{\eta\rest
k}]: \lh(\mrot(T))\leq k<\omega\rangle=\infty)$.

For a local tree creating pair $(K,\Sigma)$, $\bQ^\tree_D(K,\Sigma)$ is the
forcing notion 
\[\bQ^\tree_{\cC(D)}(K,\Sigma)=\{p\in\bQ^\tree_\emptyset(K,\Sigma):\langle
t^p_\eta: \eta\in T^p\rangle\mbox{ satisfies }\cC^\tree(D)\}.\] 
\end{itemize}
\end{enumerate}
\end{definition}

\begin{remark}
Strongly linked (and especially strongly$^+$ linked) creating pairs resemble
omittory creating pairs of \cite[Def.\ 2.1.1]{RoSh:470} -- in both cases the
practical examples are such that we may ``omit'' some of the creatures from
a condition $p$. Here the ``omitting'' is done by replacing $t^p_i$ by the
suitable $t^{\min}_\ell$ (see \ref{omitex}). 
\end{remark}

\begin{proposition}
Let $D$ be a non-principal ultrafilter on $\omega$. If $(K,\Sigma)$ is a
local strongly linked creating pair (tree--creating pair, respectively),
then the forcing notion $\bQ^*_D(K,\Sigma)$ ($\bQ^\tree_D(K,\Sigma)$,
respectively) is $\sigma$--centered.
\end{proposition}

\begin{proof}
Straightforward.
\end{proof}

Forcing notions $\bQ^*_D(K,\Sigma)$, though similar to the Mathias forcing
notion, do not have (in general) as many nice properties as this one. For
example ``deciding formulas without changing the root'' may easily fail,
even though we may have some kind of continuous reading of names.
 
\begin{definition}
[See {\cite[Def.~1.2.9]{RoSh:470}}]  
Let $D$ be a non-principal ultrafilter on $\omega$, $(K,\Sigma)$ be a
local strongly linked creating pair and $\dot{\tau}$ be a
$\bQ^*_D(K,\Sigma)$--name for an ordinal. We say that a condition
$p\in\bQ^*_D(K,\Sigma)$ {\em approximates $\dot{\tau}$ at $t^p_n$} (or
at $n$) whenever the following demand is satisfied:
\begin{enumerate}
\item[$(*)$] for each $w_1\in\pos(w^p,t^p_0,\ldots,t^p_{n-1})$, if there
is a condition $r\in\bQ^*_D(K,\Sigma)$ stronger than $p$ and such that
$w^r=w_1$ and $r$ decides the value of $\dot{\tau}$, then the condition
$(w_1,t^p_n, t^p_{n+1},\ldots)$ decides the value of $\dot{\tau}$ 
\end{enumerate}
\end{definition}

\begin{proposition}
Assume that $D$ is a Ramsey ultrafilter on $\omega$ and $(K,\Sigma)$ is a
local, really finitary and strongly$^+$ linked creating pair. Then for each
$p\in\bQ^*_D(K,\Sigma)$ and a name $\dot{\tau}$ for an ordinal, there is a
condition $q \geq p$ which approximates $\dot{\tau}$ at every $n$ and such
that $w^p=w^q$. 
\end{proposition}

\begin{proof}
The proof follows the lines of the appropriate proof for the Mathias forcing
notion (see e.g.~\cite[\S 7.4]{BaJu95}). Let $\langle t^{\min}_\ell:\ell<
\omega\rangle$ witness that $(K,\Sigma)$ satisfies $(\otimes)^{\rm sl}$ of
\ref{defult}(1),  $(\forall u\in\dom(\val[t^{\min}_\ell]))(\exists !v\in
\rng(\val[t^{\min}_\ell]))(u\vartriangleleft v)$. For simplicity, we assume
that $\nor[t^{\min}_\ell]\leq 1$ (for $\ell<\omega$). 

For a condition $p\in\bQ^*_D(K,\Sigma)$ and $n\in\omega$ let 
\[{\rm supp}^n(p)\stackrel{\rm def}{=}\{m^{t^p_i}_\dn: i<\omega\ \&\
\nor[t^p_i]>n+1\}\in D.\]

Choose inductively conditions $p_n\in\bQ^*_D(K,\Sigma)$ such that for each
$n<\omega$:
\begin{enumerate}
\item $p_0=p$, $p_n\leq p_{n+1}$, $w^{p_n}=w^p$, and $t^{p_{n+1}}_i=
t^{p_n}_i$ for $i<n$, 
\item if $w\in\pos(w^p,t^p_0,\ldots,t^p_{n-1})$ and there is $p^*\geq p$
such that $w^{p^*}=w$ and $p^*$ decides $\dot{\tau}$, then
$(w,t^{p_{n+1}}_n,t^{p_{n+1}}_{n+1},\ldots)$ decides $\dot{\tau}$.
\end{enumerate}
(Note that ``strongly linked'' implies that if $w^{q_0}=w^{q_1}$, then
$q_0,q_1$ are compatible; also remember that $(K,\Sigma)$ is full.) Since
$D$ is Ramsey, we may choose an increasing sequence $\langle i_n:n<\omega
\rangle\subseteq\omega\setminus\lh(w^p)$ such that 
\[\{i_n:n<\omega\}\in D\quad\mbox{ and }\quad (\forall n\in\omega)(i_n+2<
i_{n+1}\in {\rm supp}^{n+1}(p_{i_n-\lh(w^p)+2})).\]
For $j<\omega$ let 
\[t^q_j=\left\{\begin{array}{ll}
t^{p_{i_n-\lh(w^p)+2}}_j&\mbox{if }j=i_{n+1}-\lh(w^p),\ n\in\omega,\\
t^{\min}_j     &\mbox{if }j+\lh(w^p)\notin\{i_{n+1}:n<\omega\}.
	       \end{array}\right.\]  
It should be clear that $q\stackrel{\rm def}{=}(w^p,t^q_0,t^q_1,\ldots)\in
\bQ^*_D(K,\Sigma)$ is a condition stronger than $p$ and for every $w\in
\pos(w^p,t^p_0,\ldots,t^p_{i_n-\lh(w^p)})$, $n\in\omega$ we have
\[(w,t^q_{i_n-\lh(w^p)+1}, t^q_{i_n-\lh(w^p)+2},\ldots)\geq (w,t^{p_{i_n-
\lh(w^p)+2}}_{i_n-\lh(w^p)+1},t^{p_{i_n-\lh(w^p)+2}}_{i_n-\lh(w^p)+2},
\ldots).\]
Hence easily $q$ approximates $\dot{\tau}$ at all points of the form
$i_{n+1}-\lh(w^p)+1$ (for $n<\omega$), and by the additional demand on
$t^{\min}_\ell$ (in ``strongly$^+$'') we conclude that $q$ approximates
$\dot{\tau}$ at all $n<\omega$. 
\end{proof}

\begin{proposition}
\label{omitcoh}
Suppose that $(K,\Sigma)$ is a strongly$^+$ linked local creating pair (with
$t^{\min}_n$ witnessing this). Assume that $\nor[t^{\min}_n]\leq 1$ (for
$n<\omega$) and 
\begin{enumerate}
\item[$(*)$] for each $n\in\omega$ there are disjoint sets $A_n,B_n\subseteq
\bH(n)$ such that 
\begin{itemize}
\item if $\langle u,v\rangle\in\val[t^{\min}_n]$, then $v(n)\notin A_n\cup
B_n$,
\item if $t\in K$, $\nor[t]>1$, $u\in\dom(\val[t])$, $\lh(u)=n$

\noindent then there are $v_0,v_1$ such that $\langle u,v_0\rangle, \langle
u,v_1\rangle\in\val[t]$ and $v_1(n)\in A_n$ and $v_0(n)\in B_n$.
\end{itemize}
\end{enumerate}
Let $D$ be an ultrafilter on $\omega$. Then the forcing notion $\bQ^*_D(K,
\Sigma)$ adds a Cohen real.
\end{proposition}

\begin{proof}
Let $\dot{W}$ be the name for $\bQ^*_D(K,\Sigma)$--generic real and let
$\dot{K}=\{\dot{k}_n:n<\omega\}$ be a name for a subset of $\omega$ such
that  
\[\dot{K}=\{k\in\omega:\dot{W}(k)\in A_k\cup B_k\}.\] 
(Clearly $\dot{K}$ is infinite.) Let $\dot{c}\in\can$ be given by $\dot{c}
(n)=0$ if and only if $\dot{W}(\dot{k}_n)\in A_n$. It should be clear that
$\dot{c}$ is a name for a Cohen real over $\bV$.
\end{proof}

Now we will give some norm conditions that can be used in the context of
local and forgetful creating pairs. Note that if $(K,\Sigma)$ is of that
type, then for each $t\in K$ we have (unique) set $P_t\subseteq \bH(
m^t_\dn)$ such that 
\[\langle u,v\rangle\in \val[t]\quad\mbox{if and only if}\quad v(m^t_\dn)
\in P_t,\ u\vartriangleleft v\ \mbox{ and }\ \lh(v)=\lh(u)+1\]
(for some, equivalently all, $u\in\dom(\val[t])$). The set $P_t$ corresponds
to $\pos(t)$ in the tree--creatures context, and below we will use the
notation $\pos(t)$ for it (hoping that this does not cause any confusion). 
Our next definition is a variant of \ref{ccccomplete}(a,b) for the case of
local forgetful creating pairs.

\begin{definition}
\label{complete}
A local forgetful creating pair $(K,\Sigma)$ for $\bH$ is {\em complete\/}
if 
\begin{enumerate}
\item[(a)] for each $i\in\omega$ and a nonempty set $A\subseteq\bH(i)$,
there is a unique creature $t^i_A\in K$ such that $m^{t^i_A}_\dn=i$ and
$\pos(t^i_A)=A$, 
\item[(b)] if $i\in\omega$, $A\subseteq B\subseteq\bH(i)$, then $t^i_A\in
\Sigma(t^i_B)$ and $\nor[t^i_A]\leq\nor[t^i_B]$.
\end{enumerate}
\end{definition}

\begin{definition}
Let $\bH:\omega\longrightarrow\cHa$ be such that $(\forall n\in\omega)(|
\bH(n)|> 2^n)$ and let $(K,\Sigma)$ be a complete pair for $\bH$. 
\label{normK}
\begin{enumerate}
\item {\em A $1$--norming system\/} (for $\bH$) is a pair
$(\bar{K},\bar{g})$ such that
\begin{enumerate}
\item[$(\alpha)$] $\bar{K}=\langle K_\ell:\ell\in\omega\rangle$ is a
sequence of infinite pairwise disjoint subsets of $\omega$, $\min(K_\ell)
\geq\ell$,
\item[$(\beta)$]  $\bar{g}=\langle g_\rho:\rho\in\fs\rangle$, where for each
$\ell<\omega$:
\item[$(\gamma)$] if $\rho\in 2^{\textstyle\ell}$ then $g_\rho\in
\prod\limits_{m\in K_\ell}\bH(m)$, and
\item[$(\delta)$] for every $m\in K_\ell$, there are no repetitions in
$\langle g_\rho(m):\rho\in 2^{\textstyle \ell}\rangle$.
\end{enumerate}
\item Let $\cnor$ be a norm condition for $K$ and $(\bar{K},\bar{g})$ be a
$1$--norming system. We define {\em $(\bar{K},\bar{g})$--modified version
$\cnor^{\bar{K},\bar{g}}$ of $\cnor$} by

a sequence $\langle t_i:i<\omega\rangle$ satisfies
$\cnor^{\bar{K},\bar{g}}$\quad if and only if 

it satisfies $\cnor$ and for some $\rho_0,\ldots,\rho_k\in\can$, $k<\omega$
we have
\[(\forall i,\ell<\omega)(\forall \rho\in 2^{\textstyle \ell})([m^{t_i}_\dn
\in K_\ell\ \&\ g_\rho(m^{t_i}_\dn)\notin \pos(t_i)]\ \Rightarrow\ [\rho
\vartriangleleft\rho_0\vee\ldots\vee\rho\vartriangleleft\rho_k]).\]
\item If $\cnor$ is one of ${\mathcal C}^\infty(\nor)$, ${\mathcal
C}^\cF(\nor)$ or ${\mathcal C}^f(\nor)$ (for suitable $f$, $\cF$; see
\ref{nordef}), then the forcing notions corresponding to
$(\bar{K},\bar{g})$--modified versions of $\cnor$ will be denoted by
$\bQ^{\bar{K},\bar{g}}_\infty(K,\Sigma)$, $\bQ^{\bar{K},\bar{g}}_\cF(K,
\Sigma)$, $\bQ^{\bar{K},\bar{g}}_f(K,\Sigma)$, respectively. 
\end{enumerate}
\end{definition}

\begin{proposition}
\label{cccKg}
Let $\bH:\omega\longrightarrow\cHa$ be such that $(\forall n\in\omega)(|
\bH(n)|> 2^n)$. Assume that $(K,\Sigma)$ is a complete creating pair, and
$(\bar{K},\bar{g})$ is a $1$--norming system (for $\bH$).  
\begin{enumerate}
\item If $f:\omega\times\omega\longrightarrow\omega$ is fast and
$(K,\Sigma)$ is linked, then $\bQ^{\bar{K},\bar{g}}_\infty(K,\Sigma)$ and
$\bQ^{\bar{K},\bar{g}}_f(K,\Sigma)$ are $\sigma$-$*$--linked Souslin forcing
notions.  
\item Assume that $h:\omega\times\omega\longrightarrow\omega$ is regressive 
and $\cF\subseteq (\omega\setminus 2)^{\textstyle\omega}$ is an $h$--closed
family which is either countable or $\geq^*$--directed. Suppose $(K,\Sigma)$
is local and $h$--linked. Then the forcing notion $\bQ^{\bar{K},\bar{g}}_\cF
(K,\Sigma)$ is $\sigma$-$*$--linked, and if $\cF$ is countable and
$\geq^*$--directed, then $\bQ^{\bar{K},\bar{g}}_\cF(K,\Sigma)$ is also
Souslin.     
\end{enumerate}
\end{proposition}

\begin{proof}
Straightforward.
\end{proof}

\begin{definition}
Let $(K,\Sigma)$ be a local forgetful creating pair for $\bH$. 
\label{normU}
\begin{enumerate}
\item {\em A $2$--norming system\/} is a sequence $\bar{U}=\langle
U_{\rho,k}:\rho\in\fs\ \&\ k<\omega\rangle $ of pairwise disjoint infinite
subsets of $\omega$ such that $\lh(\rho)\leq\min(U_{\rho,k})$.
\item For a norm condition $\cnor$ and a $2$--norming system $\bar{U}$ we
define {\em $\bar{U}$--modified version $\cnor^{\bar{U}}$ of $\cnor$} by

a sequence $\langle t_i:i<\omega\rangle$ satisfies
$\cnor^{\bar{U}}$\quad if and only if 

it satisfies $\cnor$ and for some $\rho_0,\ldots,\rho_\ell\in\can$ and
$k_0,\ldots,k_\ell<\omega$, $\ell<\omega$, 

for every $i,k<\omega$ and $\rho\in\fs$ we have: 
\[[m^{t_i}_\dn\in U_{\rho,k}\ \&\ \pos(t_i)\neq\bH(m^{t_i}_\dn)]\
\Rightarrow\ [\rho \vartriangleleft\rho_0\vee\ldots\vee\rho\vartriangleleft
\rho_\ell \mbox{ and } k\in\{k_0,\ldots,k_\ell\}].\] 
We will use notation $\bQ^{\bar{U}}_\infty(K,\Sigma)$, $\bQ^{\bar{U}}_\cF(K,
\Sigma)$, $\bQ^{\bar{U}}_f(K,\Sigma)$ for the respective forcing notions
(and suitable $f$, $\cF$).
\end{enumerate}
\end{definition}

\begin{proposition}
\label{cccU}
Let $(K,\Sigma)$ be a complete creating pair, and $\bar{U}$ be a
$2$--norming system.   
\begin{enumerate}
\item If $f:\omega\times\omega\longrightarrow\omega$ is fast and
$(K,\Sigma)$ is linked, then $\bQ^{\bar{U}}_\infty(K,\Sigma)$ and
$\bQ^{\bar{U}}_f(K,\Sigma)$ are $\sigma$-$*$--linked Souslin forcing
notions.   
\item Assume that $h:\omega\times\omega\longrightarrow\omega$ is regressive 
and $\cF\subseteq (\omega\setminus 2)^{\textstyle\omega}$ is an $h$--closed
family which is either countable or $\geq^*$--directed. Suppose $(K,\Sigma)$
is local and $h$--linked. Then the forcing notion $\bQ^{\bar{U}}_\cF(K,
\Sigma)$ is $\sigma$-$*$--linked, and if $\cF$ is countable and
$\geq^*$--directed, then $\bQ^{\bar{U}}_\cF(K,\Sigma)$ is also Souslin.    
\end{enumerate}
\end{proposition}

\begin{proof}
Straightforward.
\end{proof}

\subsection{Universal forcing notions}
Here we introduce constructions involving very peculiar norm conditions. As
a matter of fact, norms are not important in that type of constructions, but
they still provide examples. Prototypes for the method described here are the
Universal Meager forcing notions (i.e., the Amoeba for Category) and forcing
notions related to variants of the PP--property (see \cite[Ch.VI,
\S2.12]{Sh:f}, \cite[Ch.7]{RoSh:470}).  

\begin{definition}
\label{finsys}
Let $(K,\Sigma)$ be a tree--creating pair for $\bH$. {\em A finite candidate
for $(K,\Sigma)$} is a system $\langle t_\eta:\eta\in \hat{S}\rangle$ such
that 
\begin{enumerate}
\item[(i)] $S\subseteq\bigcup\limits_{k\leq\lev(S)}\prod\limits_{i<k}\bH(i)$
is a tree of height $\lev(S)<\omega$, each node in $S$ has a successor at
the last level $\lev(S)$,
\item[(ii)] $\hat{S}=S\setminus\max(S)$ (i.e., non-maximal nodes of $S$),
\item[(iii)] $t_\eta\in\TCR_\eta[\bH]\cap K$ and $\pos(t_\eta)=\suc_S(\eta)$
(for $\eta\in\hat{S}$).
\end{enumerate}
The collection of finite candidates for $(K,\Sigma)$ is denoted
$\FC(K,\Sigma)$. It is equipped with the following order (similar to that of
$\bQ^\tree_\emptyset(K,\Sigma)$): 

$\langle t^0_\eta:\eta\in\hat{S}^0\rangle\leq \langle t^1_\eta:\eta\in
\hat{S}^1\rangle$\quad if and only if\quad $\lev(S^0)\leq\lev(S^1)$ and 
\[(\forall\eta\in S^1)(\lh(\eta)<\lev(S^0)\quad\Rightarrow\quad\eta\in S^0\
\&\ t^1_\eta\in\Sigma(t^0_\eta)).\]
\end{definition}

\begin{remark}
Finite candidates for tree--creating pairs correspond to that for creating
pairs (see \ref{candidates}). In general, finite candidates do not have to
be finite (just the respective tree is of finite height), but if
$(K,\Sigma)$ is finitary then they are.
\end{remark}

\begin{definition}
\label{univpar}
Let $\bH:\omega\longrightarrow\cHa$. {\em A universality parameter $\gp$ for
$\bH$} is a tuple $(K^\gp,\Sigma^\gp,\cF^\gp,\cG^\gp)=(K,\Sigma,\cF,\cG)$
such that
\begin{enumerate}
\item[$(\alpha)$] $(K,\Sigma)$ is a really finitary local tree--creating
pair for $\bH$,
\item[$(\beta)$] $\cF\subseteq\baire$ is either countable or
$\leq^*$--directed [note the direction of the inequality!],
\item[$(\gamma)$] elements of $\cG$ are quadruples $(\langle t_\eta:\eta\in
\hat{S}\rangle,n_\dn,n_\up,\bar{r})$ such that
\begin{itemize}
\item $\langle t_\eta:\eta\in\hat{S}\rangle\in\FC(K,\Sigma)$,
\item $n_\dn\leq n_\up\leq\lev(S)$,
\item $\bar{r}=\langle r_i:i\in\dom(\bar{r})\rangle$, $r_i\in\omega$,
$\dom(\bar{r})\subseteq [n_\dn,n_\up]$,
\end{itemize}
\item[$(\delta)$] {\bf if:}
\begin{itemize}
\item $(\langle t^0_\eta:\eta\in \hat{S}^0\rangle,n^0_\dn,n^0_\up,\bar{r}^0)
\in\cG$, 
\item $\bar{r}^1=\langle r^1_i:i\in\dom(\bar{r}^1)\rangle$,
$r^1_i\in\omega$, $\dom(\bar{r}^0)\subseteq\dom(\bar{r}^1)$, and $r^0_i\leq
r^1_i$ for $i\in\dom(\bar{r}^0)$, and 
\item $\langle t^1_\eta:\eta\in\hat{S}^1\rangle \in\FC(K,\Sigma)$, $\langle
t^0_\eta:\eta\in\hat{S}^0\rangle\leq\langle t^1_\eta:\eta\in\hat{S}^1
\rangle$ and  
\item $n^1_\dn\leq n^0_\dn$, $n^0_\up\leq n^1_\up\leq\lev(S^1)$ and
$\dom(\bar{r}^1)\subseteq [n^1_\dn,n^1_\up]$, 
\end{itemize}
{\bf then:}\qquad $(\langle t^1_\eta:\eta\in \hat{S}^1\rangle,n^1_\dn,
n^1_\up,\bar{r}^1)\in\cG$, 
\item[$(\varepsilon)$] for some increasing function $F=F^\cG\in\baire$, {\bf
if:} 
\begin{itemize}
\item $(\langle t^\ell_\eta:\eta\in \hat{S}^\ell\rangle,n^\ell_\dn,
n^\ell_\up,\bar{r}^\ell)\in\cG$ (for $\ell<2$), $\lev(S^0)=\lev(S^1)$,  
\item $\langle t_\eta:\eta\in\hat{S}\rangle\in\FC(K,\Sigma)$, $\langle
t_\eta:\eta\in\hat{S}\rangle\leq\langle t^\ell_\eta:\eta\in \hat{S}^\ell
\rangle$ (for $\ell<2$),   
\item $\lev(S)<n^0_\dn$, $n^0_\up<n^1_\dn$, $F(n^1_\up)<\lev(S^1)$,
\end{itemize}
{\bf then:}\qquad there is $(\langle t^*_\eta:\eta\in \hat{S}^*\rangle,
n^*_\dn,n^*_\up,\bar{r}^*)\in\cG$ such that
\begin{itemize}
\item $n^*_\dn=n^0_\dn$, $n^*_\up=n^1_\up$, $\bar{r}^*=\bar{r}^0\cup
\bar{r}^1$, $\lev(S^*)=\lev(S^0)=\lev(S^1)$,
\item $S\subseteq S^*$ and $t^*_\eta=t_\eta$ for $\eta\in\hat{S}$,
\item $\langle t^*_\eta:\eta\in \hat{S}^*\rangle\leq\langle t^\ell_\eta:
\eta\in \hat{S}^\ell\rangle$ (for $\ell<2$).
\end{itemize}
\end{enumerate}
\end{definition}

\begin{remark}
In Definition \ref{unforc} below, we may think about the forcing notion
$\bQp$ in the following way. We have a criterion for ``small trees''
provided by $\cG$ and $\cF$ (these are $(\cG,f)$--narrow trees, see
\ref{unforc}(c)). We try to add a small tree that will almost cover all
small trees from the ground model. So, naturally, a condition $p$ consists
of a small tree (it is the system $\langle t^p_\eta: \eta\in T^p\rangle$),
in which some finite part (``below $N^p$'') is declared to be fixed. Now,
when we extend the condition $p$, we cannot change the tree below the level
$N^p$, but above that level we may {\em increase\/} it. The function $f^p$
controls in some sense the ``smallness'' of the tree $T^p$. See more later.
\end{remark}

\begin{definition}
\label{unforc}
Let $\gp=(K,\Sigma,\cF,\cG)$ be a universality parameter for $\bH$. We
define a forcing notion $\bQp$ as follows.
\medskip

\noindent {\bf A condition} in $\bQp$ is a triple $p=(N^p,\langle t^p_\eta:
\eta\in T^p\rangle,f^p)$ such that
\begin{enumerate}
\item[(a)] $\langle t^p_\eta:\eta\in T^p\rangle\in\bQ^\tree_\emptyset(K,
\Sigma)$, $\mrot(T^p)=\langle\rangle$,
\item[(b)] $N^p\in\omega$, $f^p\in\cF$,
\item[(c)] the system $\langle t^p_\eta: \eta\in T^p\rangle$ is {\em $(\cG,
f^p)$--narrow,} which means: 

\noindent for infinitely many $n<\omega$, for some $(\langle t_\eta:\eta\in
\hat{S}\rangle,n_\dn,n_\up,\bar{r})\in\cG$ we have
\begin{itemize}
\item $n_\dn=n$, and $(\forall i\in\dom(\bar{r}))(r_i\leq f^p(i))$ and
\item if $\eta\in T^p$, $\lh(\eta)<\lev(S)$, then $\eta\in S$ and $t^p_\eta
\in\Sigma(t_\eta)$. 
\end{itemize}
\end{enumerate}
\noindent{\bf The relation $\leq$} on $\bQp$ is given by:\\
$(N^0,\langle t^0_\eta:\eta\in T^0\rangle,f^0)\leq (N^1,\langle t^1_\eta:
\eta\in T^1\rangle,f^1)$\quad if and only if
\begin{itemize}
\item $N^0\leq N^1$, $\langle t^0_\eta:\eta\in T^0\rangle\geq\langle
t^1_\eta:\eta\in T^1\rangle$ (in $\bQ^\tree_\emptyset(K,\Sigma)$), and
\item if $\eta\in T^1$, $\lh(\eta)<N^0$ then $\eta\in T^0$ and $t^0_\eta=
t^1_\eta$, and
\item $(\forall^\infty i\in\omega)(f^0(i)\leq f^1(i))$.
\end{itemize}
\end{definition}

\begin{proposition}
\label{unccc}
If $\gp=(K,\Sigma,\cF,\cG)$ is a universality parameter, then $\bQp$ is a
$\sigma$--centered forcing notion. If additionally $\cF$ is countable {\em
and\/} $\leq^*$--directed then $\bQp$ is Borel ccc.
\end{proposition}

\begin{proof}
It is easy to check that the relation $\leq$ of $\bQp$ is transitive (so
$\bQp$ is a forcing notion). Let us argue that it is $\sigma$--centered when
$\cF$ is countable (the case of $\leq^*$--directed $\cF$ can be handled
similarly). 

For $\langle t_\eta:\eta\in S\rangle\in \FC(K,\Sigma)$, $f\in\cF$ let
\[\begin{array}{ll}
Q^{\langle t_\eta:\eta\in S\rangle}_f=\{p\in\bQp:& N^p=\lev(S)\mbox{ and }
S\subseteq T^p\mbox{ and } f^p=f\mbox{ and}\\
&(\forall\eta\in T^p)(\lh(\eta)<N^p\ \Rightarrow\ \eta\in S\ \&\ t^p_\eta=
t_\eta)\}.
  \end{array}\]  

\begin{claim}
\label{cl18}
Each $Q^{\langle t_\eta:\eta\in S\rangle}_f$ is a directed subset of $\bQp$.
\end{claim}

\begin{proof}[Proof of the claim] 
Let $(N^\ell,\langle t^\ell_\eta:\eta\in T^\ell\rangle,f^\ell)\in Q^{\langle
t_\eta:\eta\in S\rangle}_f$ (for $\ell<2$). (Thus $N^\ell=\lev(S)$,
$f^\ell=f$.) 

Let $F^\cG\in\baire$ be the increasing function given by
\ref{univpar}($\varepsilon$). Pick a sequence 
\[\lev(S)<n^{0,0}_\dn<n^{0,0}_\up<n^{1,0}_\dn<n^{1,0}_\up<\ldots<n^{0,k}_\dn
<n^{0,k}_\up<n^{1,k}_\dn<n^{1,k}_\up<\ldots\]
such that $F^\cG(n^{1,k}_\up)+1<n^{0,k+1}_\dn$ and (for $\ell<2$ and $k\in
\omega$) 
\[(\langle t^\ell_\eta:\eta\in T^\ell\ \&\ \lh(\eta)<n^{\ell,k}_\up\rangle,
n^{\ell,k}_\dn,n^{\ell,k}_\up, f\rest [n^{\ell,k}_\dn, n^{\ell,k}_\up])\in
\cG\]
(possible by the definition of the forcing $\bQp$ and
\ref{univpar}($\delta$)). Now build inductively a system $\langle t^*_\eta:
\eta\in T^*\rangle\in \bQ^\tree_\emptyset(K,\Sigma)$ as follows.

We declare that $\mrot(T^*)=\mrot(S)=\langle\rangle$, and if $\eta\in T^*$,
$\lh(\eta)<\lev(S)$, then $t^*_\eta=t_\eta$ and $\suc_{T^*}(\eta)=
\pos(t^*_\eta)$. 

Suppose we have defined $T^*$ up to the level $n^{0,k}_\dn-1$, so we know
$t^*_\eta$ for $\lh(\eta)<n^{0,k}_\dn-1$. Let $S_k$ be the tree of height
$n^{0,k}_\dn-1$ built from these $t^*_\eta$ (so it is the respective
``initial part'' of our future $T^*$), and assume that $\langle t^*_\eta:
\eta\in \hat{S}_k\rangle\leq \langle t^\ell_\eta:\eta\in T^\ell\rangle$ (for
$\ell=0,1$). Apply \ref{univpar}($\varepsilon$) to get $\langle t^*_\eta:
\eta\in\hat{S}_{k+1}\rangle$ such that
\[(\langle t^*_\eta:\eta\in\hat{S}_{k+1}\rangle,n^{0,k}_\dn,n^{1,k}_\up,f
\rest [n^{0,k}_\dn,n^{1,k}_\up])\in\cG,\]
and $S_k\subseteq S_{k+1}$, $\lev(S_{k+1})=n^{0,k+1}_\dn-1$ and $\langle
t^*_\eta:\eta\in\hat{S}_{k+1}\rangle\leq \langle t^\ell_\eta:\eta\in T^\ell
\rangle$ (for $\ell<2$). We declare that $T^*$ up to the level
$n^{0,k+1}_\dn-1$ is $S_{k+1}$ (and the respective $t^*_\eta$ are as chosen
above). 

Plainly, $(\lev(S),\langle t^*_\eta:\eta\in T^*\rangle,f)\in 
Q^{\langle t_\eta:\eta\in S\rangle}_f$ is a condition stronger than both 
$(N^0,\langle t^0_\eta:\eta\in T^0\rangle,f^0)$ and $(N^1,\langle t^1_\eta:
\eta\in T^1\rangle,f^1)$. 
\end{proof}
The rest should be clear. 
\end{proof}

\subsection{Examples}

\begin{example}
\label{ranult}
Let $\bH(i)=\omega$. Suppose that $\bD,\bB,\bS$ are functions such that
\begin{itemize}
\item $\dom(\bD)\subseteq\fseo$, $\bD(\eta)$ is a non-principal ultrafilter
on $\omega$ (for $\eta\in\dom(\bD)$), 
\item $\dom(\bS)=\dom(\bB)=\fseo\setminus\bD(\eta)$ and for each $\eta\in
\dom(\bB)$: 
\[2\leq\bB(\eta)\in\omega\cup\{\omega\},\quad \bS(\eta)=\langle s^\eta_k:k
\in\bB(\eta)\rangle\subseteq (0,1),\ \mbox{ and }\ \sum\limits_{k\in
\bB(\eta)}s^\eta_k=1.\]  
\end{itemize}
We build a ccc--complete (see \ref{ccccomplete}) mixing triple
$(K_{\ref{ranult}},\Sigma_{\ref{ranult}},\bF_{\ref{ranult}})$ for $\bH$ (for
the parameters $\bB,\bD,\bS$). 
\end{example}

\begin{proof}[Construction]
Let $K_{\ref{ranult}}$ consist of all tree creatures $t\in\TCR[\bH]$ such
that 
\begin{itemize}
\item $\dis[t]=(n_t,\eta_t,A_t)$ for some $n_t\in\omega$, $\eta_t\in
\prod\limits_{i<n_t}\bH(i)$ and $A_t\subseteq \omega$ such that
\[A_t\in\bD(\eta_t)\quad \mbox{if }\eta_t\in\dom(\bD),\quad \mbox{ and
}\quad\emptyset\neq A_t\subseteq \bB(\eta_t)\quad \mbox{if }\eta_t\in
\dom(\bB),\]
\item $\nor[t]=n_t$,
\item $\val[t]=\{\langle\eta_t,\nu\rangle:\eta_t\vartriangleleft\nu\in
\prod\limits_{i\leq n_t}\bH(i)\ \&\ \nu(n_t)\in A_t\}$.
\end{itemize}
The operation $\Sigma_{\ref{ranult}}$ is natural:
\[\Sigma_{\ref{ranult}}(t)=\{s\in K_{\ref{ranult}}: n_s=n_t\ \&\ \eta_s=
\eta_t\ \&\ A_s\subseteq A_t\}.\]
For $\eta\in \prod\limits_{i<n}\bH(i)$ let $t^*_\eta\in K_{\ref{ranult}}\cap
\TCR_\eta[\bH]$ be such that 
\[A_{t^*_\eta}=\left\{\begin{array}{ll}
\omega&\mbox{if }\eta\in\dom(\bD),\\
\bB(\eta)&\mbox{if }\eta\in\dom(\bB),
		      \end{array}\right.\]
and for $\langle r_\nu:\nu\in\pos(t^*_\eta)\rangle\in [0,1]^{\textstyle
\pos(t^*_\eta)}$ let
\[F_\eta^{\ref{ranult}}(r_\nu:\nu\in\pos(t^*_\eta))=\left\{
\begin{array}{ll}
\sum\limits_{k\in\bB(\eta)}s^\eta_k\cdot r_{\eta\conc\langle k\rangle}&
\mbox{if }\eta\in\dom(\bB),\\
\lim_{\bD(\eta)} \langle r_{\eta\conc\langle k\rangle}:k\in\omega\rangle&
\mbox{if }\eta\in\dom(\bD).
\end{array}\right.\]
Let $\bF_{\ref{ranult}}=\langle F_\eta^{\ref{ranult}}:\eta\in
\bigcup\limits_{n\in\omega}\prod\limits_{i<n}\bH(i)\rangle$. Check that
$(K_{\ref{ranult}},\Sigma_{\ref{ranult}},\bF_{\ref{ranult}})$ is a
ccc--complete mixing triple for $\bH$. 
\end{proof}

\begin{conclusion}
Let $\bH,\bD,\bB$ and $\bS$ be as in \ref{ranult}. Then the forcing notion
$\bQ^\mtree_\emptyset(K_{\ref{ranult}},\Sigma_{\ref{ranult}},
\bF_{\ref{ranult}})$ (for the parameters $\bD,\bB,\bS$) is ccc (and
non-trivial). 
\end{conclusion}

\begin{remark}
If $\dom(\bD)=\fseo$ then $\bQ^\mtree_\emptyset(K_{\ref{ranult}},
\Sigma_{\ref{ranult}},\bF_{\ref{ranult}})$ is the ``Laver with
ultrafilters'' forcing notion.

If $\dom(\bD)=\emptyset$ and $\bB(\eta)$ is finite (for each $\eta$) then 
$\bQ^\mtree_\emptyset(K_{\ref{ranult}},\Sigma_{\ref{ranult}},
\bF_{\ref{ranult}})$ is the random real forcing (with weights determined by
$\bS$ in an obvious way). 

Between these two extremes we have cases of ``mixtures of random with
ultrafilters'' and our next observation applies to most of them. It could be
formulated with a larger generality (e.g.~regarding $\dom(\bD)$), what
should be clear after reading the proof.
\end{remark}

\begin{proposition}
Let $\bH(n)=\omega$ for $n\in\omega$, and let $X\subseteq\omega$ be an
infinite co-infinite set. Suppose that 
\begin{enumerate}
\item[(a)] $\bD$ is a function such that $\dom(\bD)=\{\eta\in\fseo:\lh(\eta)
\in X\}$ and $\bD(\eta)$ is a non-principal Ramsey ultrafilter on $\omega$
(for $\eta\in\dom(\bD)$),
\item[(b)] $\bar{n}=\langle n_\ell: \ell\in\omega\setminus X\rangle$ is a
sequence of integers, $n_\ell\geq 1$,
\item[(c)] $\bB,\bS$ are functions such that $\dom(\bB)=\dom(\bS)=\{\eta\in
\fseo:\lh(\eta)\notin X\}$, $\bB(\eta)=n_{\lh(\eta)}$ and $\bS(\eta)=\langle
\frac{1}{n_{\lh(\eta)}}:k<n_{\lh(\eta)}\rangle$.
\end{enumerate}
Let $(K_{\ref{ranult}},\Sigma_{\ref{ranult}},\bF_{\ref{ranult}})$ be the
mixing triple built for the parameters $\bD,\bB,\bS$ as in \ref{ranult}, and
let $f\in\baire$ be such that $f(n)\geq 2$ (for $n\in\omega$).

Then for every $\bQ^\mtree_\emptyset(K_{\ref{ranult}},\Sigma_{\ref{ranult}},
\bF_{\ref{ranult}})$--name $\dot{\tau}$ for a real in
$\prod\limits_{i\in\omega}f(i)$, there are an increasing sequence $\langle
m_j:j\in\omega\rangle\subseteq\omega$ and a function $g\in
\prod\limits_{i\in\omega}f(i)$ such that 
\[\forces_{\bQ^\mtree_\emptyset(K_{\ref{ranult}},\Sigma_{\ref{ranult}},
\bF_{\ref{ranult}})}\mbox{`` }(\forall^\infty j\in\omega)(\dot{\tau}\rest
[m_j,m_{j+1})\neq g\rest [m_j,m_{j+1}))\mbox{ ''.}\]
Hence, in particular, forcing with $\bQ^\mtree_\emptyset(K_{\ref{ranult}},
\Sigma_{\ref{ranult}}, \bF_{\ref{ranult}})$ does not add Cohen reals (but it
clearly adds a dominating real).
\end{proposition}

\begin{proof}
For notational convenience, let $(K_{\ref{ranult}},\Sigma_{\ref{ranult}},
\bF_{\ref{ranult}})=(K,\Sigma,\bF)$. 

Note that we may assume that $f(n)>2^{2^{n+2}}$ (as we may work with the
mapping $n\mapsto f\rest [k_n,k_{n+1})$ for some increasing $\langle k_n:
n<\omega\rangle$ instead). Let $X=\{x_m:m<\omega\}$ be the increasing
enumeration, and let $m_k$ be defined by: $m_0=0$ and $m_{k+1}=m_k+2^k\cdot(
1+\prod\limits_{\ell\in x_k\setminus X}n_\ell)^k$ (for $k\in\omega$). [Here
we keep the convention that if $x_k\setminus X=\emptyset$, then
$\prod\limits_{\ell\in x_k\setminus X}n_\ell=1$; or just assume that
$x_0>0$.] 

Let $\dot{\tau}$ be a $\bQ^\mtree_\emptyset(K,\Sigma,\bF)$--name for an
element of $\prod\limits_{i\in\omega}f(i)$, and let $p$ be a condition in
$\bQ^\mtree_\emptyset(K,\Sigma,\bF)$. We may assume that $\lh(\mrot(p))\in
X$, and just for simplicity let $\lh(\mrot(p))=x_0$.

We define inductively a tree $T\subseteq T^p$, mappings $Y:T\cap
\bigcup\limits_{i<\omega}\omega^{\textstyle x_i+1}\longrightarrow\iso$,
$\pi:T\cap\bigcup\limits_{i<\omega}\omega^{\textstyle x_i+1}\longrightarrow
\omega$, a function $g\in\prod\limits_{i\in\omega}f(i)$, and a system
$\langle q_\eta: \eta\in T\cap\bigcup\limits_{i<\omega}\omega^{\textstyle
x_i+1}\rangle$ of conditions from $\bQ^\mtree_\emptyset(K,\Sigma,\bF)$.

We declare $\mrot(T)=\mrot(p)$. Using \ref{semidec} and \ref{mixcos} we
choose an increasing sequence of conditions $\langle q_{\mrot(T)}^k:
k<\omega\rangle$ and values for $g(m_k)<f(m_k)$ (thus defining $g\rest\{m_k:
k\in\omega\}$) such that 
\begin{itemize}
\item $q_{\mrot(T)}^0\geq p$, $\mrot(q_{\mrot(T)}^k)=\mrot(T)$ and
$\mu^\bF(q_{\mrot(T)}^k)>\frac{3}{4}\mu^\bF(p)$,
\item $q_{\mrot(p)}^k\forces\dot{\tau}(m_k)\neq g(m_k)$.
\end{itemize}
Since $\bD(\mrot(T))$ is a Ramsey ultrafilter we may choose a set $\{a_k:
k<\omega\}\in\bD(\mrot(T))$ (the increasing enumeration) such that
\[\mrot(T)\conc \langle a_k\rangle\in T^{q^k_{\mrot(T)}}\quad\mbox{and}
\quad\mu^\bF((q^k_{\mrot(T)})^{[\mrot(T)\conc\langle a_k\rangle]})>
\frac{3}{4} \mu^\bF(p).\]
We declare that $\mrot(T)\conc a_k\in T$ for all $k<\omega$ and we let 
\[\pi(\mrot(T)\conc\langle a_k\rangle)=k,\quad q_{\mrot(T)\conc\langle
a_k\rangle}= (q^k_{\mrot(T)})^{[\mrot(T)\conc\langle a_k\rangle]}.\]
Next we choose pairwise disjoint sets $Y(\mrot(T)\conc\langle a_k\rangle)
\subseteq \omega\setminus\{m_\ell:\ell\in\omega\}$ such that $m_{k+1}\cap
Y(\mrot(T)\conc\langle a_k\rangle)=\emptyset$ and 
\[(\forall\ell<\omega)(|Y(\mrot(T)\conc\langle a_k\rangle)\cap (m_{\ell
+k+1},m_{\ell+k+2})|=1).\]
Suppose now that we have already defined $T\cap \omega^{\textstyle x_i+1}$,
$i<\omega$, together with $Y(\eta)\subseteq\omega\setminus\{m_k:k\in\omega
\}$, $\pi(\eta)\in\omega$, $q_\eta\in\bQ^\mtree_\emptyset(K,\Sigma,\bF)$
(for $\eta\in T\cap \omega^{\textstyle x_j+1}$, $j\leq i$) and 
\[g_i\stackrel{\rm def}{=}g\rest (\{m_k:k<\omega\}\cup\bigcup \{Y(\eta):
\eta\in T\cap \omega^{\textstyle x_j+1}, j<i\})\]
so that the following conditions are satisfied:
\begin{enumerate}
\item[$(\alpha)$] if $\nu\vartriangleleft\eta$, and $\pi(\nu),\pi(\eta)$ are 
defined, then $\pi(\nu)<\pi(\eta)$,
\item[$(\beta)$]  if $\nu\in T\cap\omega^{\textstyle x_j}$, $j\leq i$, then
there are no repetitions in the sequence $\langle \pi(\eta):\nu
\vartriangleleft\eta\in T\cap\omega^{\textstyle x_j+1}\rangle$, 
\item[$(\gamma)$] if $\eta\in T\cap\omega^{\textstyle x_j+1}$, $j\leq i$,
then 
\[Y_\eta\cap m_{\pi(\eta)+1}=\emptyset\quad\mbox{and}\quad (\forall\ell<
\omega)(|Y_\eta\cap [m_{\pi(\eta)+\ell+1},m_{\pi(\eta)+\ell+2})|=1),\] 
\item[$(\delta)$] the (defined) $Y(\eta)$'s are pairwise disjoint,
\item[$(\varepsilon)$] for each $\eta\in T\cap\omega^{\textstyle x_i+1}$ we
have: $\mrot(q_\eta)=\eta$ and
\[q_\eta\forces (\forall k\leq\pi(\eta))(g_i\rest [m_k,m_{k+1})
\nsubseteq \dot{\tau}),\]
\item[$(\zeta)$] if $q\in\bQ^\mtree_\emptyset(K,\Sigma,\bF)$ is a condition
such that $\mrot(q)=\mrot(T)$, $T^q\cap \omega^{\textstyle x_i+1}=T\cap
\omega^{\textstyle x_i+1}$ and $q^{[\eta]}=q_\eta$ for all $\eta\in\omega^{
\textstyle x_i+1}$, then $\mu^\bF(q)\geq(\frac{1}{2}+\frac{1}{2^{i+2}})
\mu^\bF(p)$.
\end{enumerate}
[Check that these conditions are satisfied at the first stage of the
construction.] 

\noindent Note that it follows from clauses $(\alpha)$ and $(\beta)$ that
for each $k$ we have
\[|\{\eta\in T\cap\omega^{\textstyle \leq x_i+1}:\pi(\eta)\leq k\}|\leq 2^k
(1+\prod_{\ell\in x_i\setminus X}n_\ell)^k,\]
and hence (by clause $(\gamma)$)
\begin{enumerate}
\item[$(\boxtimes_k)$]\quad $|[m_{k+1},m_{k+2})\cap \bigcup\{Y(\eta):\eta\in
T\cap\omega^{\textstyle \leq x_i+1}\}|\leq 2^k (1+\prod\limits_{\ell\in
x_i\setminus X}n_\ell)^k$.
\end{enumerate}
Fix $\eta\in T\cap\omega^{\textstyle x_i+1}$ (and note that $\pi(\eta)\geq
i$). Choose an increasing sequence $\langle q^k_\eta:\pi(\eta)<k<\omega
\rangle$ of conditions in $\bQ^\mtree_\emptyset(K,\Sigma,\bF)$ and values
for $g\rest Y_\eta$ such that 
\begin{itemize}
\item $q_\eta^k\geq q_\eta$, $\mrot(q_\eta^k)=\eta$ and $\mu^\bF(q_\eta^k)>
(1-2^{-(i+4)})\mu^\bF(q_\eta)$,
\item if $j\in Y_\eta\cap [m_k,m_{k+1})$, $\pi(\eta)<k<\omega$, then 
$q_\eta^k\forces\dot{\tau}(j)\neq g(j)$,
\item the sequence $\langle T^{q_\eta^k}\cap \omega^{\textstyle x_{i+1}}: 
\pi(\eta)<k<\omega\rangle$ is constant (and let $\{\nu_\ell: \ell<\ell^*\}$
be the enumeration of $T^{q^k_2}\cap\omega^{\textstyle x_i+1}$; necessarily
$\ell^*\leq\prod\limits_{k=x_i+1}^{x_{i+1}-1}n_k$),
\item $r_\ell\stackrel{\rm def}{=}\lim\limits_{k\to\infty}\mu^\bF((
q_\eta^k)^{[\nu_\ell]})>0$ for each $\ell<\ell^*$ (note that the
limit exists as the sequence is non-increasing). 
\end{itemize}
[Why possible? Use \ref{semidec} and \ref{mixcos}, remember our additional
assumption on $f$.] Choose $a_k^\ell\in\omega$ (for $\ell<\ell^*$ and $k>
\pi(\eta)$) so that: 
\begin{itemize}
\item $\{a_k^\ell:\pi(\eta)<k<\omega\}\in \bD(\nu_\ell)$,
\item $\nu_\ell\conc\langle a^\ell_k\rangle\in T^{q^k_\eta}$ and 
$\mu^\bF((q^k_\eta)^{[\eta\conc\langle a_k^\ell\rangle]})>
(1-2^{-(i+4)})r_\ell$
\end{itemize}
(remember that each $\bD(\nu_\ell)$ is a Ramsey ultrafilter). Now we declare
that $\nu_\ell\conc a_k^\ell\in T$ for all $\ell<\ell^*$ and $\pi(\eta)<k
<\omega$ and we let 
\[\pi(\nu_\ell\conc\langle a_k^\ell\rangle)=k,\quad q_{\nu_\ell\conc\langle
a_k^\ell\rangle}= (q^k_\eta)^{[\nu_\ell\conc\langle a_k^\ell\rangle]}.\]
This finishes the definitions of $T\cap\omega^{\textstyle \leq x_{i+1}\!+\!
1}$ and of $\pi(\nu),q_\nu$ for $\nu\in T\cap\omega^{\textstyle x_{i+1}\!+\!
1}$. It should be clear that (the respective variants of) clauses
$(\alpha)$, $(\beta)$, $(\varepsilon)$ and $(\zeta)$ are satisfied. Using
$(\boxtimes_k)$ we may easily choose sets $Y(\nu)$ (for $\nu\in T\cap
\omega^{\textstyle x_{i+1}\!+\!1}$) so that the demands $(\gamma),(\delta)$
hold. The construction is finished.

The tree $T$ is perfect and it determines a condition $q\in
\bQ^\tree_\emptyset(K,\Sigma)$. 
\begin{claim}
\label{cl13}
$\mu^\bF(q)>0$ (so $q\in\bQ^\mtree_\emptyset(K,\Sigma,\bF)$).
\end{claim}

\begin{proof}[Proof of the claim] 
First note that the clause $(\zeta)$ is not enough to show this, as there
are fronts in $T^q$ which are not included in any $T^q\cap\omega^{\textstyle
\leq x_i+1}$. However, we may use the specific way the construction was
carried out to build a semi--$\bF$--measure $\mu:T^q\longrightarrow [0,1]$
such that $\mu(\mrot(q))=\frac{3}{8}\mu^\bF(p)>0$ (what is enough by
\ref{usesemi}). So, if $\eta\in T\cap\omega^{\textstyle x_i+1}$, $i<\omega$, 
then we let $\mu(\eta)=(1-\frac{1}{2^{i+1}})\mu^\bF(q_\eta)$; if $\eta\in
T$, $x_i+1<\lh(\eta)\leq x_{i+1}$ then $\mu(\eta)=F_\eta(\mu(\nu):\nu\in
\suc_T(\eta))$. Now check. 
\end{proof}
Thus $q$ is a condition stronger than $p$ and it forces that $(\forall
k\in\omega)(\dot{\tau}\rest [m_k,m_{k+1})\neq g\rest [m_k,m_{k+1}))$. Since
$\bQ^\mtree_\emptyset(K,\Sigma,\bF)$ satisfies the ccc, we may easily
finish.
\end{proof}

Our next example is a small modification of \ref{loctree}. In a similar way
we may modify other examples from the previous section to produce more
strongly linked creating pairs.

\begin{example}
\label{omitex}
Let $\bH\in\baire$ be a strictly increasing function such that $\bH(0)>2$.
We construct a really finitary, strongly$^+$ linked creating pair
$(K_{\ref{omitex}},\Sigma_{\ref{omitex}})$ for $\bH$ which is satisfies the
demands of \ref{omitcoh} (in particular $(*)$ there).
\end{example}

\begin{proof}[Construction]
The family $K_{\ref{omitex}}$ consists of creatures $t\in\CR[\bH]$
such that
\begin{itemize}
\item $\dis[t]=(m_t,E_t)$ such that $m_t<\omega$ and $\emptyset\neq E_t
\subseteq\bH(m_t)\setminus\{0\}$,
\item $\val[t]=\{\langle u,v\rangle\in\prod\limits_{i<m_t}\bH(i)\times
\prod\limits_{i\leq m_t}\bH(i): u\vartriangleleft v\ \&\ v(m_t)\notin
E_t\}$, 
\item $\nor[t]=\log_4(\frac{\bH(m_t)}{|E_t|})$.
\end{itemize}
The operation $\Sigma_{\ref{omitex}}$ is natural: 
\[\Sigma_{\ref{omitex}}(t)=\{s\in K_{\ref{omitex}}: m_s=m_t\ \&\
E_t\subseteq E_s\}.\]
Now check.
\end{proof}
 
\begin{conclusion}
Suppose that $\bH\in\baire$ is strictly increasing, $\bH(0)>2$ and
$(K_{\ref{omitex}},\Sigma_{\ref{omitex}})$ is built as in \ref{omitex} for
$\bH$. Let $D$ be a Ramsey ultrafilter on $\omega$. Then the forcing notion
$\bQ^*_D(K_{\ref{omitex}},\Sigma_{\ref{omitex}})$ is $\sigma$--centered,
adds a Cohen real and adds a dominating real.
\end{conclusion}

Now we turn to universality parameters. As said before, one of the
prototypes here is the Universal Meager forcing notion. Let us represent it
as $\bQp$ (for a suitable $\gp$).

\begin{example}
\label{UM}
We construct a universality parameter $\gp$ such that $\bQp$ is the
Universal Meager forcing notion. 
\end{example}

\begin{proof}[Construction]
Let $\bH:\omega\longrightarrow\omega\setminus 2$ and let $K$ consists of
tree creatures $t$ for $\bH$ such that
\begin{itemize}
\item $\dis[t]=(m_t,\eta_t,A_t)$ for some $m_t<\omega$, $\eta_t\in
\bigcup\limits_{i<m_t}\bH(i)$ and $\emptyset\neq A_t\subseteq\bH(m_t)$,
\item $\nor[t]=|A_t|$,
\item $\val[t]=\{\langle\eta_t,\eta_t\conc\langle a\rangle\rangle:a\in
A_t\}$.
\end{itemize}
The operation $\Sigma$ is natural, so $s\in\Sigma(t)$ if and only if
$\eta_s=\eta_t$ and $A_s\subseteq A_t$. Let $\cF=\{f\}$, $f(i)=0$. 

$\cG$ consists of quadruples $(\langle t_\eta:\eta\in
\hat{S}\rangle,n_\dn,n_\up,\bar{r})$ such that
\begin{itemize}
\item $\langle t_\eta:\eta\in\hat{S}\rangle\in\FC(K,\Sigma)$,
\item $n_\dn\leq n_\up\leq\lev(S)$,
\item $\bar{r}=\langle r_i:i\in\dom(\bar{r})\rangle$, $r_i<\omega$,
$\dom(\bar{r})\subseteq [n_\dn,n_\up]$,
\item if $\eta\in S$, $\lh(\eta)=n_\dn$ then for some $\nu\in S$ we have
$\eta\trianglelefteq\nu$, $\lh(\nu)<n_\up$ and $\nor[t_\nu]<\bH(\lh(\nu))$. 
\end{itemize}
Easily $\gp=(K,\Sigma,\cF,\cG)$ is a universality parameter and $\bQp$ is
the Universal Meager forcing notion,
\end{proof}

\begin{remark}
Our next example \ref{PPex} captures a number of constructions related to
the PP--property. Under the assumptions on $(K,\Sigma)$ as there, we may
think that we have a way to measure how large splittings are, and this fully  
determines what are the tree-creatures in $K$ (and what are the norms). The
function $\bF$ is used to define (possibly totally not related) norms of
sets of nodes of the same length. Thus $\bF$ may just count how many
elements are in $\max(S)$ (in this case the universality parameter given by
\ref{PPex} is related to the PP--property). Other possibilities for $\bF$
include taking the maximum value of $\nor[t_\eta]$, or taking the product of
all relevant $\nor[t_\eta]$'s. 
\end{remark}

\begin{example}
\label{PPex}
Assume that $\bH:\omega\longrightarrow\omega\setminus 2$ is strictly
increasing and a family $\cF\subseteq\baire$ is either countable or
$\leq^*$--directed [note the direction of the inequality]. Let $(K,\Sigma)$
be a local, really finitary, tree--creating pair for $\bH$ such that 
\begin{itemize}
\item for each $\eta\in\prod\limits_{i<n}\bH(i)$, $n<\omega$ and a nonempty
$A\subseteq\bH(n)$ there is a unique $t_{\eta,A}\in\TCR_\eta[\bH]\cap K$
with $\pos(t_{\eta,A})=\{\eta\conc \langle a\rangle:a\in A\}$, and
\item if $|A|=1$, then $\nor[t_{\eta,A}]\leq 1$, and
\item if $\emptyset\neq B\subseteq A\subseteq\bH(n)$, then $t_{\eta,B}\in
\Sigma(t_{\eta,A})$ and $\nor[t_{\eta,B}]\leq \nor[t_{\eta,A}]$. 
\end{itemize}
Furthermore, let $\bF:\FC(K,\Sigma)\longrightarrow\mbR^{{\geq}0}$ be
such that 
\begin{itemize}
\item if $\langle t_\eta:\eta\in\hat{S}\rangle\in\FC(K,\Sigma)$,
$\nor[t_\eta]\leq 1$ (for $\eta\in\hat{S}$), then $\bF(\langle t_\eta:\eta
\in\hat{S}\rangle)\leq 1$, and
\item if $\langle t^\ell_\eta:\eta\in\hat{S}\rangle\in FC(K,\Sigma)$ (for
$\ell<2$), and $\langle t^0_\eta:\eta\in\hat{S}\rangle\leq \langle t^1_\eta:
\eta\in\hat{S}\rangle$, then $\bF(\langle t^1_\eta:\eta\in\hat{S}\rangle)
\leq\bF(\langle t^0_\eta:\eta\in\hat{S}\rangle)$. 
\end{itemize}
We construct $\cG=\cG^{K,\Sigma}_\bF$ such that $(K,\Sigma,\cF,\cG)$ is a
universality parameter.
\end{example}

\begin{proof}[Construction]
For a nonempty set $Y\subseteq\prod\limits_{j<i}\bH(j)$, $i<\omega$, we
define ${\rm NOR}(Y)$ as the value of $\bF(\langle t_\eta:\eta\in\hat{S}
\rangle)$ for $S$ such that $\max(S)=Y$ and $\mrot(S)=\langle\rangle$. (Note
that under our assumptions on $(K,\Sigma)$ there is exactly one such
$\langle t_\eta:\eta\in\hat{S}\rangle \in\FC(K,\Sigma)$.) 

Let $\cG$ consists of all quadruples $(\langle t_\eta:\eta\in\hat{S}\rangle,
n_\dn,n_\up,\bar{r})$ satisfying the demands of \ref{univpar}($\gamma$) and
such that 

for some sequence $\langle Y_i:i\in\dom(\bar{r})\rangle$ we have
\begin{itemize}
\item $Y_i\subseteq\prod\limits_{j<i}\bH(j)$, ${\rm NOR}(Y_i)\leq r_i$,
\item $(\forall\eta\in S\cap\prod\limits_{j<n_\up}\bH(j))(\exists i\in\dom(
\bar{r}))(\eta\rest i\in Y_i)$. 
\end{itemize}
Now check.
\end{proof}

\begin{example}
\label{ucmz}
A universality parameter $\gp$ such that $\bQp$ is the ``universal closed
measure zero'' forcing notion.
\end{example}
 
\begin{proof}[Construction]
Let $\bH$, $(K,\Sigma)$ and $\cF$ be as defined in the construction for
\ref{UM}.

Let $\cG$ consists of all quadruples $(\langle t_\eta:\eta\in\hat{S}\rangle,
n_\dn,n_\up,\bar{r})$ satisfying the demands of \ref{univpar}($\gamma$) and
such that 
\[\frac{|S\cap\prod\limits_{i<n_\up}\bH(i)|}{|\prod\limits_{i<n_\up}\bH(i)|}
\leq \sum_{i\in\dom(\bar{r})}\frac{1}{(i+1)^2}.\]
Let $\gp=(K,\Sigma,\cF,\cG)$. Note that the forcing notion $\bQp$ is
equivalent to $\bQ$ defined as follows.

\noindent{\bf A condition} in $\bQ$ is a pair $(N,T)$ such that $N<\omega$
and $T\subseteq\bigcup\limits_{n<\omega}\prod\limits_{i<n}\bH(i)$ is a tree
such that $[T]$ is a measure zero subset of $\prod\limits_{i<\omega}
\bH(i)$;

\noindent{\bf the order} of $\bQ$ is the natural one: $(N_0,T_0)\leq
(N_1,T_1)$ if and only if $N_0\leq N_1$, $T_0\subseteq T_1$ and $T_1\cap
\prod\limits_{i<N_0}\bH(i)\subseteq T_0$.
\end{proof}

\section{Interlude: ideals} 
Here we introduce $\sigma$--ideals determined by forcing notions discussed
in this paper. Most of the content of this part is well known and belongs to
folklore (some of this material is presented in Judah and Ros{\l}anowski
\cite{JuRo92}, \cite{JuRo}).   

\subsection{Generic ideals}
We will show how a Souslin ccc forcing notion adding one real produces a
ccc Borel $\sigma$--ideal on some Polish space. While we could do this in a
larger generality (e.g., considering any name for a real, not only the ones
of the form specified in \ref{cont}(3) below, compare \cite[\S4]{Sh:666} and
\cite[\S6, \S7]{Sh:630}), we have decided to use the specific form of the 
forcing notions we want to deal with and simplify the notation and arguments
(loosing slightly on generality, but it will be clear how possible
generalizations go). 

\begin{context}
\label{cont}
\begin{enumerate}
\item $\bH:\omega\longrightarrow\cHa$, $|\bH(n)|\geq 2$ for all $n\in
\omega$; $\cT_n=\prod\limits_{i<n}\bH(i)$ (for $n\in\omega$) and
$\cT=\bigcup\limits_{n<\omega}\cT_n$. Let $\cX=\prod\limits_{n\in\omega}
\bH(n)=[\cT]$ be equipped with the natural product (Polish) topology.  
\item $\bP$ is a Souslin ccc forcing notion with a parameter $r\in\can$
(which also encodes $\bH$), so we have $\Sigma^1_1$--formulas
$\varphi_0(\cdot,r),\varphi_1(\cdot,\cdot,r),\varphi_2(\cdot,\cdot,r)$
defining $\bP,\leq_\bP$ and $\incomp_\bP$, respectively.
\item $\dot{W}=\langle p_\eta:\eta\in\cT\rangle\subseteq\bP$ is such that
\begin{enumerate}
\item[$(\alpha)$] $p_{\langle\rangle}=\emptyset_{\bP}$, and if $\eta
\vartriangleleft\nu\in\cT$ then $p_\eta\leq p_\nu$, 
\item[$(\beta)$] $\langle p_{\eta\conc\langle a\rangle}:a\in\bH(\lh(\eta))
\rangle$ is a maximal antichain above $p_\eta$,
\item[$(\gamma)$] for each $p\in\bP$ there is $n<\omega$ such that
\[|\{\eta\in\cT_n:p,p_\eta\mbox{ are compatible }\}|\geq 2,\]
\item[$(\delta)$] if $p,q\in\bP$ are incompatible, then there is $\eta\in
\cT$ such that $p,p_\eta$ are compatible but $q,p_\eta$ are incompatible.
\end{enumerate}
[We will treat $\dot{W}$ as a $\bP$--name for a real in $\cX$ such that
$p_\eta\forces\eta\vartriangleleft\dot{W}$, and so
$\forces\dot{W}\notin\bV$.]  
\end{enumerate}
\end{context}

\begin{definition}
\label{ideal}
For $\bP,\dot{W},\bH,\cX$ as in \ref{cont}, let $\Ipw$ be the collection of
all Borel subsets $B$ of $\cX$ such that
\[\forces_{\bP}\mbox{`` }\dot{W}\notin B\mbox{ ''.}\]
\end{definition}

\begin{proposition}
\label{basideal}
\begin{enumerate}
\item $\Ipw$ is a ccc $\sigma$--ideal of Borel subsets of $\cX$.
\item $\Ipw$ contains all singletons.
\item Let $\psi(p,c)$ be the formula 
\[\mbox{`` }p\in\bP\mbox{ and }c\in\can\mbox{ is a Borel code (for a set
$\sharp c\subseteq\cX$) and }p\forces\dot{W}\in\sharp c\mbox{ ''.}\]
\begin{enumerate}
\item[(a)] If $M$ is a transitive model of ${\rm ZFC}^*$, $p,r,c,\dot{W}\in
M$, then   
\[\psi(p,c)\quad\Leftrightarrow\quad M\models \psi(p,c).\]
\item[(b)] There are a $\Sigma^1_2$--formula $\psi_0$ and a
$\Pi^1_2$--formula $\psi_1$ (both with the parameter $r$) such that 
\[\psi(p,c)\equiv\psi_0(p,c)\equiv\psi_1(p,c)\]
(i.e., the equivalences are provable in {\rm ZFC}).
\end{enumerate}
\item The formula ``$c\in\can$ is a Borel code (for a set $\sharp
c\subseteq\cX$) and $\sharp c\in\Ipw$'' is absolute between transitive
models of ${\rm ZFC}^*$ (containing $r,\dot{W},c$). 
\end{enumerate}
\end{proposition}

\begin{proof}
(1), (2)\quad  Straightforward.

\noindent (3) See (the proof of) \cite[Lemma 3.6.12]{BaJu95}.

\noindent (4) Follows from (3) and the definition of $\Ipw$ (remember that
``being a maximal antichain of $\bP$'' is absolute; see
\cite[3.6.4]{BaJu95}). 
\end{proof}

\begin{definition}
\label{generic}
Let $\bP,\dot{W},\bH,\cX$ be as in \ref{cont}.
\begin{enumerate}
\item Let $x\in\cX$ and $M$ be a transitive model of ${\rm ZFC}^*$,
$r,\dot{W}\in M$. We say that $x$ is {\em $\Ipw$--generic over $M$}, if
$x\notin B$ for every Borel set $B$ from $\Ipw$ coded in $M$. 
\item For a condition $p\in\bP$ let 
\[S(p)=S_{\dot{W}}(p)\stackrel{\rm def}{=}\{\eta\in\cT: p_\eta,p\mbox{ are 
compatible }\}.\] 
\item For a maximal antichain $\cA\subseteq\bP$ let
\[B_\cA=B_{\dot{W},\cA}\stackrel{\rm def}{=}\{x\in\cX:(\forall p\in\cA)( 
\exists n<\omega)(x\rest n\notin S(p))\}.\]
Let $\Ipwz$ be the family of all subsets of $\cX$ that can be covered by a
set of the form $B_\cA$ (for a maximal antichain $\cA\subseteq\bP$).
\end{enumerate}
\end{definition}

\begin{proposition}
\label{element}
\begin{enumerate}
\item[(a)] For each $p\in\bP$, $S(p)$ is a perfect subtree of $\cT$. If
$p\leq q$, then $S(q)\subseteq S(p)$. If $p,q\in\bP$ and $S(q)\subseteq
S(p)$, then $q\forces p\in \Gamma_{\bP}$.  
\item[(b)] If $G\subseteq\bP$ is a generic filter over $\bV$, then $\bV[G]= 
\bV[\dot{W}^G]$. 
\item[(c)] $\Ipwz$ is an ideal of subsets of $\cX$; sets $B_\cA$ (for a
maximal antichain $\cA\subseteq\bP$) are $\Pi^0_2$.
\item[(d)] Let $x\in\cX$ and $M$ be a transitive model of ${\rm ZFC}^*$, $r,
\dot{W}\in M$. Then $x$ is $\Ipw$--generic over $M$ if and only if there is
a $\bP^M$--generic filter $G\subseteq\bP^M$ over $M$ such that
$\dot{W}^G=x$. 
\item[(e)] $\Ipw$ is the $\sigma$--ideal generated by $\Ipwz$. Every set
from $\Ipw$ can be covered by a $\Sigma^0_3$ set from $\Ipw$.
\end{enumerate}
\end{proposition}

\begin{proof}
Straightforward (or see \cite[\S 2]{JuRo92}).
\end{proof}

\begin{conclusion}
\label{quocon}
The quotient algebra $\borel(\cX)/\Ipw$ is a ccc complete Boolean
algebra. The mapping 
\[\pi:\bP\longrightarrow\borel(\cX)/\Ipw:p\mapsto[S(p)]_{\Ipw}\]
satisfies:
\begin{enumerate}
\item $\rng(\pi)$ is a dense subset of the algebra $\borel(\cX)/\Ipw$,
\item $(\forall p,q\in\bP)(p\incomp q\ \Leftrightarrow\ \pi(p)\cap\pi(q)={
\bf 0})$,
\item $(\forall p,q\in\bP)(q\forces p\in\Gamma_\bP\ \Leftrightarrow
\pi(q)\subseteq \pi(p))$.
\end{enumerate}
Consequently, ${\rm RO}(\bP)\cong\borel(\cX)/\Ipw$. Moreover, $\pi$ maps
$\dot{W}$ onto the canonical name for the generic real in
$\borel(\cX)/\Ipw$, so for a Borel code $c$ we have $\lbv\dot{W}\in\sharp
c\rbv_{{\rm RO}(\bP)}= [\sharp c]_{\Ipw}$. 
\end{conclusion}
  
\begin{remark}
\label{namrem}
It follows from \ref{quocon} that we have nice description of names for
reals in the extensions via $\bP$. 
\begin{enumerate}
\item If $\dot{\tau}$ is a $\bP$--name for an element of $\cX$, then there
is a Borel function $f:\cX\longrightarrow\cX$ such that $\forces_\bP
f(\dot{W})=\dot{\tau}$. 
\item If $\dot{B}$ is a $\bP$--name for a Borel subset of $\cX$, then there 
is a Borel set $A\subseteq\cX\times\cX$ such that $\forces_\bP\dot{B}=
(A)_{\dot{W}}$, where $(A)_x=\{y:(x,y)\in A\}$. 
\end{enumerate}
(See \cite[Lemma 3.7.1]{BaJu95}.)
\end{remark}

\subsection{Universality ideals}
For a forcing notion $\bQp$ (where $\gp$ is a universality parameter) we may
consider the ccc ideal defined as in \ref{ideal}, however there is another
Borel $\sigma$--ideal related to $\bQp$ (justifying the term ``universality
forcing notion'').

\begin{definition}
\label{unideal}
Let $\gp=(K,\Sigma,\cF,\cG)$ be a universality parameter for $\bH$.
\begin{enumerate}
\item We say that {\em $\gp$ is suitable} whenever: 
\begin{enumerate}
\item[(a)] for every $f\in\cF$ and $n<\omega$ there is $N>n$ such that 

\noindent {\bf if} $(\langle t_\eta:\eta\in\hat{S}\rangle,n_\dn,n_\up,f
\rest [n_\dn,n_\up])\in\cG$, $N\leq n_\dn$ and $\eta\in\prod\limits_{i<n}
\bH(i)$,  

\noindent {\bf then} $(\exists\nu\in\prod\limits_{i<\lev(S)}\bH(i))(\eta
\vartriangleleft\nu\ \&\ \nu\notin S)$,  
\item[(b)] for every $f\in\cF$ and  $n<\omega$ there is $N>n$ such that 

\noindent {\bf if} $\langle t_\eta:\eta\in\hat{S}\rangle\in\FC(K,\Sigma)$,
$\lev(S)=n$, $\eta\in\prod\limits_{i<N}\bH(i)$ and $\eta\rest n\in S$,

\noindent {\bf then} there is $(\langle t^*_\eta:\eta\in\hat{S}^*\rangle,
n_\dn,n_\up,f\rest [n_\dn,n_\up])\in\cG$ such that $n<n_\dn\leq n_\up<N$,
$S\subseteq S^*$ and $t_\eta=t^*_\eta$ for $\eta\in S$ and $\mrot(S)=
\mrot(S^*)$ and $\eta\in S^*$. 
\end{enumerate}
\item $\cI^0_\gp$ is the collection of subsets $A$ of
$\prod\limits_{i<\omega}\bH(i)$ such that for some $f\in\cF$ and a $(\cG,
f)$--narrow system $\langle t_\eta:\eta\in T\rangle\in\bQ^\tree_\emptyset
(K,\Sigma)$ (see \ref{unforc}(c)) we have $A\subseteq [T]$.\\
Trees $T$ as above will be called $(\cG,f)$--narrow.
\item $\cI_\gp$ is the $\sigma$--ideal of subsets of
$\prod\limits_{i<\omega} \bH(i)$ generated by $\cI^0_\gp$.
\item $\dot{T}_\gp$ is a $\bQp$--name such that 
\[\forces_{\bQp}\dot{T}_\gp=\bigcup\{T^p\cap\prod_{i<N^p}\bH(i): p\in
\Gamma_\bQp\}.\]
\end{enumerate}
\end{definition}

\begin{proposition}
Let $\gp=(K,\Sigma,\cF,\cG)$ be a suitable universality parameter for $\bH$.
\begin{enumerate}
\item Every set in $\cI^0_\gp$ is nowhere dense (in the product topology of 
$\prod\limits_{i<\omega} \bH(i)$); all singletons belong to $\cI^0_\gp$.
\item $\cI_\gp$ is a proper Borel $\sigma$--ideal of subsets of
$\prod\limits_{i<\omega} \bH(i)$. 
\item If $\cF$ is $\leq^*$--directed then $\cI^0_\gp$ is an ideal.
\item In $\bV^\bQp$, $\dot{T}_\gp$ is a tree with no maximal branches which
is $(\cG,h)$--narrow for some function $h$ (possibly $h\notin\cF$). If $\cF$
is a singleton, then $h\in\cF$.
\item Suppose that $\langle\bP_\alpha,\dbQ_\alpha:\alpha<\delta\rangle$ is
finite support iteration of ccc forcing notions such that for some
increasing sequence $\alpha_0<\alpha_1<\alpha_2<\ldots<\delta$,
$\dbQ_{\alpha_n}$ is (forced to be) $\bQp$. Let $\dot{T}_n$ be the name for
the tree $\dot{T}_\gp$ added at stage $\alpha_n$. Then, in
$\bV^{\bP_\delta}$, if $T\in\bV$ is $(\cG,f)$--narrow for some $f\in\cF$ then
$T\subseteq \dot{T}_n$ for some $n<\omega$.
\item If in 5 above we additionally assume that $\cF=\{f\}$ then, in
$\bV^{\bP_\delta}$, there is a set from $\cI_\gp$ which contains all Borel
sets from $\cI_\gp$ coded in $\bV$. 
\end{enumerate}
\end{proposition}

\subsection{Baire Property}
Let us discuss the Baire property related to our ideals and remind how we
may build models in which all projective sets of reals are regular in this
respect.

\begin{definition}
Let $\cI$ be a $\sigma$--ideal of subsets of a Polish space $\cX$. A set
$A\subseteq\cX$ is said to have $\cI$--Baire property if for some Borel set
$B\in\borel(\cX)$ the symmetric difference $A\vartriangle B$ is in $\cI$.
\end{definition}

Sets with $\cI$--Baire property constitute a $\sigma$--field of subsets of
$\cX$. Clearly, every Borel set  has $\cI$--Baire property. But even more,
using the method of Category Base (see Morgan \cite{Mo77}, \cite{Mo90}) we
may show that, for a number of $\sigma$--ideals $\cI$ on $\cX$, the family
of sets with $\cI$--Baire property is closed under Souslin operation $\cA$
(so it includes all $\Sigma^1_1$ subsets of $\cX$; this applies to, e.g.,
all ccc Borel ideals on $\cX$). However, it may happen that there are
$\Sigma^1_2$ subsets of $\cX$ that do not have the $\cI$--Baire
property. Below we present a tool that can be used to build a model in
which all projective subsets of $\cX$ have the $\Ipw$--Baire property for
some pairs $(\bP,\dot{W})$ as in \ref{cont}. So let us work in the context
of \ref{cont}.   

Let $\BB$ be a complete ccc Boolean algebra. Every $\BB$--name $\dot{\tau}$
for a real in $\cX$ determines a complete subalgebra $\BB_{\dot{\tau}}$ of
$\BB$. This subalgebra is generated by the family $\{\lbv
\eta\vartriangleleft\dot{\tau}\rbv_{\BB}: \eta\in\cT\}$ (as a complete
subalgebra). If $\dot{\tau}$ is a name for $\Ipw$--generic real such that
$\lbv\dot{\tau}\in B\rbv_{\BB}\neq {\bf 0}$ for every Borel set $B\in\borel
(\cX)\setminus\Ipw$, then $\BB_{\dot{\tau}}$ is a copy of the algebra ${\rm
RO}(\bP)$ (and the isomorphism maps $\dot{W}$ to the name $\dot{\tau}$).

\begin{definition}
\label{homo}
Let $\BB$ be a complete ccc Boolean algebra, and let $\BB'\lesdot\BB$ be its 
complete subalgebra. We say that $\BB$ is {\em $(\bP,\dot{W})$--homogeneous
over $\BB'$} if 
\begin{enumerate}
\item[$(\oplus)_1$] for each complete embeddings ${\bf e}_1,{\bf e}_2: \BB'*
{\rm RO}(\bP)\longrightarrow\BB$ such that ${\bf e}_1\rest\BB'={\bf e}_2
\rest\BB'={\rm id}_{\BB'}$ there is an automorphism $\pi$ of $\BB$ such that
${\bf e}_2=\pi\comp{\bf e}_1$, and
\item[$(\oplus)_2$] if ${\bf e}:\BB'*{\rm RO}(\bP)\longrightarrow\BB$ is a
complete embedding, $b\in\BB$ and $\pi_0$ is an automorphism of the complete
subalgebra of $\BB$ generated by ${\bf e}[\BB'*{\rm RO}(\bP)]\cup\{b\}$,
$\pi_0\rest {\bf e}[\BB'*{\rm RO}(\bP)]={\rm id}$, then there is an
automorphism $\pi$ of $\BB$ extending $\pi_0$.
\end{enumerate}
\end{definition}

\begin {theorem}[Solovay]
\label{Solovay}
Assume that $\BB$ is a complete ccc Boolean algebra, $\BB'$ is its complete
subalgebra, and $\dot{t}$ is a $\BB'$--name for a real. Suppose that 
\begin{enumerate}
\item $\forces_{\BB}$`` the union of all Borel sets from the ideal $\Ipw$
coded in $\bV^{\BB'}$ belongs to the ideal $\Ipw$ '',
\item $\BB$ is $(\bP,\dot{W})$--homogeneous over $\BB'$,
\item if $\dot{\tau}$ is a $\BB$--name for a real in $\cX$,
$\forces_{\BB}$`` $\dot{\tau}$ is $\Ipw$--generic over $\bV^{\BB'}$ '',\\
then $\dot{\tau}$ is a $\BB''$--name for some $\BB''$ such that
$\BB'\lesdot\BB'' \lesdot\BB$ and $\BB''*{\rm RO}(\bP)$ embeds in $\BB$ (by
an embedding that is identity on $\BB''$).   
\end{enumerate}
Then $\BB\forces$`` any subset of $\cX$ definable with $\dot{t}$ and with
(real) parameters from $\bV$ has the $\Ipw$--Baire property ''.
\end{theorem}

\begin{proof}
Suppose that $G\subseteq\BB$ is a generic filter over $\bV$ and $x\in\cX
\cap\bV[G]$ is $\Ipw$--generic over $\bV[G\cap\BB']$. It follows from the
third assumption that there is a $\BB$--name $\dot{\tau}$ for $x$ such that 
$\lbv\dot{\tau}\in\dot{A}\rbv_{\BB}\neq{\bf 0}$ for every $\BB'$--name
$\dot{A}$ for an element of $\borel(\cX)\setminus \Ipw$. Consequently, we
find a complete embedding ${\bf e}:\BB'*{\rm RO}(\bP)\longrightarrow\BB$
such that ${\bf e}\rest\BB'={\rm id}_{\BB'}$ and ${\bf
e}(\dot{W})=\dot{\tau}$.   

Now continue like in the proof of \cite[Theorem 2.3]{JuRo}; see also the
original Solovay's article \cite{So70}.
\end{proof}

\section{Sweet and Sour}

\subsection{On sweetness}
Solovay's Theorem \ref{Solovay} shows that to build a model with ($\Ipw$--)
Baire property of projective reals we may construct a suitably homogeneous   
complete Boolean algebra adding enough generic reals. The main tool for
building homogeneous Boolean algebras is amalgamation. 

\begin{definition}
\begin{enumerate}
\item Suppose $\bP,\bQ$ are forcing notions such that $\bP\lesdot{\rm
RO}(\bQ)$. Then $(\bQ:\bP)$ is a $\bP$--name for a forcing notion which is a
suborder of $\bQ$,
\[\begin{array}{l}
p\forces_{\bP}\mbox{`` }q\in(\bQ:\bP)\mbox{ ''\qquad\qquad if and only if}\\
\mbox{every $p'\in\bP$ stronger than $p$ is compatible with $q$ in $\bQ$.} 
  \end{array}\]
\item Suppose that $\bP,\bQ_0,\bQ_1$ are forcing notions and $f_\ell:\bP
\longrightarrow{\rm RO}(\bQ_\ell)$ (for $\ell<2$) are complete
embeddings. {\em The amalgamation of $\bQ_0,\bQ_1$ over $f_0,f_1$\/} is  
\[\begin{array}{l}
\bQ_0\times_{f_0,f_1}\bQ_1=\\
\{(q_1,q_2)\in\bQ_0\times\bQ_1: (\exists p\in\bP)(p\forces\mbox{`` }q_0 \in
(\bQ_0:f_0[\bP])\ \&\ q_1\in (\bQ_1:f_1[\bP])\mbox{''})\}
  \end{array}\]  
ordered in the natural way (so $(q_0,q_1)\leq(q_0',q_1')$ if and only if
$q_0 \leq q_0'$ and $q_1\leq q_1'$).  
\end{enumerate}
\end{definition}
On how the repeated use of amalgamations produces homogeneous Boolean
algebras see \cite{JuRo}. The main problem with amalgamating is that quite
often it produces forcing notions that collapse $\omega_1$ (this effect is
related to sourness discussed later). While the ccc might not be preserved
in amalgamations, some strong variants of it are.

\begin{definition}
[Shelah {\cite[Def.~7.2]{Sh:176}}; see Judah and Shelah
{\cite[Def.~2.5]{JdSh:446}} too] 
\label{sweet}
\ \\
A triple $(\bP,\cD,\bar{E})$ is {\em model of sweetness\/} (on $\bP$)
whenever: 
\begin{enumerate}
\item[(i)]    $\bP$ is a forcing notion, $\cD$ is a dense subset of $\bP$,
\item[(ii)]   $\bar{E}=\langle E_n:n<\omega\rangle$, each $E_n$ is an
equivalence relation on $\cD$ such that $\cD/E_n$ is countable,
\item[(iii)]  equivalence classes of each $E_n$ are directed,
$E_{n+1}\subseteq E_n$,
\item[(iv)]   if $\{p_i:i\leq\omega\}\subseteq \cD$, $p_i\; E_i\; p_\omega$
(for $i\in\omega$) then
\[(\forall n\in\omega)(\exists q\geq p_\omega)(q\; E_n\; p_\omega\ \&\
(\forall i\geq n)(p_i\leq q)),\]
\item[(v)]    if $p,q\in\cD$, $p\leq q$ and $n\in\omega$ then there is
$k\in\omega$ such that
\[(\forall p'\in [p]_{E_k})(\exists q'\in [q]_{E_n})(p'\leq q').\]
\end{enumerate}
If there is a model of sweetness on $\bP$, then we say that $\bP$ is {\em
sweet}. 
\end{definition}

\begin{definition}
[Stern {\cite[Def.~1.2]{St85}}]
\label{topsweet}
Let $\bP$ be a forcing notion and $\tau$ be a topology on $\bP$. We say that
$(\bP,\tau)$ is {\em a model of topological sweetness\/} whenever the
following conditions are satisfied:
\begin{enumerate}
\item[(i)]   the topology $\tau$ has a countable basis,
\item[(ii)]  $\emptyset_\bP$ is an isolated point in $\tau$,
\item[(iii)] if a sequence $\langle p_n:n<\omega\rangle\subseteq\bP$
is $\tau$--converging to $p\in\bP$, $q\geq p$ and $W$ is a
$\tau$--neighbourhood of $q$, then there is a condition $r\in\bP$ such that
\begin{enumerate}
\item[(a)] $r\in W$, $r\geq q$,
\item[(b)] the set $\{n\in\omega: p_n\leq r\}$ is infinite.
\end{enumerate}
\end{enumerate}
A forcing notion $\bP$ is {\em topologically sweet}, if there is a topology
$\tau$ on $\bP$ such that $(\bP,\tau)$ is a model of topological sweetness.
\end{definition}

\begin{definition}
\begin{enumerate}
\item Let $(\bP_\ell,\cD_\ell,\bar{E}_\ell)$ (for $\ell<2$) be models of
sweetness. We say that $(\bP_1,\cD_1,\bar{E}_1)$ {\em extends the model
$(\bP_0,\cD_0,\bar{E}_0)$} if
\begin{itemize}
\item $\bP_0\lesdot\bP_1$, $\cD_0\subseteq\cD_1$ and $E_n^0=E_n^1
\rest\cD_0$ for each $n\in\omega$, 
\item if $p\in\cD_0$, $n\in\omega$, then $[p]_{E_n^1}\subseteq\cD_0$, 
\item if $p\leq q$, $p\in\cD_1,q\in\cD_0$, then $p\in\cD_0$.
\end{itemize}
\item Let $(\bP_\ell,\tau_\ell)$ (for $\ell<2$) be models of topological
sweetness. We say that $(\bP_1,\tau_1)$ {\em extends the model
$(\bP_0,\tau_0)$} if $\bP_0\lesdot\bP_1$, $\bP_0$ is a $\tau_1$--open
subset of $\bP_1$ and $\tau_1\rest\bP_0=\tau_0$.
\end{enumerate}
\end{definition}

\begin{theorem}
\label{sweetamal}
\begin{enumerate}
\item {\em (Shelah \cite[7.5]{Sh:176})} Suppose that $(\bP_\ell,\cD_\ell,
\bar{E}_\ell)$ (for $\ell<2$) are models of sweetness, and $f_\ell:\bP
\longrightarrow{\rm RO}(\bP_\ell)$ (for $\ell<2$) are complete embeddings. 
Then there is a model of sweetness $(\bP_0\times_{f_0,f_1}\bP_1,\cD^*,
\bar{E}^*)$ based on the amalgamation $\bP_0\times_{f_0,f_1}\bP_1$ and
extending each  $(\bP_\ell,\cD_\ell,\bar{E}_\ell)$ for $\ell<2$.
\item {\em (Stern \cite[\S 2.2]{St85})} Suppose that $(\bP_\ell,\tau_\ell)$
(for $\ell<2$) are models of topological sweetness, and $f_\ell:\bP
\longrightarrow{\rm RO}(\bP_\ell)$ (for $\ell<2$) are complete embeddings. 
Then there is a model of topological sweetness $(\bP_0\times_{f_0,f_1}\bP_1,
\tau^*)$ based on the amalgamation $\bP_0\times_{f_0,f_1}\bP_1$ and
extending each $(\bP_\ell,\tau_\ell)$ for $\ell<2$.
\end{enumerate}
\end{theorem}

\begin{definition}
\label{itersweet}
Let $\cB$ be a countable basis of a topology on a forcing notion $\bP$. We
say that $(\bP,\cB)$ is {\em a model of iterable sweetness\/} if 
\begin{enumerate}
\item[(i)] $\cB$ is closed under finite intersections,
\item[(ii)] each $U\in \cB$ is directed and $p\leq q\in U\ \Rightarrow\
p\in U$,  
\item[(iii)] if $\langle p_n:n\leq \omega\rangle\subseteq U$ and the
sequence $\langle p_n:n<\omega\rangle$ converges to $p_\omega$ (in the
topology generated by $\cB$), then there is a condition $p\in U$ such that
$(\forall n\leq\omega)(p_n\leq p)$. 
\end{enumerate}
\end{definition}

\begin{proposition}
\label{swimitsw}
Assume that $(\bP,\cD,\bar{E})$ is a sweetness model on $\bP$, $\bar{E}=
\langle E_n:n<\omega\rangle$, $\cD=\bP$. Furthermore, suppose that any two
compatible elements of $\bP$ have a least upper bound, i.e., if $p_0,p_1\in
\bP$ are compatible, then there is $q\geq p_0,p_1$ such that
$(\forall r\in\bP)(r\geq p_0\ \&\ r\geq p_1\ \Rightarrow\ r\geq q)$. For
$\bar{p}=\langle p_\ell:\ell\leq k\rangle\subseteq\bP$ and $\bar{n}=\langle
n_\ell:\ell\leq k\rangle\subseteq\omega$ let  
\[U(\bar{p},\bar{n})\stackrel{\rm def}{=}\{q\in\bP: (\forall\ell\leq k)
(\exists q'\in [p_\ell]_{E_{n_\ell}})(q\leq q')\},\]
and let $\cB$ be the collection of all sets $U(\bar{p},\bar{n})$ (for
$\bar{p}\subseteq\bP$, $\bar{n}\subseteq\omega$, $\lh(\bar{p})=
\lh(\bar{n})$). Then $(\bP,\cB)$ is a model of iterable sweetness.
\end{proposition}

\begin{proof}
Plainly, $\cB$ is closed under finite intersections. Since the equivalence
classes of each $E_n$ are directed and any two compatible members of $\bP$
have a least upper bound, we may conclude that the elements of $\cB$ are
directed and downward closed.

Before we verify the demand \ref{itersweet}(iii), let us first note that if  
$p\in U(\bar{p},\bar{n})$, $\bar{p}=\langle p_\ell:\ell\leq k\rangle$,
$\bar{n}=\langle n_\ell:\ell\leq k\rangle$, then for some $N$ we have
$U(\langle p\rangle,\langle N\rangle)\subseteq U(\bar{p},\bar{n})$. [Why?
Use \ref{sweet}(v) to choose $N$ such that
\[(\forall\ell\leq k)(\forall q\in [p]_{E_N})(\exists q'\in
[p_\ell]_{E_{n_\ell}})(q\leq q').\]
Clearly this $N$ is as required.]

Now suppose that a sequence $\langle p_n:n<\omega\rangle$ converges to
$p_\omega$ in the topology generated by $\cB$, and $p_\omega,p_n\in
U(\bar{q},\bar{n})\in\cB$ for all $n<\omega$. Take $N$ such that $U(\langle
p_\omega\rangle,\langle N\rangle)\subseteq U(\bar{q},\bar{n})$. Choose an
increasing sequence $\langle m_i:i<\omega\rangle$ such that
\[(\forall i<\omega)(\forall n\geq m_i)(p_n\in U(\langle p_\omega\rangle,
\langle N+1+i\rangle)).\]
Next pick conditions $p^*_i\in [p_\omega]_{E_{N+1+i}}$ such that
\[(\forall i<\omega)(\forall n\in [m_i,m_{i+1}))(p_n\leq p^*_i)\]
(remember that each $[p_\omega]_{E_{N+1+i}}$ is directed). It follows from
\ref{sweet}(iv) that we may find a condition $q'\geq p_\omega$ such that
$q'\in [p_\omega]_{E_N}$ and $(\forall i<\omega)(p^*_i\leq q')$. Then 
\[(\forall n\geq m_0)(p_n\leq q')\quad\mbox{ and }\quad q'\in U(\bar{q},
\bar{n}).\]
Since $U(\bar{q},\bar{n})$ is directed and $q',p_0,\ldots, p_{m_0}\in
U(\bar{q},\bar{n})$, the conditions $q',p_0,\ldots, p_{m_0}$ have an upper
bound in $U(\bar{q},\bar{n})$ -- let $q$ be such an upper bound. 
\end{proof}

\begin{lemma}
\label{bastoplem}
Assume that $(\bP,\tau)$ is a model of topological sweetness.
\begin{enumerate}
\item If $p,q\in\bP$, $p\leq q$ and $q\in U\in\tau$, then there is an open
neighbourhood $V$ of $p$ such that
\[(\forall r\in V)(\exists r'\in U)(r\leq r').\]
\item If $m\in\omega$, $p\in U\in \tau$, then there is an open neighbourhood 
$V$ of $p$ such that any $p_0,\ldots,p_m\in V$ have a common upper bound in
$U$.
\end{enumerate}
\end{lemma}

\begin{proof}
Straightforward.
\end{proof}

\begin{theorem}
\label{getiter}
Suppose that $(\bP,\tau)$ is a model of topological sweetness and
$\dot{\cB},\dot{\bQ}$ are $\bP$--names such that
\[\forces_{\bP}\mbox{`` $(\dot{\bQ},\dot{\cB})$ is a model of iterable
sweetness ''}.\]
Then there is dense subset $\mbR$ of the iteration $\bP*\dot{\bQ}$ and a
topology $\tau^*$ on $\mbR$ such that $\bP\subseteq\mbR$ and $(\mbR,\tau^*)$
is a model of topological sweetness extending the model $(\bP,\tau)$.
\end{theorem}

\begin{proof}
Let $\dot{V}_n$ be $\bP$--names such that
\[\forces_{\bP}\mbox{`` }\dot{V}_0=\{\emptyset_{\dot{\bQ}}\}\ \mbox{ and }\
\{\dot{V}_n:n<\omega\}\mbox{ enumerates }\dot{\cB}\setminus\{\emptyset\}
\mbox{ ''}\]
(note that $\forces\mbox{`` }(\dot{\bQ},\{\dot{V}_n:n<\omega\})$ is a model
of iterable sweetness '' and also $\forces\mbox{`` }\emptyset_{\dbQ}\in 
\dot{V}_n$ for all $n$ '', remember \ref{itersweet}(ii)). Let $\cU$ be a
countable basis of the topology $\tau$ and let
\[\mbR\stackrel{\rm def}{=}\{(p,\dot{q})\in\bP*\dbQ: p\neq\emptyset_{\bP}
\mbox{ and } (\exists n<\omega)(p\forces\dot{q}\in\dot{V}_n)\}\cup
\{(\emptyset_\bP,\dot{\emptyset}_{\dbQ})\}.\]  
For $U\in\cU$, $\bar{U}=\langle U_0,\ldots, U_{m-1}\rangle\subseteq\cU$ and
$\bar{n}=\langle n_0,\ldots, n_{m-1}\rangle\subseteq\omega$ and $\bar{k}=
\langle k_0,\ldots, k_M\rangle\subseteq\omega$ ($m,M<\omega$) we put 
\[\begin{array}{ll}
U^*(U,\bar{U},\bar{n},\bar{k})\stackrel{\rm def}{=}\{(p,\dot{q})\in\mbR:&
p\in U\ \mbox{ and }\ (\forall\ell\leq M)(p\forces\dot{q}\in
\dot{V}_{k_\ell})\ \mbox{ and}\\
&(\forall\ell<m)(\exists p'\in U_\ell)(p\leq p'\ \&\ p'\forces\dot{q}\in
\dot{V}_{n_\ell})\}.
  \end{array}\]
Let $\cC$ be the collection of all sets of the form $U^*(U,\bar{U},\bar{n},
\bar{k})$ (for suitable parameters $U,\bar{U},\bar{n},\bar{k}$).

\begin{claim}
\label{cl14}
\begin{enumerate}
\item The family $\cC$ forms a countable basis of a topology on $\mbR$;
we will denote this topology by $\tau^*$. $\emptyset_{\mbR}$ is an isolated
point in $\tau^*$.   

\item If $(p,\dot{q}),(p',\dot{q}')\in\mbR$, $(p,\dot{q})\leq (p',\dot{q}')$
and $(p',\dot{q}')\in U^*\in\cC$, then there is $V^*\in\cC$ such that 
$(p,\dot{q})\in V^*$ and $(\forall r\in V^*)(\exists r'\in U^*)(r\leq r')$. 
\end{enumerate}
\end{claim}

\begin{proof}[Proof of the claim] 
1)\quad Should be clear.

\noindent 2)\quad Let $U^*=U^*(U,\langle U_0,\ldots,U_{m-1}\rangle,\langle
n_0,\ldots,n_{m-1}\rangle,\langle k_0,\ldots, k_M\rangle)$ and let
$p_\ell'\in U_\ell$ (for $\ell<m$) be such that $p'\leq p_\ell'$ and
$p_\ell'\forces\dot{q}'\in\dot{V}_{n_\ell}$. Choose $U_\ell'\in\cU$ such
that $p_\ell'\in U_\ell'\subseteq U_\ell$ and every two members of $U_\ell'$
have a common upper bound in $U_\ell$ (possible by \ref{bastoplem}(2)). Next
pick $U',U''\in\cU$ such that $p'\in U''\subseteq U'\subseteq U$ and 
\begin{itemize}
\item every member of $U'$ has an upper bound in each of the sets $U_0',
\ldots,U_{m-1}'$, 
\item every $M+1$ elements of $U''$ have a common upper bound in $U'$.
\end{itemize}
Let $U^+\in\cU$ be such that $p\in U^+$ and each element of $U^+$ has an
upper bound in $U''$, and let $k$ be such that
$p\forces\dot{q}\in\dot{V}_k$. Put 
\[V^*=U^*(U^+,\langle U_0',\ldots,U_{m-1}',U'',\ldots, U''\rangle,\langle
n_0,\ldots,n_{m-1},k_0,\ldots,k_M \rangle,\langle k\rangle)\in\cC.\]
First we show that $(p,\dot{q})\in V^*$. By our choices, $p\in U^+$ and
$p\forces\dot{q}\in\dot{V}_k$. For $\ell<m$, $p_\ell'\in U_\ell'$ is a
condition stronger than $p'\geq p$ and $p_\ell'\forces \dot{q}\leq\dot{q}'
\in\dot{V}_{n_\ell}$, so $p_\ell'\forces\dot{q}\in\dot{V}_{n_\ell}$
(remember \ref{itersweet}(ii)). Next, for $i\leq M$, we have $p'\in U''$ and
$p'\forces\dot{q}\leq\dot{q}'\in\dot{V}_{k_i}$, so $p'\forces\dot{q}\in
\dot{V}_{k_i}$.  

Now, suppose that $(p^*,\dot{q}^*)\in V^*$. Then $p^*\in U^+$ and we have
conditions $p^*_\ell\in U_\ell'$ (for $\ell<m$) and conditions $p^{**}_i\in
U''$ (for $i\leq M$) such that 
\begin{itemize}
\item $p^*_\ell\geq p^*$ and $p^{**}_i\geq p^*$, 
\item $p^*_\ell\forces \dot{q}^*\in \dot{V}_{n_\ell}$ and $p^{**}_i\forces
\dot{q}^*\in\dot{V}_{k_i}$. 
\end{itemize}
Pick a condition $p^+\in U'$ stronger than all $p^{**}_i$ (for $i\leq
M$). We claim that $(p^+,\dot{q}^*)\in U^*$. Clearly $p^+\in U$ and $p^+
\forces\dot{q}^*\in\dot{V}_{k_i}$ (for $i\leq M$). Fix $\ell<m$. By the
choice of $U'$, we find $p^+_\ell\in U_\ell'$ stronger than $p^+$. By the
choice of $U_\ell'$, we find a condition $p_\ell\in U_\ell$ stronger than
both $p^+_\ell$ and $p^*_\ell$. Then $p_\ell\forces\dot{q}^*\in
\dot{V}_{n_\ell}$ and we are done.
\end{proof}

\begin{claim}
\label{cl15}
Suppose that a sequence $\langle (p_k,\dot{q}_k):k<\omega\rangle\subseteq
\mbR$ is $\tau^*$--converging to $(p_\omega,\dot{q}_\omega)$, and 
$(p_\omega,\dot{q}_\omega)\in U^*(U,\langle U_0,\ldots,U_{m-1}\rangle,
\langle n_0,\ldots,n_{m-1}\rangle,\bar{k})\in\cC$. Then there are an
infinite set $X\subseteq\omega$ and conditions $p^*_\ell\geq p_\omega$ (for
$\ell<m$) such that for each $\ell<m$ and $k\in X$:
\begin{enumerate}
\item[(i)]  $p^*_\ell\in U_\ell$, $p^*_\ell\geq p_k$,
\item[(ii)] $p^*_\ell\forces\dot{q}_k\in\dot{V}_{n_\ell}$.
\end{enumerate}
\end{claim}

\begin{proof}[Proof of the claim] 
Pick $p^\ell_\omega\in U_\ell$, $p^\ell_\omega\geq p_\omega$ such that
$p^\ell_\omega\forces\dot{q}_\omega\in \dot{V}_{n_\ell}$ (for $\ell<m$). Fix
sequences $\langle W^\ell_j:j<\omega\rangle\subseteq \cU$ (for $\ell<m$)
such that
\begin{itemize}
\item $p^\ell_\omega\in W^\ell_{j+1}\subseteq W^\ell_j$,
\item $\{W^\ell_j:j<\omega\}$ forms a basis of neighbourhoods of
$p_\omega^\ell$. 
\end{itemize}
Clearly $(p_\omega,\dot{q}_\omega)\in U^*(U,\langle W^0_j,\ldots,W^{m-1}_j
\rangle, \langle n_0,\ldots,n_{m-1}\rangle,\bar{k})$ for every $j<\omega$,
so we may pick an increasing sequence $\langle k_j:j<\omega\rangle
\subseteq\omega$ such that
\[(\forall j<\omega)((p_{k_j},\dot{q}_{k_j})\in U^*(U,\langle W^0_j,
\ldots,W^{m-1}_j \rangle, \langle n_0,\ldots,n_{m-1}\rangle,\bar{k})).\]
Let $p^\ell_{k_j}\in W^\ell_j$ be such that $p_{k_j}\leq p^\ell_{k_j}$,
$p_{k_j}^\ell\forces\dot{q}_{k_j}\in\dot{V}_{n_\ell}$ (for $\ell<m$,
$j<\omega$). Each sequence $\langle p^\ell_{k_j}:j<\omega\rangle$
$\tau$--converges to $p^\ell_\omega$ so we may find an infinite set
$A\subseteq\omega$ and conditions $p^*_\ell\in U_\ell$ (for $\ell<m$) such
that 
\[p^*_\ell\geq p^\ell_\omega\geq p_\omega\quad\mbox{and}\quad (\forall j\in
A)(\forall \ell<m)(p^*_\ell\geq p^\ell_{k_j}).\]
Let $X=\{k_j:j\in A\}$. 
\end{proof}

\begin{claim}
\label{cl16}
Suppose $\langle (p_n,\dot{q}_n):n<\omega\rangle\subseteq\mbR$ is
$\tau^*$--converging to $(p_\omega,\dot{q}_\omega)$. Then there is
$X\in\iso$ such that:
\begin{enumerate}
\item[$(\boxtimes)$] {\em if\/} $(p_\omega,\dot{q}_\omega)\in U^*(U,\langle 
U_0, \ldots,U_{m-1}\rangle,\langle n_0,\ldots,n_{m-1}\rangle,\bar{k})\in
\cC$,\\ 
{\em then\/} for some $N\in\omega$ and $p^*_\ell\in U_\ell$ (for $\ell<m$)
we have:
\begin{enumerate}
\item[(i)]  $p^*_\ell\geq p_\omega$, $(\forall n\in X\setminus N)(p^*_\ell
\geq p_n)$, 
\item[(ii)] $(\forall n\in X\setminus N)(\forall \ell<m)(p^*_\ell\forces
\dot{q}_n\in \dot{V}_{n_\ell})$.
\end{enumerate}
\end{enumerate}
\end{claim}

\begin{proof}[Proof of the claim] 
Let $\langle U^*_i:i<\omega\rangle$ enumerate all sets $U^*\in\cC$ to which
$(p_\omega,\dot{q}_\omega)$ belongs. Apply \ref{cl15} to choose inductively
a decreasing sequence $\langle X_i:i<\omega\rangle\subseteq\iso$ such that
for each $i<\omega$:
\begin{quotation}
if $U^*_i=U^*(U,\langle U_0, \ldots,U_{m-1}\rangle,\langle n_0,\ldots,
n_{m-1}\rangle,\bar{k})$, then there are conditions $p^i_\ell\geq p_\omega$
(for $\ell<m$) satisfying 
\[(\forall\ell<m)(\forall n\in X_i)(p^i_\ell\in U_\ell\ \&\ p^i_\ell\geq
p_n\ \&\ p^i_\ell\forces\dot{q}_n\in \dot{V}_{n_\ell}).\]
\end{quotation}
Next pick an infinite set $X\subseteq\omega$ almost included in all $X_n$'s.
\end{proof}

\begin{claim}
\label{cl17}
Suppose that $\langle (p_n,\dot{q}_n):n<\omega\rangle\subseteq\mbR$
$\tau^*$--converges to $(p_\omega,\dot{q}_\omega)\in U^*\in\cC$. Then there
is a condition $(p^*,\dot{q}^*)\in U^*$ stronger than $(p_\omega,
\dot{q}_\omega)$ and such that
\[(\exists^\infty n\in\omega)((p_n,\dot{q}_n)\leq (p^*,\dot{q}^*)).\]
\end{claim}

\begin{proof}[Proof of the claim] 
Let $U^*=U^*(U,\langle U_0,\ldots\!,U_{m-1}\rangle,\langle n_0,\ldots\!,
n_{m-1}\rangle,\langle k_0,\ldots,k_M\rangle)$, and let $p^\ell_\omega\in
U_\ell$ (for $\ell<m$) be such that $p^\ell_\omega\geq p_\omega$ and
$p^\ell_\omega\forces\dot{q}_\omega\in\dot{V}_{n_\ell}$. Pick $U',U_0',
\ldots,U_{m-1}'\in\cU$ such that 
\begin{itemize}
\item $p_\omega\in U'\subseteq U$, $p^\ell_\omega\in U_\ell'\subseteq
U_\ell$,
\item any 3 elements of $U_\ell'$ have a common upper bound in $U_\ell$,
\item every element of $U'$ has an upper bound in each of $U_0',\ldots,
U_{m-1}'$. 
\end{itemize}
Apply \ref{cl16} to choose $X\in\iso$ such that the condition $(\boxtimes)$
of \ref{cl16} holds. Note that then
\[p_\omega\forces\mbox{`` the sequence }\langle \dot{q}_n:n\in X\rangle
\mbox{  is $\dot{\cB}$--convergent to $\dot{q}_\omega$ ''.}\]
[Why? If not, then we may pick a condition $r\geq p_\omega$ and an integer
$N$ such that 
\[ r\forces\mbox{`` }\dot{q}_\omega\in\dot{V}_N\ \&\ (\exists^\infty
n\in X)(\dot{q}_n\notin \dot{V}_N)\mbox{ ''.}\]
Pick $W\in\cU$ such that $r\in W$ and any 2 members of $W$ are compatible,
and apply $(\boxtimes)$ of \ref{cl16} to $U^*(U,\langle W\rangle,\langle
N\rangle,\langle k_0,\ldots,k_M\rangle)$. We get a condition $r^*\in W$ such
that  
\[(\forall^\infty n\in X)(r^*\geq p_n\ \&\ r^*\forces\dot{q}_n\in
\dot{V}_N).\] 
Since $r^*,r$ are compatible, we get a contradiction.]

Since $\langle p_n:n\in X\rangle$ $\tau$--converges to $p_\omega$, we may
find an infinite $X'\subseteq X$ and a condition $p^*\in U'$ such that 
\[p^*\geq p_\omega\ \&\ (\forall n\in X')(p^*\geq p_n).\]
Next use $(\boxtimes)$ of \ref{cl16} to pick conditions $p_\ell'\in U_\ell'$
and $N\in \omega$ such that
\[(\forall n\in X'\setminus N)\big(p_\ell'\geq p_n\ \mbox{ and }\
p_\ell'\forces \dot{q}_n\in\dot{V}_{n_\ell}\quad\mbox{and}\quad (\forall
i\leq M)(p_n\forces \dot{q}_n\in\dot{V}_{k_i})\big).\]
By the choice of $U',U_0',\ldots,U_{m-1}'$ we get conditions $p^*_\ell\in
U_\ell$ such that $p^*_\ell\geq p^*$, $p^*_\ell\geq p_\ell'$, $p^*_\ell\geq
p^\ell_\omega$ (for $\ell<m$). 

Now we are going to define a $\bP$--name $\dot{q}^*$ for a condition in
$\dot{\bQ}$. Let $\cA$ be a maximal antichain of $\bP$ such that for each
$\ell<m$ and $r\in\cA$:
\begin{itemize}
\item either $r\geq p^*_\ell$ or $r,p^*_\ell$ are incompatible, and
\item either $r\geq p^*$ or $r,p^*$ are incompatible.
\end{itemize}
Fix $r\in\cA$. If $r,p^*$ are incompatible, then let $\dot{q}_r$ be
$\dot{\emptyset}_{\dot{\bQ}}$. Assume $r\geq p^*$ and let 
\[I=\{n_\ell: \ell<m\ \&\ r\geq p^*_\ell\}\cup\{k_i:i\leq M\}\quad
(\neq\emptyset).\]
The condition $r$ forces that the sequence $\langle\dot{q}_n:n\in X'
\setminus N\rangle$ converges to $\dot{q}_\omega$, and $\dot{q}_\omega 
\in\bigcap\limits_{j\in I}\dot{V}_j$, and $\dot{q}_n\in\bigcap\limits_{j\in 
I}\dot{V}_j$ (for all $n\in X'\setminus N$). Applying \ref{itersweet}(i+iii)
we find a $\bP$--name $\dot{q}_r$ for an element of $\dot{\bQ}$ such that 
\[r\forces\mbox{`` }(\forall n\in X'\setminus N)(\dot{q}_n\leq\dot{q}_r\ \&\
\dot{q}_\omega\leq\dot{q}_r\ \&\ \dot{q}_r\in\bigcap_{j\in I}\dot{V}_j) 
\mbox{ ''.}\]
Now, let $\dot{q}^*$ be a $\bP$--name such that $r\forces\dot{q}^*=
\dot{q}_r$ (for $r\in\cA$). 

Look at the condition $(p^*,q^*)\in\mbR$. Clearly $p^*\forces\dot{q}^*\in
\dot{V}_{k_i}$ (for all $i\leq M$) and $p^*_\ell\geq p^*$, $p^*_\ell\forces
\dot{q}^*\in\dot{V}_{n_\ell}$ (for $\ell<m$), so $(p^*,\dot{q}^*)\in U^*$. 
Moreover, if $n\in (X'\setminus N)\cup\{\omega\}$, then $p^*\geq p_n$ and
$p^*\forces\dot{q}_n\leq \dot{q}$, so $(p^*,\dot{q}^*)\geq (p_n,\dot{q}_n)$
and we are done. 
\end{proof}

\begin{claim}
\label{clxyz}
If $(p,\dot{q})\in U^*\in \cC$, then there is $V^*\in\cC$ such that
$(p,\dot{q})\in V^*$ and any two conditions $(p_0,\dot{q}_0),(p_1,\dot{q}_1)
\in V^*$ have a common upper bound in $U^*$.
\end{claim}

\begin{proof}[Proof of the claim] 
Let $U^*=U^*(U,\langle U_0,\ldots,U_{m-1}\rangle, \langle n_0,\ldots,n_{m-1}
\rangle, \langle k_0,\ldots, k_M\rangle)$ and let $p_\ell\in U_\ell$ (for
$\ell<m$) be such that $p\leq p_\ell$ and $p_\ell\forces\dot{q}\in
\dot{V}_{n_\ell}$. Pick $U'_\ell\in\cU$ such that $p_\ell\in
U'_\ell\subseteq U_\ell$ and every three members of $U'_\ell$ have a common
upper bound in $U_\ell$. Also choose $U',U''\in\cU$ such that $p\in
U''\subseteq U'\subseteq U$, and each member of $U'$ has an upper bound in
every $U'_\ell$ (for $\ell<m$) and every two members of $U''$ have a common
upper bound in $U'$. Put
\[V^*=U^*(U'',\langle U'_0,\ldots,U'_{m-1}\rangle,\langle n_0,\ldots,n_{m-1}
\rangle, \langle k_0,\ldots,k_M\rangle)\in\cC.\]
Clearly $(p,\dot{q})\in V^*$. Suppose now that $(p_0,\dot{q}_0),(p_1,
\dot{q}_1)\in V^*$ and let $p^i_\ell\in U'_\ell$ be such that $p^i_\ell
\forces \dot{q}_i\in \dot{V}_{n_\ell}$ and $p_i\leq p^i_\ell$ (for $i=0,1$
and $\ell<m$). Also let $p^*\in U'$ be stronger than both $p_0$ and $p_1$,
and for each $\ell<m$ let $p^*_\ell\in U_\ell$ be stronger than both $p^*$
and $p^0_\ell$ and $p^1_\ell$. Now, like in the proof of \ref{cl17}, choose
a $\bP$--name $\dot{q}^*$ for a condition in $\dot{\bQ}$ such that 
\begin{quotation}
$p^*\forces$`` $\dot{q}^*\geq \dot{q}_0\ \&\ \dot{q}^*\geq \dot{q}_1$ '',
and\\
$p^*_\ell\forces$`` $\dot{q}^*\in \dot{V}_{n_\ell}$ '' for $\ell<m$, and\\
$p^*\forces$`` $\dot{q}^*\in \dot{V}_{k_\ell}$ '' for $\ell\leq M$.
\end{quotation}
(Remember that the sets $\bigcap\limits_{j\in I}\dot{V}_j$ are forced to be
directed by \ref{itersweet}(i+ii).) Then $(p^*,\dot{q}^*)\in U^*$ is a
condition stronger than both $(p_0,\dot{q}_0)$ and $(p_1,\dot{q}_1)$.  
\end{proof}

Now we may put together \ref{cl17}, \ref{cl14}(2) and \ref{clxyz} to
conclude that the topology $\tau^*$ satisfies the demand
\ref{topsweet}(iii), finishing the proof of the theorem.

Note that $(p,\dot{\emptyset}_{\dot{\bQ}})\in\mbR$ for each $p\in\bP$ and
the mapping $p\mapsto (p,\dot{\emptyset}_{\dot{\bQ}})$ is a homeomorphic
embedding of $(\bP,\tau)$ into $(\mbR,\tau^*)$, so we may think that $\tau$
is the restriction of $\tau^*$ to $\bP\subseteq\mbR$. Moreover, under this
interpretation, $\bP$ is an open subset of $\mbR$.
\end{proof}

Let us state (without a proof) a result that shows how \ref{getiter}
(together with \ref{sweetamal}) can be applied. 

\begin{context}
\label{contbis}
Let $\cK$ be a collection of tuples $\bar{a}=\langle\varphi_0,\varphi_1,
\varphi_2,h,W\rangle$, where $\varphi_i$ are $\Sigma^1_1$ formulas and $h,W$
are Borel functions such that for each real $r$:
\begin{itemize}
\item $\varphi_0(\cdot,r),\varphi_1(\cdot,r),\varphi_2(\cdot,r)$ define a
Souslin ccc forcing notion $\bP^{\bar{a},r}$ (and $\leq_{\bP^{\bar{a},r}}$,
$\incomp_{\bP^{\bar{a},r}}$),   
\item $h(r)$ is a function $\bH^{\bar{a},r}:\omega\longrightarrow\cHa$ (and 
we have $\cT^{\bar{a},r}_n=\prod\limits_{i<n}\bH^{\bar{a},r}(i)$ and
$\cT^{\bar{a},r}=\bigcup\limits_{n<\omega}\cT_n$),
\item $W(r)$ is a sequence $\dot{W}^{\bar{a},r}=\langle p_\eta:\eta\in
\cT^{\bar{a},r}\rangle\subseteq\bP^{\bar{a},r}$,
\item $\bP^{\bar{a},r},\bH^{\bar{a},r},\dot{W}^{\bar{a},r}$ are as in
\ref{cont}. 
\end{itemize}
$\cK_0,\cK_1$ are subsets of $\cK$ such that 
\begin{itemize}
\item for each $\bar{a}\in\cK_0$ and every real $r$, the forcing notion
$\bP^{\bar{a},r}$ is iterably sweet,
\item for each $\bar{a}\in\cK_1$ and a real $r$, there is an iterably sweet
Souslin ccc forcing notion $\bQ_{\bar{a},r}$ such that 
\[\begin{array}{ll}
\forces_{\bQ_{\bar{a},r}}&\mbox{`` the union of all Borel sets from the
ideal ${\mathcal I}_{\bP^{\bar{a},r},\dot{W}^{\bar{a},r}}$}\\
&\mbox{ \ coded in $\bV$ belongs to the ideal ${\mathcal I}_{\bP^{\bar{a},
r},\dot{W}^{\bar{a},r}}$ ''.}
  \end{array}\]
\end{itemize}
\end{context}

\begin{theorem}
Assume GCH. Let $\cK,\cK_0,\cK_1$ be as in \ref{contbis} and
$\kappa=\kappa^{<\kappa}$. Then there is a forcing notion $\bQ$ preserving
cofinalities and cardinalities and such that for every generic
$G\subseteq\bQ$ over $\bV$ we have (in $\bV[G]$): 
\begin{enumerate}
\item $\con=\kappa$,
\item for every real $r$ and $\bar{a}\in\cK_1$, every projective subset of
$\prod\limits_{n<\omega}\bH^{\bar{a},r}(n)$ has ${\mathcal I}_{\bP^{\bar{a},
r},\dot{W}^{\bar{a},r}}$--Baire property,
\item for each real $r$, a sequence $\langle r_\alpha:\alpha<\mu\rangle$ of
reals, $\mu<\kappa$, and $\bar{a}\in\cK_0$ there is an ${\mathcal
I}_{\bP^{\bar{a},r},\dot{W}^{\bar{a},r}}$--generic real $r^*$ over
$\bV[r,\langle r_\alpha:\alpha<\mu\rangle]$. 
\end{enumerate}
\end{theorem}

\begin{proof}
See \cite{Sh:F380}.
\end{proof}

\begin{proposition}
\label{TSforCREA}
\begin{enumerate}
\item In the cases discussed in \ref{cccthm}(1,2), the forcing notion
under considerations is topologically sweet, provided $K$ is countable. 
\item If $(K,\Sigma)$ is a really finitary (see \ref{finitary}) linked
tree--creating pair, and $f$ is a fast function, then the forcing notion
$\bQ^\tree_f(K,\Sigma)$ is topologically sweet. 
\item Let $\gp=(K,\Sigma,\cF,\cG)$ be a universality parameter for
$\bH$. Assume that 
\begin{enumerate}
\item[(a)] $\cF$ is either countable or ${<}\omega_1$---$\leq^*$--directed,
and 
\item[(b)] for each $\eta\in\bigcup\limits_{n<\omega}\prod\limits_{i<n}
\bH(i)$ there is $t^{\max}_\eta\in\TCR_\eta[\bH]\cap K$ such that 
\[(\forall t\in\TCR_\eta[\bH]\cap K)(t\in \Sigma(t^{\max}_\eta)).\] 
\end{enumerate}
Then the forcing notion $\bQp$ is sweet (and thus iterably sweet, provided
elements compatible in $\bQp$ have the least upper bound). 
\end{enumerate}
\end{proposition}

\begin{proof}
1)\quad Let $(K,\Sigma,\Sigma^\bot)$ be a $\otimes$--creating triple
for $\bH:\omega\longrightarrow\cHa$. Suppose that $(K,\Sigma,\Sigma^\bot)$
is linked, gluing and has the cutting property.

For ${\bf c}=(w,t_0,\ldots,t_n)\in\FC(K,\Sigma,\Sigma^\bot)$ and $N<\omega$
let 
\[\begin{array}{lr}
U({\bf c},N)=\{p\in\qinf:& w^p=w\quad \&\quad (\forall k\leq n)(t^p_k=t_k)
\quad \&\ \ \\ 
&(\forall k>n)(\nor[t^p_k]\geq N)\}.
  \end{array}\]
Let $\tau$ be the topology generated by the sets $U({\bf c},N)$ (for ${\bf
c}\in\FC(K,\Sigma,\Sigma^\bot)$ and $N<\omega$) and
$\{\emptyset_{\qinf}\}$. It is straightforward to check that $(\qinf,\tau)$
is a model of topological sweetness.
\medskip

\noindent Other instances of 1) and 2) can be handled similarly.
\medskip

\noindent 3)\quad We consider the case when $\cF$ is countable only (if
$\cF$ is ${<}\omega_1$---$\leq^*$--directed the proof is similar). We put
$\cD=\bQp$ and we define relations $E_n$ (for $n<\omega$) on $\cD$ as
follows: 

$p_0\;E_n\;p_1$\qquad if and only if

$N^{p_0}=N^{p_1}$, $f^{p_0}=f^{p_1}$ and 
\[(\forall\eta\in T^{p_0})(\lh(\eta)<N^{p_0}+n\quad\Rightarrow\quad\eta\in 
T^{p_1}\ \&\ t^{p_0}=t^{p_1}).\]
We claim that $(\bQp,\cD,\langle E_n:n<\omega\rangle)$ is a model of
sweetness. Plainly, each $E_n$ is an equivalence relation with countably
many equivalence classes, $E_{n+1}\subseteq E_n$. Similarly as in
\ref{cl18} one can show that the equivalence classes of each $E_n$ are
directed and that the clause \ref{sweet}(v) is satisfied. Let us show that
the demand \ref{sweet}(iv) holds.

So suppose $p_i\;E_i\;p_\omega$ (for $i\geq n$, $n<\omega$). Thus $f^{p_i}=
f^{p_\omega}=f$, $N^{p_i}=N^{p_\omega}=N$ (for $i\geq n$). We build
inductively a system $\langle t_\eta:\eta\in T\rangle\in\bQ^\tree_\emptyset
(K,\Sigma)$ as follows.

First we let $M_0=N+n$ and we declare 
\[\eta\in T^{p_\omega}\ \&\ \lh(\eta)<M_0\quad\Rightarrow\quad \eta\in T\
\&\ t_\eta=t^{p_\omega}_\eta.\]
Suppose that we have defined $\langle t_\eta:\eta\in T\ \&\ \lh(\eta)<
M_k\rangle\in\FC(K,\Sigma)$ already. Pick $M_k'>M_k$ such that for some
$n^k_\dn,n^k_\up$ we have $M_k<n^k_\dn\leq n^k_\up\leq F^\cG(n^k_\up)<M_k'$
and   
\[(\langle t^{p_\omega}_\eta:\eta\in T^{p_\omega}\ \&\ \lh(\eta)<M_k'
\rangle,n^k_\dn,n^k_\up, f\rest [n^k_\dn,n^k_\up])\in\cG.\] 
Next, choose $M_{k+1}$ and $n^{k,i}_\dn,n^{k,i}_\up$ (for $i\in [n,M'_k]$)
such that $M'_k<n^{k,n}_\dn$, $n^{k,i}_\dn\leq n^{k,i}_\up\leq
F^\cG(n^{k,i}_\up)<n^{k,i+1}_\dn$-2, $F^\cG(n^{k,M'_k}_\up)<M_{k+1}$ and  
\[(\langle t^{p_i}_\eta:\eta\in T^{p_i}\ \&\ \lh(\eta)<M_{k+1}\rangle,
n^{k,i}_\dn,n^{k,i}_\up, f\rest [n^{k,i}_\dn,n^{k,i}_\up])\in\cG.\]  
Let $\langle t_\eta^k:\eta\in \hat{S}_k\rangle\in\FC(K,\Sigma)$ be such that
$\mrot(S_k)=\langle\rangle$, $\lev(S_k)=M_{k+1}$, and $t^k_\eta=
t^{p_\omega}_\eta$ when $\lh(\eta)<M_k'$ and $t^k_\eta=t^{\max}_\eta$ when
$M_k'\leq\lh(\eta)<\lev(S_k)$. Apply repeatedly \ref{univpar}($\varepsilon$)
to get $\langle t_\eta:\eta\in T\ \&\ \lh(\eta)<M_{k+1}\rangle\in\FC(K,
\Sigma)$ such that 
\[(\langle t_\eta:\eta\in T\ \&\ \lh(\eta)<M_{k+1}\rangle,n^k_\dn,
n^{k,M'_k}_\up, f\rest [n^k_\dn,n^{k,M'_k}_\up])\in\cG,\]
and $\langle t_\eta:\eta\in T\ \&\ \lh(\eta)<M_{k+1}\rangle\leq\langle
t^k_\eta:\eta\in\hat{S}_k\rangle$, and 
\[\langle t_\eta:\eta\in T\ \&\ \lh(\eta)<M_{k+1}\rangle\leq \langle
t^{p_i}_\eta:\eta\in T^{p_i}\rangle\qquad\mbox{ for all }i\in [n,M_k'].\]
Note that then $\langle t_\eta:\eta\in T\ \&\ \lh(\eta)<M_{k+1}\rangle\leq
\langle t^{p_i}_\eta:\eta\in T^{p_i}\rangle$ for all $i\in [n,\omega]$.

After the construction is carried out one easily checks that $q=(N,\langle
t_\eta:\eta\in T\rangle,f)\in\bQp$ is a condition stronger than all $p_i$'s
(for $i\in [n,\omega]$) and $q\;E_n\;p_\omega$. 
\end{proof}

\begin{definition}
\label{plusF}
Suppose that a function $h:\omega\times\omega\longrightarrow\omega$ is
regressive, $\cF\subseteq(\omega\setminus 2)^{\textstyle\omega}$, and
$(K,\Sigma)$ is a local creating pair for $\bH$. We define $\bQ^h_\cF(K,
\Sigma)$ as the suborder of $\bQ^*_\cF(K,\Sigma)$ consisting of conditions
$p\in\bQ^*_\cF(K,\Sigma)$ such that 
\[(\exists f\in\cF)(\forall k<\omega)(\forall^\infty n<\omega)(h^{(k)}(
m^{t^p_n},\nor[t^p_n])\geq f(m^{t^p_n})),\]
where $h^{(k+1)}(i,j)=h^{(k)}(i,h(i,j))$, $h^{(0)}(i,j)=j$. 
\end{definition}

\begin{remark}
Note that the norm condition introduced in \ref{plusF} is in many cases
nothing new. If for each $f\in\cF$ there is $f^+\in\cF$ such that 
\[(\forall k<\omega)(\forall^\infty n<\omega)(h^{(k)}(n,f(n))\geq f^+(n)),\]
then clearly $\bQ^*_\cF(K,\Sigma)=\bQ^h_\cF(K,\Sigma)$.
\end{remark}

\begin{proposition}
\label{Ftopsweet}
Let $h:\omega\times\omega\longrightarrow\omega$ be regressive and 
$\cF\subseteq (\omega\setminus 2)^{\textstyle\omega}$ be a countable
family. Assume that $(K,\Sigma)$ is a local $h$--linked creating pair for
$\bH$ such that $K$ is countable and 
\begin{enumerate}
\item[$(*)$] if $s,t\in K$ and $s\in \Sigma(t)$, then $\nor[s]\leq
\nor[t]$. 
\end{enumerate}
Then the forcing notion $\bQ^h_\cF(K,\Sigma)$ is topologically sweet.
\end{proposition}

\begin{proof}
For a finite candidate ${\bf c}=(w,t_0,\ldots,t_m)\in\FC(K,\Sigma)$ and
sequences $\bar{f}=\langle f_\ell:\ell\leq\ell^*\subseteq\cF$ and $\bar{k}
=\langle k_\ell:\ell\leq\ell^*\rangle\subseteq\omega$ we let
\[\begin{array}{r}
U({\bf c},\bar{f},\bar{k})=\{p\in\bQ^h_\cF(K,\Sigma): w^p=w\quad \&\quad
(\forall n\leq m)(t^p_n=t_n) \quad \&\quad \\ 
(\forall\ell\leq\ell^*)(\forall n>m)(h^{(k_\ell)}(m^{t^p_n}_\dn,\nor[
t^p_n])\geq f_\ell(m^{t^p_n}))\quad \&\quad \\ 
(\forall\ell\leq\ell^*)(\forall k<\omega)(\forall^\infty n<\omega)
(h^{(k)}(m^{t^p_n}_\dn,\nor[t^p_n])\geq f_\ell(m^{t^p_n}))\ \}.
  \end{array}\]
Plainly the sets $U({\bf c},\bar{f},\bar{k})$ (for suitable ${\bf
c},\bar{f},\bar{k}$) and $\{\emptyset_{\bQ^h_\cF(K,\Sigma)}\}$ constitute a
basis of a topology $\tau$ on $\bQ^h_\cF(K,\Sigma)$. it is not difficult to
check that $(\bQ^h_\cF(K,\Sigma),\tau)$ is a model of topological sweetness.
\end{proof}

\subsection{The sour part of the spectrum}
The main point of the sweet properties is that amalgamations of sweet
forcing notions are ccc, see \ref{sweetamal}. A kind of opposite behaviour is
when amalgamating results in a forcing notion that collapses $\omega_1$. This
effect will be called sourness, and we have a number of variants of it (see
\ref{sour} below). 

In this part, whenever we use $\bQ$--names for elements of a space $\cX=
\prod\limits_{i<\omega}\bH(i)$, we assume that they are in {\em a standard
form}. Thus, a $\bQ$--name $\dot{\tau}$ for a real in $\cX$ is a system 
\[\langle q^n_{\eta,k}:n<\omega\ \&\ \eta\in\prod_{i<n}\bH(i)\ \&\ k<N_\eta
\rangle,\]
where $q^n_{\eta,k}\in\bQ$, $N_\eta\leq\omega$ and for each $n<\omega$
\[\langle q^n_{\eta,k}:\eta\in\prod_{i<n}\bH(i)\ \&\ k<N_\eta\rangle\]
is a maximal antichain of $\bQ$, and 
\[\nu\vartriangleleft\eta\quad \Rightarrow\quad (\forall k<N_\eta)(\exists
m<N_\nu)(q^{\lh(\nu)}_{\nu,m}\leq q^{\lh(\eta)}_{\eta,k}).\]
(The intension is that $q^n_{\eta,k}\forces\eta\vartriangleleft\dot{\tau}$.)
If ($\bQ$ is ccc and) $\dot{\tau}_0$ is a $\bQ$--name for an element of
$\cX$ then there is a standard name $\dot{\tau}$ such that $\forces
\dot{\tau}_0=\dot{\tau}$. The point of using standard names is that their
specific forms allows us to bound the complexity of some formulas; see,
e.g., \cite[3.6.12]{BaJu95}.
 
\begin{context}
\label{soucon}
Let $\bP,\dot{W},\cX$ be as in \ref{cont}, and let $\bQ_0,\bQ_1$ be Souslin
ccc forcing notions. Suppose that, for $\ell<2$, $\dot{\tau}_\ell$ is a
standard $\bQ_\ell$--name for an element of $\cX$. Furthermore, suppose that
there are isomorphisms $f_\ell:{\rm RO}(\bP)\stackrel{\rm onto}{
\longrightarrow}({\rm RO}(\bQ_\ell))_{\dot{\tau}_\ell}$ mapping $\dot{W}$
onto $\dot{\tau}_\ell$ (i.e., $f_\ell(\lbv\eta\vartriangleleft\dot{W}
\rbv_{{\rm RO}(\bP)})=\lbv\eta\vartriangleleft\dot{\tau}_\ell\rbv_{{\rm
RO}(\bQ_\ell)}$).\\ 
Let $r$ be a real encoding all parameters required for the definitions of
the above objects (including partial orders and the respective
incompatibility relations). 
\end{context}

\begin{definition}
\label{sour}
Let $\bP,\dot{W},\cX,\bQ_i,\dot{\tau}_\ell,f_\ell,r$ be as in \ref{soucon}. 
\begin{enumerate}
\item[(a)] The amalgamation $\bQ_0\times_{f_0,f_1}\bQ_1$ will be also
denoted $\bQ_0\times_{\dot{\tau}_0=\dot{\tau}_1}\bQ_1$. 
\end{enumerate}
We say that 
\begin{enumerate}
\item[(b)] {\em $(\bQ_0,\dot{\tau}_0)$ is weakly sour to $(\bQ_1,
\dot{\tau}_1)$} if the amalgamation $\bQ_0\times_{\dot{\tau}_0=\dot{\tau}_1}
\bQ_1$ fails the ccc;
\item[(c)] {\em $(\bQ_0,\dot{\tau}_0)$ is sour to $(\bQ_1,\dot{\tau}_1)$}
whenever the following condition holds:
\begin{enumerate}
\item[$(\boxplus)$] if $\bV\subseteq\bV'$ are universes of {\rm ZFC},
$r\in\bV$, $G_0,G_1\in\bV'$, $G_\ell\subseteq\bQ_\ell^\bV$ is generic
over $\bV$ and $\dot{\tau}^{G_0}_0=\dot{\tau}^{G_1}_1$, \quad then
$\bV'\models$`` $\omega_1^\bV$ is countable ''; 
\end{enumerate}
\item[(d)] {\em $(\bQ_0,\dot{\tau}_0)$ is explicitly sour to
$(\bQ_1,\dot{\tau}_1)$} if there are sequences $\langle E_m:m<\omega\rangle$
and $\langle q^\ell_{\alpha,n}:\alpha<\omega_1\ \&\ n<\omega\rangle$ (for
$\ell<2$) such that 
\begin{enumerate}
\item[(i)]    each $E_m$ is an equivalence relation on $\omega_1$ with at
most countably many equivalence classes, 
\item[(ii)]   $\{q^\ell_{\alpha,n}:n<\omega\}\subseteq\bQ_\ell$ is predense
in $\bQ_\ell$ (for each $\alpha<\omega_1$, $\ell<2$), 
\item[(iii)]  if $\alpha<\beta<\omega_1$, $m<\omega$, $\alpha\; E_m\;\beta$
and $n_0,n_1<m$ and both $(q^0_{\alpha,n_0},q^1_{\alpha,n_1})$ and
$(q^0_{\beta,n_0},q^1_{\beta,n_1})$ are in the amalgamation $\bQ_0 
\times_{\dot{\tau}_0=\dot{\tau}_1} \bQ_1$, \quad then the conditions
$(q^0_{\alpha,n_0},q^1_{\alpha,n_1}),(q^0_{\beta,n_0},q^1_{\beta,n_1})$ are 
incompatible (in $\bQ_0\times_{\dot{\tau}_0=\dot{\tau}_1}\bQ_1$);
\end{enumerate}
\item[(e)] {\em $(\bQ_0,\dot{\tau}_0)$ is very explicitly sour to
$(\bQ_1,\dot{\tau}_1)$} if there are Borel functions $g_\ell:\cX\times\omega
\longrightarrow\bQ_\ell$ such that 
\begin{enumerate}
\item[(i)] $\{g_\ell(x,n):n<\omega\}$ is predense in $\bQ_\ell$ (for each
$x\in\cX$, $\ell<2$),
\item[(ii)] if $x_0,x_1\in\cX$ are distinct, $x_0\rest m= x_1\rest m$, $k_0,
k_1<m$, then there are $n<\omega$ and disjoint sets
$A_0,A_1\subseteq\prod\limits_{i<n}\bH(i)$ such that for $\ell<2$:
\[(\forall q\in\bQ_\ell)([q\geq g_\ell(x_0,k_\ell)\ \&\ q\geq
g_\ell(x_1,k_\ell)]\quad\Rightarrow\quad q\forces\dot{\tau}_\ell\rest n\in
A_\ell).\]
\end{enumerate}
\end{enumerate}
We say that {\em $\bQ_0$ is sour to $\bQ_1$ over $(\bP,\dot{W})$} if there
are names $\dot{\tau}_0,\dot{\tau}_1$ (as above) such that
$(\bQ_0,\dot{\tau}_0)$ is sour to $(\bQ_1,\dot{\tau}_1)$ (and similarly for
the variants).\\ 
A forcing notion $\bQ$ is {\em sour over $(\bP,\dot{W})$} if it is sour to
itself over $(\bP,\dot{W})$ (and similarly for the other notions).\\
We may skip $\dot{W}$ and say ``over $\bP$'' if it is clear what $\dot{W}$
we consider.
\end{definition}

\begin{remark}
Sourness (\ref{sour}(c)) is a strong way to say that the amalgamation $\bQ_0
\times_{\dot{\tau}_0=\dot{\tau}_1}\bQ_1$ collapses $\omega_1$. Explicit
sourness (\ref{sour}(d)) guarantees that we have a nice witness for the
collapse, see \ref{souimp} below. What is the point of ``very explicitly
sour''? On one hand that condition implies that the amalgamation collapses
the continuum, and on the other hand the name for the collapsing function is
encoded in a nice way by a real. Note that the properties of the functions
$g_0,g_1$ required in \ref{sour}(e)(i,ii) are $\Pi^1_1$ (remember that
$\dot{\tau}_\ell$ are standard names), and thus we have suitable
absoluteness.  
\end{remark}

\begin{proposition}
\label{souimp}
Let $(\bQ_0,\dot{\tau}_0),(\bQ_1,\dot{\tau}_1),(\bP,\dot{W})$ and $r$ be as
in \ref{soucon}.
\begin{enumerate}
\item If $(\bQ_0,\dot{\tau}_0)$ is explicitly sour to $(\bQ_1,\dot{\tau}_1)$
in every universe $\bV$ of {\rm ZFC} containing $r$, then
$(\bQ_0,\dot{\tau}_0)$ is sour to $(\bQ_1,\dot{\tau}_1)$.
\item If $(\bQ_0,\dot{\tau}_0)$ is very explicitly sour to
$(\bQ_1,\dot{\tau}_1)$ with functions $g_0,g_1$ witnessing this, then
$(\bQ_0,\dot{\tau}_0)$ is explicitly sour to $(\bQ_1,\dot{\tau}_1)$ in every
universe of {\rm ZFC} containing $r$ and (the Borel codes for) $g_0,g_1$.
\end{enumerate}
\end{proposition}

\begin{proof}
1)\quad Suppose that $\bV\subseteq\bV'$ are universes of {\rm ZFC}, $r\in\bV$,
$G_0,G_1\in\bV'$, $G_\ell\subseteq\bQ_\ell^\bV$ is generic over $\bV$ and
$\dot{\tau}_0^{G_0}=\dot{\tau}_1^{G_1}$. 

Assume $\bV'\models\omega_1^\bV=\omega_1$.

\begin{claim}
\label{cl19}
\begin{enumerate}
\item If $q_\ell\in\bQ_\ell$ and $(q_0,q_1)\notin\bQ_0\times_{\dot{\tau}_0
=\dot{\tau}_1}\bQ_1$, then for some disjoint Borel sets $B_0,B_1\subseteq\cX$ 
we have $q_0\forces\dot{\tau}_0\in B_0$ and $q_1\forces\dot{\tau}_1\in
B_1$. Consequently, if $q_0\in G_0$ and $q_1\in G_1$, then $(q_0,q_1)\in 
\bQ_0\times_{\dot{\tau}_0=\dot{\tau}_1}\bQ_1$.
\item If $(q_0,q_1),(q_0',q_1')\in \bQ_0\times_{\dot{\tau}_0=\dot{\tau}_1}
\bQ_1$ are incompatible in $\bQ_0\times_{\dot{\tau}_0=\dot{\tau}_1}\bQ_1$,
then for some disjoint Borel sets $B_0,B_1\subseteq\cX$ we have (for
$\ell=0,1$): 
\[(\forall p\in\bQ_\ell)(p\geq q_\ell\ \&\ p\geq q_\ell'\quad\Rightarrow
\quad p\forces\dot{\tau}_\ell\in B_\ell).\] 
\end{enumerate}
\end{claim}

\begin{proof}[Proof of the claim] 
Straightforward if you remember \ref{quocon}.
\end{proof}

Now, let $\langle E_m:m<\omega\rangle,\langle q^\ell_{\alpha,k}:\alpha<
\omega_1^\bV, k<\omega\rangle\in\bV$ witness that $(\bQ_0,\dot{\tau}_0)$ is
explicitly sour to $(\bQ_1,\dot{\tau}_1)$ (in $\bV$). By \ref{sour}(d)(ii)
we know that, in $\bV'$,
\[(\forall\alpha<\omega_1)(\exists k^0_\alpha,k^1_\alpha<\omega)(q^0_{
\alpha,k^0_\alpha}\in G_0\ \&\ q^1_{\alpha,k^1_\alpha}\in G_1).\]
For some $k^0,k^1<\omega$ the set
\[Y=\{\alpha<\omega_1:k^0_\alpha=k^0\ \&\ k^1_\alpha=k^1\}\]
is uncountable. Let $m=k^0+k^1+1$. It follows from \ref{sour}(d)(i) that
there are distinct $\alpha,\beta\in Y$ such that $\alpha\; E_m\;\beta$. By
\ref{cl19}(1) we know $(q^0_{\alpha,k^0},q^1_{\alpha,k^1}),(q^0_{\beta,k^0},
q^1_{\beta,k^1})\in \bQ_0\times_{\dot{\tau}_0=\dot{\tau}_1}\bQ_1$, so by
\ref{sour}(d)(iii) these two conditions are incompatible in $\bQ_0
\times_{\dot{\tau}_0=\dot{\tau}_1}\bQ_1$. But then, using \ref{cl19}(2), we
find disjoint Borel sets $B_0,B_1\subseteq\cX$ such that $\dot{\tau}_0^{G_0}
\in B_0$ and $\dot{\tau}_1^{G_1}\in B_1$, a contradiction to
$\dot{\tau}_0^{G_0}=\dot{\tau}_1^{G_1}$.  
\medskip 

\noindent 2)\quad Let $\bV$ contain $r$ and $g_0,g_1$. Working in $\bV$,
pick a sequence $\langle x_\alpha:\alpha<\omega_1\rangle$ of pairwise
distinct members of $\cX$. Put:
\begin{itemize}
\item $q^0_{\alpha,n}=g_0(x_\alpha,n)$, $q^1_{\alpha,n}=g_1(x_\alpha,n)$,
\item $\alpha\; E_m\;\beta$ if and only if $x_\alpha\rest m= x_\beta\rest
m$. 
\end{itemize}
Now check.
\end{proof}

\begin{proposition}
\label{noBaire}
Let $r,\bP,\dot{W},\bQ_\ell,\dot{\tau}_\ell$ (for $\ell<2$) be as in
\ref{soucon}. Moreover, let $\bQ_\ell,\dot{W}_\ell,\cX_\ell,\bH_\ell$ (for
$\ell<2$) be as in \ref{cont} (with the real $r$ encoding needed parameters).
Assume that 
\begin{enumerate}
\item[(a)] $\omega_1$ is not an inaccessible cardinal in ${\bf L}$,
\item[(b)] $(\bQ_0,\dot{\tau}_0)$ is sour to $(\bQ_1,\dot{\tau}_1)$,
\item[(c)] for every real $s$ and a Borel set $B\subseteq\cX_\ell$ coded in
${\bf L}[r,s]$, if $B\notin {\mathcal I}_{\bQ_\ell,\dot{W}_\ell}$ then there
is an ${\mathcal I}_{\bQ_\ell,\dot{W}_\ell}$--generic real over ${\bf
L}[r,s]$ belonging to the set $B$.
\end{enumerate}
Then there is a $\Sigma^1_3$ set which does not have $\Ipw$--Baire property.
\end{proposition}

\begin{proof}
Since $\bQ_\ell,\dot{W}_\ell$ are as in \ref{cont}, we have a nice
description of $\bQ_\ell$--names for reals, see \ref{namrem}. So we have
Borel functions $h_\ell:\cX_\ell\longrightarrow\cX$ such that
$\forces_{\bQ_\ell}\dot{\tau}_\ell=h_\ell(\dot{W}_\ell)$. It follows from
the assumed properties of $\dot{\tau}_\ell$ that
\begin{enumerate}
\item[$(*)_1$] for every Borel set $A\subseteq\cX$, if $A\notin \Ipw$ then
$h_\ell^{-1}[A]\notin{\mathcal I}_{\bQ_\ell,\dot{W}_\ell}$. 
\end{enumerate}
(Note that $h_\ell$ is coded by the real $r$ as well.) Since $\omega_1$ is
not inaccessible in ${\bf L}$, for some real $a$ we have $\omega_1^{{\bf
L}[a]}=\omega_1$. Let
\[X_\ell=\{x\!\in\!\cX:\mbox{for some ${\mathcal I}_{\bQ_\ell,
\dot{W}_\ell}$--generic real $y\in\cX_\ell$ over ${\bf L}[a,r]$ we have
$h_\ell(y)=x$}\}.\]
(Note that if $y$ is ${\mathcal I}_{\bQ_\ell,\dot{W}_\ell}$--generic over
${\bf L}[a,r]$, then it determines a $\bQ_\ell^{{\bf L}[a,r]}$--generic
filter $G\subseteq\bQ_\ell^{{\bf L}[a,r]}$ over ${\bf L}[a,r]$ and
$h_\ell(y)=\dot{\tau}_\ell^G$; remember \ref{element} and the choice of
$h_\ell$.) 
\begin{enumerate}
\item[$(*_2)$] $X_\ell$ is a $\Sigma^1_3$ subset of $\cX_\ell$.
\end{enumerate}
[Why? Note that, by \ref{basideal}(4), the formula ``$c$ is a Borel code for 
a set $\sharp c\subseteq\cX_\ell$ and $\sharp c\notin{\mathcal I}_{\bQ_\ell,
\dot{W}_\ell}$'' is (equivalent in {\rm ZFC} to) a $\Pi^1_2$
formula. Hence easily the formula ``$y$ is ${\mathcal
I}_{\bQ_\ell,\dot{W}_\ell}$--generic over ${\bf L}[a,r]$'' is $\Pi^1_2$.]
\begin{enumerate}
\item[$(*_3)$] For every Borel set $B\subseteq\cX$, if $B\notin\Ipw$ then
$X_\ell\cap B\neq\emptyset$.
\end{enumerate}
[Immediate by $(*_1)$ and the assumption (c).]
\begin{enumerate}
\item[$(*_4)$] $X_0\cap X_1=\emptyset$.
\end{enumerate}
[It follows from the sourness and the assumption that $\omega_1^{{\bf L}[a]}
=\omega_1$.]

Finally note that $(*_3)+(*_4)$ implies that both sets $X_0,X_1$ do not have
the $\Ipw$--Baire property.   
\end{proof}

\begin{corollary}
\label{CornoBaire}
Suppose that $r,\bP,\dot{W},\bQ_\ell,\dot{\tau}_\ell,\dot{W}_\ell$ are as in 
\ref{noBaire}, and clauses (a)--(c) there hold, and $\bQ_0=\bQ_1$,
$\dot{W}_0=\dot{W}_1$. Assume additionally that, for some $k<\omega$, there
are $\bQ_0^{(k)}$--names $\dot{\rho}_\ell^0,\dot{\rho}_\ell^1$ (for $\ell<2$)
and a $\bQ_0$--name $\dot{\tau}$ such that 
\begin{enumerate}
\item[(d)] $\forces_{\bQ_0}$`` $\dot{\tau}\in\cX$ is $\Ipw$--generic over
$\bV$ '',
\item[(e)] $\forces_{\bQ^{(k)}_0}$`` $\dot{\rho}_\ell^0,\dot{\rho}_\ell^1\in 
\cX_0$ are ${\mathcal I}_{\bQ_0,\dot{W}_0}$--generic over $\bV$ '', and for
$\ell<2$ and every $\cI_{\bQ_0,\dot{W}_0}$--positive Borel set $B\subseteq
\cX_0$ there is $q\in \bQ^{(k)}_0$ such that $q\forces$``
$\dot{\rho}^0_\ell\in B$ '', 
\item[(f)] $\forces_{\bQ^{(k)}_0}$`` $\dot{\tau}[\dot{\rho}_\ell^0]=
\dot{\tau}_\ell[\dot{\rho}_\ell^1]$ ''.
\end{enumerate}
(Above, $\bQ^{(k)}_0$ stands for the iteration of length $k$ of the forcing
notion $\bQ_0$.)\\ 
Then there is a projective subset of $\cX_0$ that does not have the
${\mathcal I}_{\bQ_0,\dot{W}_0}$--Baire property.
\end{corollary}

Let us turn to getting sourness for some of the forcing notions discussed in
this paper. Of course, because of \ref{sweetamal} the forcing notions
covered by \ref{TSforCREA} are not sour (they are in the sweet kingdom,
after all). However, there are sour examples around. Let us introduce them
starting with exotic norm conditions which were chosen specially with the
sourness in mind.

\begin{proposition}
\label{getsour1}
Let $(K,\Sigma)$ be a local forgetful and complete (see \ref{complete})
creating pair for $\bH$.
\begin{enumerate}
\item Assume $(K,\Sigma)$ is linked and
\begin{enumerate}
\item[$(\alpha)$] $(\forall n<\omega)(|\bH(n)|>2^n)$,
\item[$(\beta)$]  $(\bar{K},\bar{g})$ is a 1--norming system for $\bH$ (see
\ref{normK}), $\bar{K}=\langle K_\ell:\ell<\omega\rangle$, $\bar{g}=\langle
g_\rho: \rho\in\fs\rangle$,
\item[$(\gamma)$] if $A\subseteq\bH(n)$, $a\in A$ and $\nor[t^n_A]>1$ then
$\nor[t^n_{A\setminus\{a\}}]\geq \nor[t^n_A]-1$,
\item[$(\delta)$] letting $A_n=\bH(n)\setminus\{g_\rho(n):\rho\in
2^{\textstyle\ell}\}$ for $n\in K_\ell$, and $A_n=\bH(n)$ for $n\notin
\bigcup\limits_{\ell<\omega} K_\ell$, we have $\lim\limits_{n\to\infty}
\nor[t^n_{A_n}]=\infty$.
\end{enumerate}
(Above, for $n\in\omega$ and $A\subseteq\bH(n)$, $t^n_A$ is the unique
creature $t\in K$ with $m^t_\dn=n$ and $\pos(t)=A$; see \ref{complete}.)\\
Then the forcing notion $\bQ^{\bar{K},\bar{g}}_\infty(K,\Sigma)$ (see
\ref{normK}) is very explicitly sour over Cohen.
\item If $(K,\Sigma), (\bar{K},\bar{g})$ and $\bH$ are as in (1), and
$f:\omega\times\omega\longrightarrow\omega$ is fast and 
\[(\forall n<\omega)(\forall i\geq n)(f(n,i)\leq\nor[t^i_{\bH(i)}]),\]
then $\bQ^{\bar{K},\bar{g}}_f(K,\Sigma)$ is very explicitly sour over Cohen.
\item Assume that $h:\omega\times\omega\longrightarrow\omega$ is regressive,
$\cF\subseteq(\omega\setminus 2)^{\textstyle\omega}$ is $h$--closed,
countable and  $\geq^*$--directed. Suppose that $(K,\Sigma)$ is $h$--linked
and clauses $(\alpha)$---$(\gamma)$ hold true, and  
\begin{enumerate}
\item[$(\delta^+)$] for $A_n$ as in $(\delta)$, for each $f\in\cF$ we have
$(\forall^\infty n\in\omega)(f(n)\leq\nor[t^n_{A_n}])$.
\end{enumerate}
Then the forcing notion $\bQ^{\bar{K},\bar{g}}_\cF(K,\Sigma)$ is very
explicitly sour over Cohen.
\item Assume that $(K,\Sigma)$ is linked, satisfies the demand 1($\gamma$),
and $\lim\limits_{n\to\infty}\nor[t^n_{\bH(n)}]=\infty$. Let 
$\bar{U}=\langle U_{\rho,k}: \rho\in\fs\ \&\ k<\omega\rangle$ be a
2--norming system (see \ref{normU}). Then the forcing notion
$\bQ^{\bar{U}}_\infty(K,\Sigma)$ (see \ref{normU}) is very explicitly sour
over Cohen. 
\item Similarly for forcing notions $\bQ^{\bar{U}}_f(K,\Sigma)$,
$\bQ^{\bar{U}}_\cF(K,\Sigma)$ (under assumptions parallel to that in 2,3
above, with clause ($\delta^+$) replaced by just $f(n)\leq\nor[
t^n_{\bH(n)}]$ for $n<\omega$, $f\in\cF$). 
\end{enumerate}
\end{proposition}

\begin{proof}
1)\quad Let us think about the Cohen forcing $\bC$ as the set of partial
functions from $\fs$ to 2 with the relation of extension; $\bH_\bC(\rho)=2$
for $\rho\in\fs$ (so we interpret $\omega$ as $\fs$), and $\dot{W}_\bC$ is
the natural name for the $\bC$--generic real in $2^{2^{{<}\omega}}$. 

Let $\dot{W}$ be the name for $\bQ^{\bar{K},\bar{g}}_\infty(K,
\Sigma)$--generic real, i.e., $\forces\dot{W}=\bigcup\{w^q:q\in\Gamma_{
\bQ^{\bar{K},\bar{g}}_\infty(K,\Sigma)}\}$. Let $\dot{\tau}_0,\dot{\tau}_1$
be standard $\bQ^{\bar{K},\bar{g}}_\infty(K,\Sigma)$--names for functions
from $\fs$ to 2 such that 
\begin{itemize}
\item $\dot{\tau}_0(\rho)=1$ if and only if $(\exists m\in K_{\lh(\rho)+1})
(\dot{W}(m)\in\{g_{\rho\conc\langle 0\rangle}(m),g_{\rho\conc\langle 1
\rangle}(m)\})$, 
\item $\dot{\tau}_1(\rho)=1-\dot{\tau}_0(\rho)$.
\end{itemize}

\begin{claim}
\label{cl20}
\[\forces_{\bQ^{\bar{K},\bar{g}}_\infty(K,\Sigma)}\mbox{`` }\dot{\tau}_0,
\dot{\tau}_1\mbox{ are $\bC$--generic over $\bV$ '';}\]
moreover $\lbv\dot{\tau}_\ell\in B\rbv_{\bQ^{\bar{K},\bar{g}}_\infty(K,
\Sigma)}\neq {\bf 0}$ for any non-meager Borel set $B\subseteq 2^{2^{{<}
\omega}}$. Hence $(\bC,\dot{W}_\bC), (\bQ^{\bar{K},\bar{g}}_\infty(K,
\Sigma),\dot{\tau}_0),(\bQ^{\bar{K},\bar{g}}_\infty(K,\Sigma),\dot{\tau}_1)$
are as in \ref{soucon}.
\end{claim}

\begin{proof}[Proof of the claim] 
Let $q\in \bQ^{\bar{K},\bar{g}}_\infty(K,\Sigma)$ and let $\rho_0,\ldots,
\rho_k\in\can$ be such that
\[\rho\in 2^{\textstyle\ell}\ \&\ m^{t^q_i}_\dn\in K_\ell\ \&\ g_\rho(
m^{t^q_i}_\dn)\notin\pos(t^q_i)\quad\Rightarrow\quad \rho\vartriangleleft
\rho_0\ \vee\ldots\vee\ \rho\vartriangleleft\rho_k.\]
Pick $N\geq\lh(w^q)$ such that $\rho_0\rest N,\ldots,\rho_k\rest N$ are
pairwise distinct and
\begin{enumerate}
\item[(a)] $m^{t^q_i}_\dn\geq N\quad\Rightarrow\quad \nor[t^q_i]\geq 2$, and
\item[(b)] $n\geq N\quad\Rightarrow\quad\nor[t^n_{A_n}]\geq 2$, where $A_n$ is as
in the assumption $(\delta)$.
\end{enumerate}
Suppose $M>N$ and $h:\{\rho\in\fs:N<\lh(\rho)\leq M\}\longrightarrow
2$. Build a condition $p=(w^p,t^p_0,t^p_1,\ldots)$ such that
\begin{enumerate}
\item[(c)] $w^q\vartriangleleft w^p$, $w^p(m)\in\pos(t^q_j)$ if
$m=m^{t^q_j}_\dn<\lh(w^p)$, and $t^p_i\in\Sigma(t^q_j)$ whenever
$m^{t^q_j}_\dn=m^{t^p_i}_\dn$,  
\item[(d)] if $N<\lh(\rho)\leq M$, $h(\rho)=1$, then for some $m=m(\rho)\in 
K_{\lh(\rho)+1}$ we have $m<\lh(w^p)$ and $w^p(m)\in \{g_{\rho\conc\langle
0\rangle}(m),g_{\rho\conc\langle 1\rangle}(m)\}\cap\pos(t^q_i)$, where
$m^{t^q_i}_\dn=m$, 
\item[(e)] if $m^{t^q_i}_\dn\geq\lh(w^p)$ then $\nor[t^q_i]\geq 2^{M+3}$,
\item[(f)] if $m<\lh(w^p)$ is not any of the $m(\rho)$'s from clause (d)
above (for $h(\rho)=1$), $m\in K_\ell$, $\ell>N$, then $w^p(m)\in\pos(t^q_i)
\setminus\{g_\rho(m):\rho\in 2^{\textstyle\ell}\}$, where $m^{t^q_i}_\dn=m$
(remember: $(K,\Sigma)$ is linked and (a)+(b)),
\item[(g)] if $m=\lh(w^p)+i\in K_\ell$, $N<\ell\leq M+1$ and $m^{t^q_j}_\dn
=m$, then $t^p_i\in\Sigma(t^q_j)$ is such that 
\[\nor[t^p_i]\geq\nor[t^q_j]-2^{M+1}\quad\mbox{ and }\quad \pos(t^p_i)\cap
\{g_\rho(m):\rho\in 2^{\textstyle\ell}\}=\emptyset\]
(remember assumption $(\gamma)$ and clause (e)),
\item[(h)] if $m=\lh(w^p)+i\notin\bigcup\limits_{\ell=N+1}^M K_\ell$,
$m=m^{t^q_j}_\dn$, then $t^p_i=t^q_j$. 
\end{enumerate}
It should be clear that we can build $p\in\bQ^{\bar{K},\bar{g}}_\infty(K,
\Sigma)$ satisfying the demands (c)--(h) and that then $q\leq p$ and
$p\forces h\subseteq\dot{\tau}_0$.

Now we easily conclude that $\dot{\tau}_0$ is Cohen over $\bV$; the rest
should be clear too.
\end{proof}

For $\rho\in\can$ and $k<\omega$ let $p_{\rho,k}\in\bQ^{\bar{K},
\bar{g}}_\infty(K,\Sigma)$ be such that $w^{p_{\rho,k}}=\langle\rangle$, if
$i\in K_\ell$, $i>k$, then $\pos(t^{p_{\rho,k}}_i)=\bH(i)\setminus \{g_{\rho
\rest\ell}(i)\}$, and if either $i\leq k$ or $i\notin\bigcup\limits_{\ell\in
\omega} K_\ell$ then $\pos(t^{p_{\rho,k}}_i)=\bH(i)$. Note that the function
$g^*:\can\times\omega\longrightarrow\bQ^{\bar{K},\bar{g}}_\infty(K,\Sigma):
(\rho,k)\mapsto p_{\rho,k}$ is Borel.

\begin{claim}
\label{cl21}
\begin{enumerate}
\item For each $\rho\in\can$ the set $\{p_{\rho,k}:k<\omega\}$ is predense
in $\bQ^{\bar{K},\bar{g}}_\infty(K,\Sigma)$.
\item If $\rho_0,\rho_1\in\can$, $\rho_0\rest m=\rho_1\rest m=\sigma$,
$\rho_0(m)=0$, $\rho_1(m)=1$ and $k<m$, then 
\[(\forall q\in \bQ^{\bar{K},\bar{g}}_\infty(K,\Sigma))(q\geq p_{\rho_0,k}\
\&\ q\geq p_{\rho_1,k}\quad\Rightarrow\quad q\forces\dot{\tau}_0(\sigma)
=0).\] 
\end{enumerate}
\end{claim}

\begin{proof}[Proof of the claim] 
(1)\quad Straightforward.\\
(2)\quad Note that if $q\geq p_{\rho_0,k}$, $q\geq p_{\rho_1,k}$ then for
each $i$, $n=m^{t^q_i}_\dn\in K_{m+1}$ implies ($n\geq m+1>k$ and)
$g_{\rho_0\rest (m+1)}(n), g_{\rho_1\rest (m+1)}(n)\notin\pos(t^q_i)$ (and,
of course, $\sigma\conc\langle 0\rangle=\rho_0\rest (m+1)$, $\sigma\conc
\langle 1\rangle=\rho_1\rest (m+1)$). Also, if $n<\lh(w^q)$, $n\in K_{m+1}$, 
then ($n>k$ and) $w^q(n)\notin\{g_{\rho_0\rest (m+1)}(n),g_{\rho_1\rest
(m+1)}\}$.  
\end{proof}

\noindent Now one easily shows that $(\bQ^{\bar{K},\bar{g}}_\infty(K,
\Sigma), \dot{\tau}_0)$ is very explicitly sour to $(\bQ^{\bar{K},
\bar{g}}_\infty(K,\Sigma),\dot{\tau}_1)$. 
\medskip

\noindent 2), 3)\quad Similarly.
\medskip

\noindent 4)\quad Here we think about the Cohen forcing $\bC$ as the set of
finite partial functions from $\fs\times\omega$ to 2 ordered by the
extension; $\bH_\bC$, $\dot{W}_\bC$ are interpreted suitably.

For each $n<\omega$ pick $a_n\in\bH(n)$. Let $\dot{W}$ be the name for
$\bQ^{\bar{U}}_\infty(K,\Sigma)$--generic real, and take standard
$\bQ^{\bar{U}}_\infty(K,\Sigma)$--names $\dot{\tau}_0,\dot{\tau}_1$ for
functions from $\fs\times\omega$ from 2 such that 
\begin{itemize}
\item $\dot{\tau}_0(\sigma,k)=0$ if and only if $(\forall n\in U_{\sigma
\conc\langle 0\rangle,k}\cup U_{\sigma\conc\langle 1\rangle,k})(\dot{W}(n)
\neq a_n)$, 
\item $\dot{\tau}_1(\sigma,k)=1-\dot{\tau}_0(\sigma,k)$.
\end{itemize}

\begin{claim}
\label{cl22}
\[\forces_{\bQ^{\bar{U}}_\infty(K,\Sigma)}\mbox{`` $\dot{\tau}_0,
\dot{\tau}_1$ are $\bC$--generic over $\bV$ '';}\]
moreover $\lbv\dot{\tau}_\ell\in B\rbv_{\bQ^{\bar{U}}_\infty(K,\Sigma)}\neq
{\bf 0}$ for any non-meager Borel set $B\subseteq 2^{2^{{<}\omega}}$. Hence
$(\bC,\dot{W}_\bC), (\bQ^{\bar{U}}_\infty(K,\Sigma),\dot{\tau}_0),
(\bQ^{\bar{U}}_\infty(K,\Sigma),\dot{\tau}_1)$ are as in \ref{soucon}. 
\end{claim}

\begin{proof}[Proof of the claim] 
Quite similar to \ref{cl20}.
\end{proof}

Now, for $\rho\in\can$ and $n,k<\omega$ let $p^n_{\rho,k}\in
\bQ^{\bar{U}}_\infty(K,\Sigma)$ be such that $w^{p^n_{\rho,k}}=\langle
\rangle$, and if $i\in U_{\rho\rest m,k}$, $m<\omega$, $n\leq i$, then
$\pos(t^{p^n_{\rho,k}}_i)=\bH(i)\setminus\{a_i\}$, and $\pos(t^{p^n_\rho,
k}_i)=\bH(i)$ is all other cases.

\begin{claim}
\label{cl23}
\begin{enumerate}
\item For each $\rho\in\can$ and $k\in\omega$, the set $\{p^n_{\rho,k}:n<
\omega\}$ is predense in $\bQ^{\bar{U}}_\infty(K,\Sigma)$.
\item If $k<\omega$, $\rho_0,\rho_1\in\can$, $\rho_0\rest m=\rho_1\rest m=
\sigma$, $\rho_0(m)\neq\rho_1(m)$ and $n<m$, then 
\[(\forall q\in \bQ^{\bar{U}}_\infty(K,\Sigma))(q\geq p^n_{\rho_0,k}\ \&\
q\geq p^n_{\rho_1,k}\quad\Rightarrow\quad q\forces\dot{\tau}_0(\sigma,k)=
0).\]
\end{enumerate}
\end{claim}

\begin{proof}[Proof of the claim] 
Straightforward.
\end{proof}

Now we easily finish.
\medskip

\noindent 5)\quad Similarly.
\end{proof}

\begin{definition}
\label{sourSystem}
Let $(K,\Sigma)$ be a local creating pair for $\bH$ and let $\cF\subseteq 
(\omega\setminus 2)^{\textstyle\omega}$. {\em A sourness system for $(K,
\Sigma,\cF)$} is a pair $(\bar{g},\bar{\ell})$ such that
\begin{enumerate}
\item[$(\alpha)$] $\bar{\ell}=\langle \ell_k:k<\omega\rangle\subseteq
\omega$ is increasing,
\item[$(\beta)$] $\bar{g}=\langle g_\rho:\rho\in\fs\rangle$, if $\rho\in
2^{\textstyle k}$ then $g_\rho\in\prod\limits_{i<\ell_k}\cP(\bH(i))$, and 
$\rho\vartriangleleft\rho'\ \Rightarrow\ g_\rho\vartriangleleft g_{\rho'}$,
\item[$(\gamma)$] for each $n\in [\ell_k,\ell_{k+1})$, $k<\omega$, the sets
$\{g_\rho(n):\rho\in 2^{\textstyle k+1}\}$ are pairwise disjoint and
non-empty, and $\bH(n)\setminus\bigcup\{g_\rho:\rho\in 2^{\textstyle k+1}\}
\neq\emptyset$,   
\item[$(\delta)$] if $f\in\cF$, $\langle t_0,t_1,\ldots\rangle\in
\PC(K,\Sigma)$ (see \ref{candidates}), $m^{t_0}_\dn=\ell_{k_0}$, $\nor[t_n]
\geq f(n+\ell_{k_0})$ for $n<\omega$ then: 
\begin{enumerate}
\item[(i)] for some $N<\omega$, for each $k\geq k_0$, we have
\[|\{\rho\in 2^{\textstyle k+1}:|\{n\in [\ell_k,\ell_{k+1}): g_\rho(n)\cap
\pos(t_{n-\ell_{k_0}})\neq\emptyset\}|<2^{k+1}\}|\leq N,\] 
\item[(ii)] for some $k_1\geq k_0$ we have
\[(\forall k\geq k_1)(\forall n\in [\ell_k,\ell_{k+1}))(\pos(t_{n-
\ell_{k_0}})\setminus\bigcup\{g_\rho(n):\rho\in 2^{\textstyle k+1}\}\neq
\emptyset)\]
\end{enumerate}
\end{enumerate}
\end{definition}

\begin{theorem}
\label{souF}
Suppose that $h:\omega\times\omega\longrightarrow\omega$ is regressive, and 
$\cF\subseteq (\omega\setminus 2)^{\textstyle\omega}$ is a countable
$h$--closed and $\geq^*$--directed family. Let $(K,\Sigma)$ be a local,
$h$--linked and complete creating pair for $\bH$ such that
\begin{enumerate}
\item[(a)] for some $f^*\in\cF$ we have that $(\forall n<\omega)(\nor
[ t^n_{\bH(n)}]\geq f^*(n))$ (see \ref{complete}),
\item[(b)] there is a sourness system $(\bar{g},\bar{\ell})$ for
$(K,\Sigma,\cF)$; let $\bar{g}=\langle g_\rho:\rho\in\fs\rangle$,
$\bar{\ell}=\langle\ell_k: k<\omega\rangle$,
\item[(c)] if $A\subseteq\bH(n)$, $\ell_k\leq n<\ell_{k+1}$, $\rho\in
2^{\textstyle k+1}$ and $A\setminus g_\rho(n)\neq\emptyset$, then
\[\nor[t^n_{A\setminus g_\rho(n)}]\geq h(n,\nor[t^n_A]).\]
\end{enumerate}
Then the forcing notion $\bQ^*_\cF(K,\Sigma)$ is very explicitly sour over
Cohen. 
\end{theorem}

\begin{proof}
Here we interpret the Cohen forcing notion $\bC$ in the standard way, i.e.,
it is $(\fs,\vartriangleleft)$ (and $\bH_\bC$, $\dot{W}_\bC$ are natural). 

For $\rho\in\can$ let $g_\rho=\bigcup\limits_{i\in\omega}g_{\rho\rest i}\in
\prod\limits_{i\in\omega}\cP(\bH(i))$ (remember \ref{sourSystem}$(\beta)$). 
Let $\dot{W}$ be the canonical name for the $\bQ^*_\cF(K,\Sigma)$--generic
real and let $\dot{\tau}_0,\dot{\tau}_1$ be standard $\bQ^*_\cF(K,
\Sigma)$--names for reals in $\can$ such that 
\begin{itemize}
\item $\dot{\tau}_0(k)=0$\quad if and only if\quad there are $\rho_0,\rho_1
\in\can$ such that $\rho_0\rest k=\rho_1\rest k$, $\rho_0(k)\neq\rho_1(k)$
and $(\forall n\geq\ell_k)(\dot{W}(n)\notin g_{\rho_0}(n)\cup
g_{\rho_1}(n))$;  
\item $\dot{\tau}_1(n)=1-\dot{\tau}_0(n)$.
\end{itemize}

\begin{claim}
\label{cl24}
$(\bC,\dot{W}_\bC),(\bQ^*_\cF(K,\Sigma),\dot{\tau}_0),(\bQ^*_\cF(K,\Sigma),
\dot{\tau}_1)$ are as in \ref{soucon}.
\end{claim}

\begin{proof}[Proof of the claim] 
We will show that $\dot{\tau}_0$ is (a name for) a Cohen real over $\bV$;
then the rest should be clear.

So suppose $p=(w^p,t^p_0,t^p_1,\ldots)\in\bQ^*_\cF(K,\Sigma)$ and we may
assume that $\lh(w^p)=\ell_{k_0}$, $k_0<\omega$. For $k\geq k_0$ let 
\[\Upsilon^p_k=\Upsilon_k\stackrel{\rm def}{=} \{\rho\in 2^{\textstyle k+1}:
|\{n\in [\ell_k,\ell_{k+1}):g_\rho(n)\cap\pos(t^p_{n-\ell_{k_0}})\neq
\emptyset\}|< 2^{k+1}\},\] 
and let $N$ be such that $(\forall k\geq k_0)(|\Upsilon_k|<N)$ (remember
\ref{sourSystem}($\delta$(i)). Let 
\[u_k=\{i\leq k:(\exists\rho_0,\rho_1\in\Upsilon_k)(\rho_0\rest i=\rho_1
\rest i\ \&\ \rho_0(i)\neq\rho_1(i))\}.\]
Clearly, for all $k\geq k_0$, $|u_k|\leq N$, and thus we may choose an
infinite set $A\subseteq\omega$ such that $\{u_k:k\in A\}$ forms a
$\Delta$--system with the heart, say, $u^*$. Take $k_1>\max(u^*)+k_0+N$ such
that
\begin{enumerate}
\item[$(*_1)$] $(\forall k\geq k_1)(\forall n\in [\ell_k,\ell_{k+1}))
(\pos(t^p_{n-\ell_{k_0}})\setminus\bigcup\{g_\rho(n):\rho\in 2^{\textstyle
k+1}\}\neq\emptyset)$
\end{enumerate}
(remember \ref{sourSystem}($\delta$(ii)) and suppose that $v\subseteq
[k_1, k_2)$, $k_2>k_1$. We are going to build a condition $q\geq p$ such
that 
\[q\forces(\forall j\in [k_1,k_2))(\dot{\tau}_0(j)=0\ \Leftrightarrow\ j\in 
v).\] 
To this end, pick $i>k_2$ such that
\begin{enumerate}
\item[$(*_2)$] $u_i\cap [k_1,k_2)=\emptyset$, and
\item[$(*_3)$] for some $f\in\cF$ we have
\[(\forall n\geq \ell_{i+1})(f(n)<h^{(k_2)}(n,\nor[t^p_{n-\ell_{k_0}}]),\]
where $h^{(k+1)}(n,m)=h(n,h^{(k)}(n,m))$.
\end{enumerate}
(Possible by the choice of $k_1$ and the assumption that $\cF$ is
$h$--closed.) Since $|\Upsilon_i|<N<2^{k_1}$, we may find $\rho^*\in\can$
such that $\rho^*\rest k_1\notin \{\sigma\rest
k_1:\sigma\in\Upsilon_i\}$. For $k\in v$ fix $\rho_k\in\can$ such that
$\rho_k\rest k=\rho^*\rest k$, $\rho_k(k)=1-\rho^*(k)$. 

Let $\langle \sigma_j:j<2^{i+1}-|v|-1\rangle$ enumerate $2^{\textstyle
i+1}\setminus\{\rho^*\rest (i+1),\rho_k\rest (i+1): k\in v\}$. By induction
on $j<2^{i+1}-|v|-1$ define $n^*_j\in [\ell_i,\ell_{i+1})\cup \{*\}$ as
follows:
\medskip

if there is $n\in [\ell_i,\ell_{i+1})\setminus \{n^*_{j'}:j'<j\}$ such that
$g_{\sigma_j}(n)\cap\pos(t^p_{n-\ell_{k_0}})\neq\emptyset$, 

then $n^*_j$ is the first such number,

otherwise $n^*_j$ is $*$.
\medskip

\noindent Now we choose $w^q,t^q_0,t^q_1,\ldots$ so that:
\begin{enumerate}
\item[(i)] $\lh(w^q)=\ell_{i+1}$, $w^p\triangleleft w^q$;
\item[(ii)] if $n\in [\ell_i,\ell_{i+1})$, $n=n^*_j$, $j<2^{i+1}-|v|-1$,
then $w^q(n)\in g_{\sigma_j}(n)\cap\pos(t^p_{n-\ell_{k_0}})$;
\item[(iii)] if $\ell_{k_0}\leq n<\ell_{k_1}$, then $w^q(n)\in\pos(
t^p_{n-\ell_{k_0}})$;
\item[(iv)] if $n\in [\ell_{k_1},\ell_i)\cup ([\ell_i,\ell_{i+1})
\setminus\{n^*_j: j<2^{i+1}-|v|-1\})$, then 
\[w^q(n)\in\pos(t^p_{n-\ell_{k_0}})\setminus\bigcup\{g_\rho(n):\rho\in
2^{\textstyle i+1}\};\]  
\item[(v)] $t^q_n\in\Sigma(t^p_{n+\ell_{i+1}-\ell_{k_0}})$ is such that
$\nor[t^q_n]\geq f(n+\ell_{i+1})$ and
\[\big(g_{\rho^*}(n+\ell_{i+1})\cup\bigcup_{k\in v}g_{\rho_k}(n+\ell_{i+1})
\big)\cap\pos(t^q_n)=\emptyset\]  
(where $f$ is given by $(*_3)$). 
\end{enumerate}
[Why is the choice possible? Demands (i)--(iii) are easy; (iv) can be
satisfied by $(*_1)$, remember \ref{sourSystem}$(\beta)$; (v) is possible
by the assumption (c) of the theorem and $(*_3)$.] One easily checks that
the demands (i)--(v) imply $q=(w^q,t^q_0,t^q_1,\ldots)$ is a condition in
$\bQ^*_\cF(K,\Sigma)$ stronger than $p$. Also, by (ii)+(iv)+(v), 
\[q\forces (\forall k\in v)(\forall n\geq \ell_k)(\dot{W}(n)\notin
g_{\rho^*}(n)\cup g_{\rho_k}(n))\]
(remember \ref{sourSystem}($\gamma$); thus $g_{\sigma_j}(n)\cap g_{\rho^*}
(n)=g_{\sigma_j}(n)\cap g_{\rho_k}(n)=\emptyset$ in clause (ii)). Hence
$q\forces (\forall k\in v)(\dot{\tau}_0(k)=0)$. Now we argue that $q\forces
(\forall k\in [k_1,k_2)\setminus v)(\dot{\tau}_0(k)=1)$. If not, then for
some $k\in [k_1,k_2)\setminus v$ we find $\rho_0^+,\rho_1^+\in 2^{\textstyle
i+1}$ such that 
\[\rho_0^+\rest k=\rho_1^+\rest k,\quad \rho_0^+(k)\neq\rho_1^+(k),\quad
\mbox{ and }\quad(\forall n\in [\ell_k,\ell_{i+1}))(w^q(n)\notin
g_{\rho_0^+}(n)\cup g_{\rho_1^+}(n)).\]
Necessarily, $\{\rho^+_0,\rho^+_1\}\nsubseteq\{\rho^*\rest (i+1),\rho_k\rest
(i+1): k\in v\}$. Moreover, if $\rho_\ell^+\in\{\sigma_j:j<2^{i+1}-|v|-1\}$,
then $\rho^+_\ell\in\Upsilon_i$ (as if $\rho^+_\ell=\sigma_j\notin
\Upsilon_i$ then $n^*_j\in [\ell_i,\ell_{i+1})$ and $w^q(n^*_j)\in g_{
\rho^+_\ell}(n^*_j)$). Since $\rho^*\rest k_1\notin\{\sigma\rest k_1:\sigma
\in\Upsilon_i\}$, we may conclude that $\rho^+_0,\rho^+_1\in \{\sigma_j:
j<2^{i+1}-|v|-1\}$ (and thus both are in $\Upsilon_i$). However, then we get
$k\in u_i$, what contradicts $(*_2)$.
\end{proof}
For $\rho\in\can$ and $n\in\omega$ let $p_{\rho,n}\in\bQ^*_\cF(K,\Sigma)$ be
such that $w^{p_{\rho,n}}=\langle\rangle$, $t^{p_{\rho,n}}_k=t^k_{\bH(k)}$
for $k<n$, and $t^{p_{\rho,n}}_k=t^k_{\bH(k)\setminus g_\rho(k)}$ for $k\geq
n$. Let  
\[g_0,g_1:\can\times\omega\longrightarrow\bQ^*_\cF(K,\Sigma):(\rho,n)\mapsto
p_{\rho,n}.\]
It is straightforward to check that $g_0,g_1$ witness $(\bQ^*_\cF(K,\Sigma),
\dot{\tau}_0)$ is very explicitly sour to $(\bQ^*_\cF(K,\Sigma),
\dot{\tau}_1)$ (note that if $k,k'<m$, $\rho_0,\rho_1\in\can$, $\rho_0\rest
m=\rho_1\rest m$, $\rho_0(m)\neq\rho_1(m)$ and $q\geq g_0(\rho_0,k),g_1(
\rho_1,k')$, then $q\forces$`` $\dot{\tau}_0(m)=0\ \&\ \dot{\tau}_1(m)=1$
''). 
\end{proof}

\begin{remark}
In \ref{souF}, in the assumptions on the family $\cF$, instead of demanding
that ``$\cF$ is countable $h$--closed and $\geq^*$--directed'', we may
require that ``$\cF$ is $h$--closed and either countable or
$\geq^*$--directed'', and then conclude that $\bQ^*_\cF(K,\Sigma)$ is
explicitly sour. The ``countable and $\geq^*$--directed'' assumption is
needed only to be in the context of \ref{soucon} (i.e., to make sure that
$\bQ^*_\cF(K,\Sigma)$ is Souslin). 
\end{remark}
 
\subsection{Conclusions}

\begin{conclusion}
\label{topsweEX}
\begin{enumerate}
\item Let $\bH:\omega\longrightarrow\cHa$. The following forcing notions are 
topologically sweet: 
\begin{enumerate}
\item[(a)] $\bQ^*_\infty(K_\bH,\Sigma_\bH,\Sigma^\bot_\bH)$ of \ref{bas628}
and $\bQ^*_\infty(K_{\ref{bas628bis}},\Sigma_{\ref{bas628bis}},
\Sigma^\bot_{\ref{bas628bis}})$,
\item[(b)] $\bQ^*_f(K_{\ref{loc628}},\Sigma_{\ref{loc628}})$ for $f$ as in
\ref{conclOne}(1), and $\bQ^*_\cF(K_{\ref{loc628}},\Sigma_{\ref{loc628}})$
for $\cF$ as defined in \ref{conclOne}(2) (the ``e.g.'' part).
\end{enumerate}
\end{enumerate}
\end{conclusion}

\begin{conclusion}
The forcing notions $\bQp$ for the universal parameters $\gp$ defined in 
\ref{UM}, \ref{PPex} and \ref{ucmz} are iterably sweet.
\end{conclusion}

\begin{conclusion}
\label{strange}
Let $\bH:\omega\longrightarrow\omega$, $\bH(i)\geq 2^{2^{i+1}}$, and let $h$ be
as in \ref{EDRF}.  
\begin{enumerate}
\item Let $f^0_k(n)=\max\{2,\bH(n)-2^{kn}\}$ and $\cF_0=\{f^0_k:k<\omega\}$. 
Then the forcing notion $\bQ^*_{\cF_0}(K_{\ref{EDRF}},\Sigma_{\ref{EDRF}})$
(constructed for $\bH$ as in \ref{EDRF}) is topologically sweet.
\item Let $f^1_k(n)=\max\{2,\bH(n)-2^k\}$, $\cF_1=\{f^1_k:k<\omega\}$. Then
the forcing notion $\bQ^*_{\cF_1}(K_{\ref{EDRF}},\Sigma_{\ref{EDRF}})$ is
very explicitly sour over Cohen.\\
Similarly if $f^2_k(n)=\max\{2,\bH(n)-(k2^n)\}$, $\cF_2=\{f^2_k:k<\omega
\}$.
\end{enumerate}
\end{conclusion}

\begin{proof}
(2)\quad We are going to apply \ref{souF}. First note that even though
$(K_{\ref{EDRF}},\Sigma_{\ref{EDRF}})$ as defined in \ref{EDRF} is not
complete we can easily make it so, or restrict our attention to the forcing
notion below some condition (the problems with completeness come from the
technical requirement in the definition of $t\in K_{\ref{EDRF}}$ that
$E_t\neq \emptyset$). 

We are going to build a sourness system $(\bar{g},\bar{\ell})$ for
$(K_{\ref{EDRF}},\Sigma_{\ref{EDRF}},\cF_1)$ such that the demand
\ref{souF}(c) holds. 

Let $\ell_0=0$, $\ell_{k+1}=\ell_k+2^{2^k}$. For $\rho\in\fs$ pick $g_\rho$
such that
\begin{enumerate}
\item[$(\oplus_1)$] if $\rho\in 2^{\textstyle k}$, then $g_\rho\in
\prod\limits_{i<\ell_k}\bH(i)$, and $\rho\vartriangleleft\rho'\ \Rightarrow\
g_\rho\vartriangleleft g_{\rho'}$, 
\item[($\oplus_2$)] if $n\in [\ell_k,\ell_{k+1})$, $k<\omega$, then there
are no repetitions in the sequence $\langle g_\rho(n):\rho\in 2^{\textstyle
k+1}\rangle$. 
\end{enumerate}
We claim that, letting $\bar{\ell}=\langle\ell_k:k<\omega\rangle$ and
$\bar{g}=\langle g_\rho:\rho\in\fs\rangle$, $(\bar{g},\bar{\ell})$ is as
required (we identify $\bH(i)$ with $[\bH(i)]^{\textstyle 1}$, of
course). Clauses \ref{sourSystem}($\alpha$)--($\gamma$) are clear.

Suppose that $\langle t_0,t_1,\ldots\rangle\in\PC(K_{\ref{EDRF}},
\Sigma_{\ref{EDRF}})$, $m^{t_n}=n$, $\dis[t_n]=(n,E_n)$, $\nor[t_n]\geq
f^1_N(n)$. Then, for large enough $n$, $|E_n|\leq 2^N$, so let
$M=\max\{|E_n|: n<\omega\}$. Assume that $\rho_i\in 2^{\textstyle k+1}$ (for
$i\leq M$) are pairwise distinct. By $(\oplus_2)$, for each $n\in
[\ell_k,\ell_{k+1})$ there is $i\leq M$ such that $g_{\rho_i}(n)\notin
E_n$. Hence for some $i\leq M$    
\[|\{n\in [\ell_k,\ell_{k+1}): g_{\rho_i}(n)\notin E_n\}|\geq
\frac{\ell_{k+1}-\ell_k}{M+1}=\frac{2^{2^k}}{M+1}.\]
Hence we easily conclude that \ref{sourSystem}($\delta$(i)) holds. The
demands \ref{sourSystem}($\delta$(ii)) and \ref{souF}(c) are even easier. 
\medskip 

For $\cF_2$ we proceed similarly, but we choose $g_\rho$ so that
$g_{\rho}(n)\in [\bH(n)]^{\textstyle 2^n}$.  
\end{proof}

\section{Epilogue}
A general problem that we have in mind in this paper is classifying ``nice''
ccc forcing notions, in particular finding dividing lines in this family, or
at least natural properties. We should explain what we mean. A forcing
notion is ``nice'' if it has a quite absolute definition, so Borel is
natural, but Souslin is more central (see \ref{verysouslin} and also
\cite{Sh:630}), but we may be happy with just ``one of the form presented in
this paper''. A dividing line is a property of such definitions, so that
both it and its negation is meaningful (that is we can prove theorems from
both). Thus a dividing line may serve as a division to cases in solving
problems. (On parallel in Model Theory see \cite{Sh:c} and \cite{Sh:702}.)

The first (possible) dividing line we considered here is determined by
``being $\omega$--nw--nep'' (or just ``being very Souslin ccc''). In some
sense, one can consider very Souslin forcing notion as those which are
really close to random and Cohen. We have examples of very Borel ccc forcing
notions (see \ref{conclOne}(2), \ref{extra}, \ref{conclThree}(3)), and
forcing notions which are not $\omega$--nw--nep (see \ref{notCon1},
\ref{notCon2}, \ref{conclOne}(1), \ref{conclThree}(1,2)). The argument for
``not being $\omega$--nw--nep'' was in all cases the same: adding a
dominating real. So we arrive to the following question.

\begin{problem}
Suppose $\bP$ is a Borel ccc forcing notion which is not equivalent to a 
$\omega$--nw--nep forcing. Does $\bP$ add a dominating real?
\end{problem}

If one looks at \ref{allsimple}(3) and \ref{extra}, then the following
(perhaps less central but still intriguing) question related to
$\omega$--nw--nep forcing notions arises. 

\begin{problem}
Assume $\bH,K,\Sigma,\cF$ and $h$ are as in \ref{allsimple}(1c) or as in
\ref{allsimple}(2b). Is the forcing notion $\bQ^*_{\cF}(K,\Sigma)$ (or
$\bQ^\tree_{\cF}(K,\Sigma)$, respectively) very Borel ccc? (Of course, we
are interested in non-finitary $(K,\Sigma)$.) 
\end{problem}

The second dividing line originates in \cite{Sh:176} and studies of the
Baire property (and measurability) of projective sets. To get a model in
which all projective sets have the Baire property, \cite{Sh:176} uses
sweetness while \cite{St85}  applies topological sweetness. However, the use
of the two variants of sweetness might be slightly confusing. What we really
need for this type of construction are two properties, say, (a)--sweetness
and (b)--sweetness such that 
\begin{enumerate}
\item[(i)]  if $\bP$ is (a)--sweet and $\dbQ$ is a $\bP$--name for a
(b)--sweet forcing notion, then $\bP*\dbQ$ is (a)--sweet,
\item[(ii)] if $\bP_0,\bP_1$ are (a)--sweet then for sufficiently many
forcing notions $\bQ$ and their two complete embeddings $f_\ell:\bQ
\longrightarrow{\rm RO}(\bP_\ell)$, the  amalgamation $\bP_0\times_{f_0, 
f_1}\bP_1$ is (a)--sweet,  
\item[(iii)] the Universal Meager forcing notion is (b)--sweet,
\item[(iv)]  (a)--sweetness implies the ccc.
\end{enumerate}
(Note that in (ii) we do not require that all amalgamations are (a)--sweet,
we just need to cover the amalgamations needed to ensure suitable
homogeneity of the Boolean algebra we construct; see \cite{JuRo}, also
\ref{Solovay} here.) To some extend this approach was materialized in
\ref{getiter}: the topological sweetness may serve as (a)--sweetness and
iterable sweetness is a good candidate for (b)--sweetness. It should be
remarked here, that it is quite surprising that compositions of
(topologically) sweet forcing notions with the Universal Meager (or the
Hechler forcing notion) are topologically sweet because {\em the second
iterand is sweet}. (The respective proofs in \cite{Sh:176}, \cite{St85} were
somewhat less general.) Still, it is very reasonable to ask
\begin{problem}
Can \ref{getiter} be improved by weakening the demands on $\dbQ$? Can you
find (a)--sweetness and (b)--sweetness satisfying (i)--(iv) and weaker then
topological sweetness and iterable sweetness, respectively?
\end{problem}

\noindent The sweet/sour division is sometimes very surprising --- compare
\ref{strange}(1) and \ref{strange}(2). The forcing notions
$\bQ^*_{\cF_0}(K_{\ref{EDRF}},\Sigma_{\ref{EDRF}})$ and $\bQ^*_{\cF_2}(
K_{\ref{EDRF}},\Sigma_{\ref{EDRF}})$ (of \ref{strange}) look very similar
and one could expect that both are like the Cohen forcing. However, the
first is topologically sweet (so not so far from Cohen) while the other is
very explicitly sour (so one could even say that worse than random).  

Topological sweetness occurs to be not so seldom (see \ref{topsweEX}),
however it does not imply that we could make real use of these forcing
notions in constructions like \cite{Sh:F380}. These forcings seem to be
quite far from the iterable sweetness, so we conjecture that that the
following has an affirmative answer. 
\begin{problem}
Let $\bQ$ be one of the forcing notion covered by \ref{TSforCREA}(1,2)
and \ref{Ftopsweet}. Is there $k<\omega$ such that the ($k$--step) iteration
$\bQ^{(k)}$ is sour over Cohen? Over $\bQ$?   
\end{problem}

Proposition \ref{TSforCREA}(3) gives sweet forcing notions, so they could be
of some use in constructions like that in \cite{Sh:F380}. However, do we
really need to force additionally with these forcings? In other particular:

\begin{problem}
Is there a universality parameter $\gp$ satisfying the requirements of
\ref{TSforCREA}(3) such that no finite iteration of the Universal Meager
forcing notion adds a $\bQp$--generic real? Does the Universal Meager
forcing add generic reals for the forcing of \ref{ucmz}? Of \ref{PPex}?
\end{problem}

An intriguing thing is that in the cases we proved sourness over the Cohen,
we actually got that the considered forcing notion is very explicitly sour
over Cohen (so in particular the amalgamation collapses $\con$). 

\begin{problem}
Let $\bQ$ be a Souslin ccc (or just nep ccc) forcing notion, and
$(\bP,\dot{W})$ is as in \ref{cont}. 
\begin{enumerate}
\item Assume that $\bQ$ is sour over $(\bP,\dot{W})$. Is it explicitly sour
over $(\bP,\dot{W})$?  Very explicitly?
\item Suppose that $\bQ$ is not topologically sweet and adds a Cohen
real. Is it (weakly) sour over Cohen?
\end{enumerate}
\end{problem}

Since sweet/sour division is related to the Baire property of projective
sets, let us finish with the following general problem.

\begin{problem}
\begin{enumerate}
\item Let $(\bP,\dot{W})$ be as in \ref{cont}, $\bP$ be of the type studied
in this paper. For a cardinal $\kappa$, let ${\mathcal I}_{\bP,
\dot{W}}^\kappa$ be a ${<}\kappa$--complete ideal generated by $\Ipw$.\\
What is the consistency strength of the statement ``every projective subset
of $\cX$ has the ${\mathcal I}_{\bP,\dot{W}}^\kappa$--Baire property'' ?\\
(We conjecture that it is always either {\rm ZFC} or ``{\rm ZFC} + there
exists an inaccessible cardinal, and we would like to characterize and/or
describe this dividing line.)
\item Similarly for (typically non-ccc) ideals ${\mathcal I}_\gp$ determined
by suitable universality parameters $\gp$ (see \ref{unideal}).
\end{enumerate}
\end{problem}


\begin{thebibliography}{10}
\makeatletter \renewcommand{\@biblabel}[1]{[#1]} \makeatother

\bibitem{BaJu95}
Tomek Bartoszy\'nski and Haim Judah.
\newblock {\em {Set Theory: On the Structure of the Real Line}}.
\newblock A K Peters, Wellesley, Massachusetts, 1995.

\bibitem{BrJd93}
J\"org Brendle and Haim Judah.
\newblock Perfect sets of random reals.
\newblock {\em Israel Journal of Mathematics}, 83:153--176, 1993.

\bibitem{GoJu92}
Martin Goldstern and Haim Judah.
\newblock {Iteration of Souslin Forcing, Projective Measurability and the Borel
  Conjecture}.
\newblock {\em Israel Journal of Mathematics}, 78:335--362, 1992.

\bibitem{Ha35}
P.~Hall.
\newblock {On representatives of subsets}.
\newblock {\em Journal of the London Mathematical Society}, 10:26--30, 1935.

\bibitem{JdSh:292}
Jaime Ihoda (Haim~Judah) and Saharon Shelah.
\newblock {Souslin forcing}.
\newblock {\em {The Journal of Symbolic Logic}}, 53:1188--1207, 1988.
{\bf [JdSh:292]}

\bibitem{JuRo}
Haim Judah and Andrzej Ros{\l}anowski.
\newblock {On Shelah's Amalgamation}.
\newblock In {\em Set Theory of the Reals}, volume~6 of {\em Israel
  Mathematical Conference Proceedings}, pages 385--414. 1992.

\bibitem{JuRo92}
Haim Judah and Andrzej Ros{\l}anowski.
\newblock {The ideals determined by Souslin forcing notions}.
\newblock unpublished notes, 1992; {\tt math.LO/9905144}\footnote{References
of the form {\tt math.XX/$\cdots$} refer to the {\tt arXiv.org/archive/math} 
Mathematics Archive}

\bibitem{JRSh:373}
Haim Judah, Andrzej Roslanowski, and Saharon Shelah.
\newblock {Examples for Souslin Forcing}.
\newblock {\em {Fundamenta Mathematicae}}, 144:23--42, 1994.
{\bf [JRSh:373]};   {\tt math.LO/9310224}

\bibitem{JdSh:446}
Haim Judah and Saharon Shelah.
\newblock {Baire Property and Axiom of Choice}.
\newblock {\em {Israel Journal of Mathematics}}, 84:435--450, 1993.
{\bf [JdSh:446]};   {\tt math.LO/9211213}

\bibitem{KeSo95}
Alexander~S. Kechris and S{\l}awomir Solecki.
\newblock {Approximation of analytic by Borel sets and definable countable
  chain conditions}.
\newblock {\em Israel Journal of Mathematics}, 89:343--356, 1995.

\bibitem{Ku84}
Kenneth Kunen.
\newblock {Random and Cohen Reals}.
\newblock In K.~Kunen and J.~E. Vaughan, editors, {\em {Handbook of
  Set--Theoretic Topology}}, pages 887--911. Elsevier Science Publishers B.V.,
  1984.

\bibitem{Mi81}
Arnold~W. Miller.
\newblock {Some properties of measure and category}.
\newblock {\em Transactions of the American Mathematical Society}, 266:93--114,
  1981.

\bibitem{Mo77}
John C.~II Morgan.
\newblock {Baire category from an abstract viewpoint}.
\newblock {\em Fundamenta Mathematicae}, 94:13--23, 1977.

\bibitem{Mo90}
John C.~II Morgan.
\newblock {\em {Point set theory}}, volume 131 of {\em {Monographs and
  Textbooks in Pure and Applied Mathematics}}.
\newblock Marcel Dekker, Inc, New York, 1990.

\bibitem{RoSh:470}
Andrzej Roslanowski and Saharon Shelah.
\newblock {Norms on possibilities I: forcing with trees and creatures}.
\newblock {\em {Memoirs of the American Mathematical Society}}, 141(671), 1999.
{\bf [RoSh:470]};   {\tt math.LO/9807172}

\bibitem{RoSh:628}
Andrzej Roslanowski and Saharon Shelah.
\newblock {Norms on possibilities II: more ccc ideals on
  $2^{\textstyle\omega}$}.
\newblock {\em {Journal of Applied Analysis}}, 3:103--127, 1997.
{\bf [RoSh:628]};   {\tt math.LO/9703222}

\bibitem{RoSh:670}
Andrzej Roslanowski and Saharon Shelah.
\newblock {Norms on possibilities III: strange subsets of the real line}.
\newblock {\em {in preparation}}.
{\bf [RoSh:670]}

\bibitem{RoSh:736}
Andrzej Roslanowski and Saharon Shelah.
\newblock {Measured creatures}.
\newblock {\em Annals of Mathematics}, submitted.
{\bf [RoSh:736]};   {\tt math.LO/0010070}

\bibitem{Sh:F380}
Andrzej Roslanowski and Saharon Shelah.
\newblock {446 revisited}.
\newblock {\em {in preparation}}.
{\bf [Sh:F380]}

\bibitem{Sh:176}
Saharon Shelah.
\newblock {Can you take Solovay's inaccessible away?}
\newblock {\em {Israel Journal of Mathematics}}, 48:1--47, 1984.
{\bf [Sh:176]}

\bibitem{Sh:480}
Saharon Shelah.
\newblock {How special are Cohen and random forcings i.e. Boolean algebras of
  the family of subsets of reals modulo meagre or null}.
\newblock {\em {Israel Journal of Mathematics}}, 88:159--174, 1994.
{\bf [Sh:480]};   {\tt math.LO/9303208}

\bibitem{Sh:666}
Saharon Shelah.
\newblock {On what I do not understand (and have something to say)}.
\newblock {\em {Fundamenta Mathematicae}}, 166:1--82, 2000.
{\bf [Sh:666]};   {\tt math.LO/9906113}

\bibitem{Sh:630}
Saharon Shelah.
\newblock {Non-elementary proper forcing notions}.
\newblock {\em {Journal of Applied Analysis}}, accepted.
{\bf [Sh:630]};   {\tt math.LO/9712283}

\bibitem{Sh:669}
Saharon Shelah.
\newblock {Non-elementary proper forcing notions II}.
\newblock {\em {In preparation}}.
{\bf [Sh:669]}

\bibitem{Sh:702}
Saharon Shelah.
\newblock {On what I do not understand (and have something to say), model
  theory}.
\newblock {\em Mathematica Japonica}, to appear.
{\bf [Sh:702]};   {\tt math.LO/9910158}

\bibitem{Sh:711}
Saharon Shelah.
\newblock {Note on $\omega$--nw--nep forcing notions}.
\newblock {\em Preprint}.
{\bf [Sh:711]}

\bibitem{Sh:c}
Saharon Shelah.
\newblock {\em {Classification theory and the number of nonisomorphic models}},
  volume~92 of {\em {Studies in Logic and the Foundations of Mathematics}}.
\newblock {North-Holland Publishing Co., Amsterdam, xxxiv+705 pp}, 1990.
{\bf [Sh:c]}

\bibitem{Sh:f}
Saharon Shelah.
\newblock {\em {Proper and improper forcing}}.
\newblock {Perspectives in Mathematical Logic}. {Springer}, 1998.
{\bf [Sh:f]}

\bibitem{So98}
S{\l}awomir Solecki.
\newblock {Analytic ideals and their applications}.
\newblock {\em Annals of Pure and Applied Logic}, 99:51--72, 1999.

\bibitem{So70}
Robert~M. Solovay.
\newblock {A model of set theory in which every set of reals is Lebesgue
  measurable}.
\newblock {\em Annals of Math.}, 92:1--56, 1970.

\bibitem{St85}
Jacques Stern.
\newblock Regularity properties of definable sets of reals.
\newblock {\em Annals of Pure and Applied Logic}, 29:289--324, 1985.
\end{thebibliography}

\end{document}